\parindent=0pt
\magnification=1200
\input amstex
\documentstyle{amsppt}
\voffset=-0.6in

\def\ddefi#1{\medskip {{\bf Definition #1.}

}}
\def\nota#1{\medskip {{\bf Notation #1.}

}}
\def\lem#1{\medskip {{\bf Lemma #1}

}}
\def\prop#1{\medskip {{\bf Proposition #1}

}}
\def\teo#1{\medskip {{\bf Theorem #1}

}}
\def\cor#1{\medskip {{\bf Corollary #1}

}}
\def\dim{\smallskip {{\bf Proof:} }}
\def\sqr#1#2{{\vcenter{\vbox{\hrule height .#2pt
                             \hbox{\vrule width .#2pt height#1pt\kern#1pt
                                   \vrule width .#2pt}
                             \hrule height .#2pt}}}}
\def\square{\mathchoice\sqr54\sqr54\sqr{4.1}3\sqr{3.5}3}
\def\finedim{{\hfill\hbox{\enspace{\bf $\square$}}} \smallskip}
\def\rem#1{\medskip {{\bf Remark #1.}

}}

\def\compo{\,{\scriptstyle\circ}\,}
\def\rapp#1#2{ { {\displaystyle {#1}} \over {\displaystyle {#2}} }}
\def\ugdef
{{\hskip 2 pt}\dot{=}_{{\hskip -5.5 pt}{\displaystyle .}}{\hskip 5
pt}}
\def\rm{\roman}

\def\vx#1{{\scriptstyle\scriptstyle\scriptstyle
\scriptstyle\circ}\!\!_{_{_{#1}}}\!\!}
\def\vn{{\scriptstyle\scriptstyle\scriptstyle
\scriptstyle\circ}\!\!\!_{_{_{n}}}\!\!}
\def\vl#1{{\scriptstyle\scriptstyle\scriptstyle
\scriptstyle\circ}\!\!\!\!\!_{_{_{#1}}}\!\!\!\!\!}
\def\vo#1{{\scriptstyle\scriptstyle\scriptstyle
\scriptstyle\circ}{#1}\!\!}

\def\edg{-\!\!\!-}
\def\dedg{\,{\scriptstyle{=\!\!\!=\!\!\!=\!\!=}}}
\def\tedg{\,{\scriptstyle{\equiv\!\!\!\equiv\!\!\!\equiv\!\!\equiv}}}
\def\qedg{\,\,^{_=}_{^=}\!\!^{_=}_{^=}\!\!^{_=}_{^=}}
\def\vedgv#1#2{\!^{^{\big|\!\!\!^{^{^{^{^{\vo{#1}\!\!\!^{#2}}}}}}}}\!\!}

\def\redg{}

\def\um{{\phantom{\vx{}-}}}
\def\dm{{\phantom{\vx{a}\edg}}}
\def\tm{{\phantom{\vx{a}\edg\vx{}-}}}
\def\qm{{\phantom{\vx{a}\edg\vx{a}\edg}}}
\def\cm{{\phantom{\vx{a}\edg\vx{a}\edg\vx{}-}}}
\def\sm{{\phantom{\vx{a}\edg\vx{a}\edg\vx{a}\edg}}}
\def\stm{{\phantom{\vx{a}\edg\vx{a}\edg\vx{a}\edg\vx{}-}}}
\def\om{{\phantom{\vx{a}\edg\vx{a}\edg\vx{a}\edg\vx{a}\edg}}}
\def\vm{{\phantom{\vx{a}\edg\vx{a}\edg\vx{a}\edg\vx{a}\edg\vx{}-}}}
\def\num{{\phantom{..}}}
\def\ndm{{\phantom{2}}}
\def\ntm{{\phantom{2..}}}
\def\nqm{{\phantom{22}}}
\def\ncm{{\phantom{22..}}}
\def\nsm{{\phantom{222}}}
\def\nstm{{\phantom{222..}}}

\def\nnm{{\phantom{2222..}}}
\def\ndcm{{\phantom{22222}}}

\def\mn{\Bbb N}
\def\mz{\Bbb Z}

\def\mr{\Bbb R}
\def\mc{\Bbb C}

\def\sl{\frak{sl}}
\def\g{\frak{g}}

\def\u{{\Cal{U}}}
\def\uq{{\Cal{U}_q}}
\def\a{\alpha}
\def\d{\delta}

\def\omc{\omega\check{}}

\def\Pc{P\check{}}
\def\autq{Aut(Q_{(\cdot|\cdot),\d})}
\def\cnt{b}
\def\dk{\tilde d}
\def\ug{\underline{\gamma}}

\def\T{{\Cal T}}
\def\E{{\Cal E}}
\def\Pp{{\Cal P}}

\def\tn{\tilde n}
\def\tW{\tilde W}
\def\tw{\tilde w}
\def\tB{\tilde{{\Cal{B}}}}
\def\tPhi{\tilde\Phi}
\def\te{\tilde E}
\def\tf{\tilde F}
\def\tnu{\tilde\nu}
\def\ts{\tilde S}
\def\tp{\tilde p}
\def\tc{\tilde C}
\def\tk{\tilde k}

\def\intr{0} 

\def\genset{1} 

\def\akma{\genset.{1}} 

\def\ntz{\genset.{2}} 
                      
\def\rlatt{\ntz.{A}} 
\def\sdfscpr{\ntz.{B}} 
\def\homz{\ntz.{C}} 
\def\wgrs{\ntz.{D}} 
\def\rsys{\ntz.{E}} 

\def\classdgr{\genset.{3}} 

\def\defstr{\genset.{4}}

\def\QNTALG{\defstr.{1}} 
\def\TTTUQ{\defstr.{2}}

\def\becc{2} 

\def\drv{\becc.{1}} 
\def\lro{\drv.{1}}
\def\fpfm{\drv.{2}}
\def\seqs{\drv.{3}}
\def\rve{\drv.{4}}
\def\rutai{\drv.{5}}
\def\uqi{\drv.{6}}

\def\dvba{\becc.{2}} 
\def\basbec{\dvba.{1}}
\def\eccofi{\dvba.{2}}
\def\firv{\dvba.{3}}
\def\commi{\dvba.{4}}

\def\adt{3} 

\def\rrv{\adt.{1}} 
\def\tgammat{\rrv.{1}} 
                    
\def\ai2a{\rrv.{2}}
\def\ngen{\rrv.{3}}
\def\stngen{\rrv.{4}}

\def\edelta{\adt.{2}} 
\def\tgd1{\edelta.{1}} 
\def\comcon{\edelta.{2}}
\def\tod1{\edelta.{3}}
\def\1coma{\edelta.{4}}
\def\dcoma{\edelta.{5}}

\def\emdelta{\adt.{3}} 
\def\eqct{\emdelta.{1}}
\def\comdind{\emdelta.{2}}

\def\cref{\adt.{4}} 
                    
\def\comefa{\cref.{1}}

\def\comed{\adt.{5}} 
              
\def\nuoed{\comed.{1}}
\def\rnorcomd{\comed.{2}}
\def\comda{\comed.{3}}
\def\comdd{\comed.{4}}

\def\2na{4} 

\def\vf{\2na.{1}} 
\def\homoinj{\vf.{1}}

\def\combrw{\2na.{2}} 
           
\def\lungho{\combrw.{1}}
\def\vudexx{\combrw.{2}}
\def\acvudexx{\combrw.{3}}
\def\tomwn{\combrw.{4}}
\def\iot{\combrw.{5}}
\def\eoeo{\combrw.{6}}
\def\redx{\combrw.{7}}
\def\othusg{\combrw.{8}}

\def\bendef{\2na.{3}} 
\def\relsl32{\bendef.{1}} 
\def\banale{\bendef.{2}} 
\def\stabom{\bendef.{3}} 
\def\m1m21{\bendef.{4}} 
\def\tconcor{\bendef.5} 
\def\am21{\bendef.{6}} 
\def\releed{\bendef.{7}} 
\def\combkkp1{\bendef.{8}} 
\def\comem2f{\bendef.{9}} 
 \def\relee5{\bendef.{10}} 
 \def\welldef{\bendef.{11}} 

\def\tvft{\2na.{4}} 
\def\fvto{\tvft.{1}} 
\def\thpnt{\tvft.{2}} 
\def\access{\tvft.{3}} 
\def\sfzt{\tvft.{4}} 
\def\susf{\tvft.{5}} 
\def\xsind{\tvft.{6}} 
\def\sfilza{\tvft.{7}}
\def\concsflz{\tvft.{8}}
\def\finpar{\tvft.{9}}
\def\rotug{\tvft.{10}}

\def\ij{5} 
\def\trivcomm{\ij.{1}}

\def\gnrl{\ij.{1}} 
\def\bancomm{\gnrl.{1}}
\def\fststp{\gnrl.{2}}

\def\prtclr{\ij.{2}}  
\def\wtlf{\prtclr.{1}}
\def\dn2{\prtclr.{2}}
\def\sndstp{\prtclr.{3}}
\def\linw{\prtclr.{4}}
\def\dadnp1{\prtclr.{5}}
\def\ladnp1{\prtclr.{6}}
\def\adnp1{\prtclr.{7}}
\def\dqt{\prtclr.{8}}
\def\esix{\prtclr.{9}}

\def\comconcl{\ij.{3}} 
\def\simlac{\comconcl.{1}}
\def\cascomm{\comconcl.{2}}

\def\secpbw{6} 

\def\acrd{\secpbw.{1}} 
\def\totord{\acrd.{1}}
\def\cov{\acrd.{2}}
\def\prb{\acrd.{3}}
\def\abnord{\acrd.{4}}

\def\pbls{\secpbw.{2}} 
\def\sqnlt{\pbls.{1}}
\def\pbwt{\pbls.{2}}
\def\lsfr{\pbls.{3}}

\def\sntt{\secpbw.{3}} 
\def\notinf{\sntt.{1}}
\def\stbwg{\sntt.{2}}
\def\descrinf{\sntt.{3}}

\def\rmat{7}

\def\rmgn{\rmat.{1}} 
\def\rta{\rmgn.{1}}
\def\dkrrv{\rmgn.{2}}
\def\dxijr{\rmgn.{3}}
\def\nondegxij{\rmgn.{4}}
\def\dkirv{\rmgn.{5}}
\def\rmaf{\rmgn.{6}}
\def\rcan{\rmgn.{7}}

\def\kirv{\rmat.{2}} 
\def\prcff{\kirv.{1}}
\def\erc{\kirv.{2}}
\def\cij{\kirv.{3}}
\def\kc{\kirv.{4}}

\def\rsnsd{\rmat.{3}} 
\def\linalg{\rsnsd.{1}}
\def\cmt{\rsnsd.{2}}
\def\rfln{\rsnsd.{3}}

\def\teal{\rmat.{4}} 
\def\pvlf{\teal.{1}}
\def\wvlf{\teal.{2}}
\def\sntw{\teal.{3}}
\def\adnd{\teal.{4}}
\def\tkndr{\teal.{5}}
\def\tkdr{\teal.{6}}
\def\sraij{\teal.{7}}

\def\kniv{\rmat.{5}} 
\def\kbef{\kniv.{1}}
\def\kft{\kniv.{2}}

\def\mfrm{\rmat.{6}} \def\emfrm{\mfrm.{1}}

\def\ref{8}

\def\beck{[1]} 
\def\bebr{[2]} 
\def\bourb{[3]} 
\def\cp{[4]} 
\def\dasl2h{[5]} 
\def\rmens{[6]} 
\def\rmcen{[7]} 
\def\dgz{[8]} 
\def\drnry{[9]} 
\def\drqybe{[10]} 
\def\drqg{[11]} 
\def\gmw{[12]} 
\def\hmph{[13]} 
\def\mapi{[14]} 
\def\jmb{[15]} 
\def\jin{[16]} 
\def\jinginv{[17]} 
\def\KM{[18]} 
\def\ktrua{[19]} 
\def\khtsa{[20]} 
\def\kirresh{[21]}
\def\lsrm{[22]} 
\def\lswg{[23]} 
\def\lss{[24]} 
\def\ld{[25]} 
\def\lusz{[26]} 
\def\liq{[27]} 
\def\lqgr{[28]} 
\def\matsu{[29]} 
\def\rRm{[30]} 
\def\tan{[31]}

{\centerline{\bf{THE $R$-MATRIX FOR (TWISTED) AFFINE QUANTUM ALGEBRAS}}} 
\vskip .2 truecm
{\centerline{Ilaria Damiani}}
\vskip .5 truecm
{\hfill{
{\it{To the victims of the NATO aggression against Yugoslavia.}}}}
\vskip .2 truecm
{\it{March $24^{th}$, 1999: NATO starts bombing the
Federal Republic of Yugoslavia, aiming at its civil and productive
structures and killing thousands of people.

May $8^{th}$, 1999: NATO missiles hit the Chinese Embassy in Belgrade,
destroying it and killing 3 people. 
\vskip .2 truecm
I want to dedicate this paper to the Yugoslav students, whose right 
to life and culture
is threatened, and to the Chinese students, who raised their voice against
this criminal war. 
}}

\vskip .5 truecm

{\centerline{\bf{Index}}}
\vskip .3 truecm
0. INTRODUCTION

\genset. GENERAL SETTING: DEFINITIONS AND PRELIMINARIES

\hskip .7 truecm\S \akma. Affine Kac-Moody algebras

\hskip .7 truecm\S \ntz. Some notations, structures and general properties

\hskip 1.4 truecm\rlatt. Root lattice

\hskip 1.4 truecm\sdfscpr. Cartan matrix and symmetric bilinear form

\hskip 1.4 truecm\homz. Weight lattices

\hskip 1.4 truecm\wgrs. Weyl and braid group

\hskip 1.4 truecm\rsys. Root system 

\hskip .7 truecm\S \classdgr. Dynkin diagrams and classification

\hskip .7 truecm\S \defstr. The quantum algebra 

\becc. COPIES OF $A_1^{(1)}$ IN
${X}_{\tn}^{(k)}$ 

\hskip .7 truecm\S \drv. Definition of the root vectors

\hskip .7 truecm\S \dvba. The homomorphisms $\varphi_i$

\adt. THE CASE OF $A_2^{(2)}$

\hskip .7 truecm\S \rrv. Real root vectors

\hskip .7 truecm\S \edelta. The imaginary root vector $\te_{\d}$

\hskip .7 truecm\S \emdelta. The imaginary root vectors $\te_{m\d}$

\hskip .7 truecm\S \cref. Commutation between positive and negative real
root vectors

\hskip .7 truecm\S \comed. The imaginary root vectors $E_{m\d}$

\2na. $A_{2n}^{(2)}$ AND A COPY OF $A_2^{(2)}$

\hskip .7 truecm\S \vf. Definition of $\varphi_1$

\hskip .7 truecm\S \combrw. Root system and Weyl group: some properties

\hskip .7 truecm\S \bendef. $\varphi_1$ is well defined

\hskip .7 truecm\S \tvft. Relations between $\varphi_1$ and the braid group
action

\ij. RELATIONS BETWEEN $\u_q^{(i)}$ AND $\u_q^{(j)}$ ($i\neq j$)

\hskip .7 truecm\S \gnrl. First level of commutation

\hskip .7 truecm\S \prtclr. Some particular computations: low twisted cases

\hskip .7 truecm\S \comconcl. General commutation formulas

\secpbw. A PBW-BASIS 

\hskip .7 truecm\S \acrd. A convex ordering

\hskip .7 truecm\S \pbls. PBW-basis and L-S formula

\hskip .7 truecm\S \sntt. Some subspaces of $\u_q^+$

\rmat. THE $R$-MATRIX

\hskip .7 truecm\S \rmgn. General strategy for the construction of the
$R$-matrix

\hskip .7 truecm\S \kirv. The Killing form on the imaginary root vectors

\hskip .7 truecm\S \rsnsd. A ``canonical'' form for the $R$-matrix

\hskip .7 truecm\S \teal. Linear triangular transformation: the $\bar
E_{\a}$'s

\hskip .7 truecm\S \kniv. The Killing form on the new imaginary root
vectors

\hskip .7 truecm\S \mfrm. The expression of $R$ in terms of the $\bar
E_{\a}$'s

\ref. REFERENCES

\vskip .5 truecm

{\bf{\intr. INTRODUCTION.}}
\vskip .2 truecm
The aim of this paper is to give an exponential multiplicative formula for
the $R$-matrix of the general affine quantum algebra $\uq$. The problem of
describing $R$ has been attacked from different viewpoints, and
different kinds of formulas for
$R$ have already been given explicitly in many (general) cases. 

For the finite type algebras, Rosso gave in \rRm\ the first
exponential multiplicative formula (for $\sl_n$): this result was
generalized in \lsrm, where the strategy of partial $R$-matrices (that
is the study of the connections between the coproduct and the braid
group action) was introduced and developed, and the general (finite) case
solved. 

This exponential approach, mainly based on the description of $R$ via the
Killing form (see \tan), has been extended 
to $A_1^{(1)}$ in \lss, and to the general non
twisted affine case in \rmens.

Other descriptions of $R$ in the affine case (both non twisted and
twisted) are also known, by means of different techniques. Of
particular relevance is the application of the theory of vertex algebras,
which turns out to be a powerful instrument for studying the
representations of the affine quantum algebras; in particular it is
important to mention the works of Khoroshkin-Tolstoy (see \khtsa\ and
\ktrua) which (together with \lss) are the first to face the study of $R$
for the affine algebras, and those of Delius-Gould-Zhang (\dgz) and
Gandenberger-MacKay-Watts (\gmw), where the
$R$-matrix was constructed in the twisted case.

The purpose of the present paper is to complete the ``exponential picture''
by including the twisted affine case in this description. 

To this aim a quite
precise analysis of the twisted affine quantum algebras is needed: since
the PBW-bases (bases of type Poincar\'e-Birkhoff-Witt), and the
possibility to make computations on them, are a fundamental tool for
working in $\uq$ (and in particular for the construction of $R$), the most
important step in this direction is producing a PBW-basis and stating its
main properties by generalizing the  well-known results in the finite and
non twisted affine case.

While approaching this question, one easily finds out that in many
key points the structure of the behaviour of the twisted affine quantum
algebras is not so different (up to some minor adjustments  involving 
no conceptual difficulty and 
requiring
just some care in the
notations; notations which are fixed in section \genset) from what is
already known in the non twisted case: the definition of the root vectors
does not present any difference from the non twisted case (it only happens
that the root system with multiplicities is slightly more complicated to
describe), and the proof of most of the commutation formulas lies on a
property 
of the root system 
($*$: if $\a$ and $\d+r\a$ are roots, then $r=\pm 1$)
which is almost always true; when this is the case the
contribution of this paper is just that of remarking the analogies and
giving references for the proofs. 
Thus, section \becc, where copies of
$\uq(A_1^{(1)})$ are embedded in the affine quantum algebra, has been
written following \beck\ (to which it refers heavily).
The same can be said about paragraphs \S \gnrl\ and \S \comconcl\ (which
are mainly based on \beck) and about section \secpbw\ (where the results
previously found are gathered together to exhibit a PBW-basis). 
Similarly, the first two paragraphs of section \rmat\ (\S \rmgn\ and \S
\kirv, with the exception of lemma \cij), show and use the validity for
the twisted case of the  results of \rmens\ (hence those
of \lsrm\ for what concerns the real root vectors) where the coproduct and
the Killing form on the root vectors are computed.

But there is a case ($A_{2n}^{(2)}$) when $\d-2\a_1$ is a root, so that
the property above denoted by $*$ is no more true. This situation must be
dealt with more carefully, what is done in section \adt\ (where the case
of $A_2^{(2)}$ is discussed in details) and in
section \2na, where $\uq(A_2^{(2)})$ is embedded in $\uq(A_{2n}^{(2)})$:
$A_2^{(2)}$ plays here the same role that $A_1$ plays in the
general frame (to each vertex there is attached a copy of
$\uq(A_1)=\uq(\sl_2)$) and $A_1^{(1)}$ in the generic affine picture (to
``almost each'' vertex of the associated finite diagram there is attached a
copy of $\uq(A_1^{(1)})$). The existence of a
copy of $\uq(A_2^{(2)})$ in $\uq(A_{2n}^{(2)})$ was already proved in
\cp\ by means of the language of vertex operators and through the Drinfeld
realization of $\u_q$ (see \jinginv\ and \jin); here we give a
different, direct proof, which has the advantage to provide a precise
description of the commutation relations that we need. 

Finally, paragraphs \S \prtclr\ (this in particular requires a deeper
understanding of the twisted root systems and Weyl groups), \S \teal\ and
\S \kniv\ are devoted to the computations which, after establishing the
general frame, are needed in order to make explicit the description of
$R$, whose final formula is given in  paragraphs \S \rsnsd\ and \S \mfrm.

I would like to thank Professors Chen Yu, Hu Naihong, Lin Zongzhu, Wang
Jianpan and the other
organizers of the International Conference on Representation Theory, held
in Shanghai (China) in June/July 1998, for giving me the opportunity to
participate in this meeting which has been for me of big interest, and for
offering the possibility to publish a contribution on the Conference
Proceedings. 

Unfortunately last spring, during the work of drawing up of
the paper, the country where I live (Italy) participated in the war against
the Federal Republic of Yugoslavia, which has brutally involved
also the People's Republic of China: 
the decision of the NATO countries to isolate the
scientific community of Yugosalvia is just one - yet hateful - of the
``collateral effects'' of this aggression. Being deeply
concerned by the evident and arrogant injustice of this 
politics, I want to express, as a mathematician, my opposition to this
choice, together with my solidariety to all those who are fighting every
day to make the science survive, even in these hard times, for the present
and future generations. 

\vskip .5 truecm 

{\bf{\genset. GENERAL SETTING: DEFINITIONS AND PRELIMINARIES.}}
\vskip .2 truecm
{\bf{\S \akma. Affine Kac-Moody algebras.}}
\vskip .1 truecm
In this paragraph we recall some basic generalities; for the
proofs and for a deeper and more detailed investigation we refer
to \KM, where these notions have been introduced and studied.

Let $\g$ be a Kac-Moody algebra of finite type.
Then an automorphism of the Dynkin diagram of $\g$ 
of order
$k$ induces an automorphism of $\g$ of the same order (it is well
known and immediate to see that 
$k=$1, 2 or 3).
This automorphism, whose eigenvalues are of course of the form $e^{{2\pi
ir\over k}}$ with
$0\leq r<k$, is diagonalizable, that is
$\g=\bigoplus_{r\in\mz/k\mz}\g^{(r)}$, $\g^{(r)}$ being the eigenspace
relative to $e^{{2\pi ir\over k}}$. It is straightforward to prove
that this decomposition of $\g$ as a direct sum is a
$\mz/k\mz$-grading, that is
$[\g^{(r)},\g^{(s)}]\subseteq\g^{(r+s)}$ $\forall r,s\in\mz/k\mz$;
in particular $\g^{(0)}$ is a Lie subalgebra of $\g$ (indeed a Kac-Moody
algebra of finite type) and
$\g^{(r)}$ is a $\g^{(0)}$-module. 

This remark allows to give a natural structure of Lie algebra to
$$
\bigoplus_{r\in\mz}\g^{([r])}\otimes \mc t^r;$$ it is easy to
construct a non trivial central extension
$$\left(\bigoplus_{r\in\mz}\g^{([r])}\otimes
\mc t^r\right)\oplus\mc c$$ of this algebra via the (non
degenerate) Killing form on $\g$ and a natural structure of Lie algebra on
$${\hat{\g}}^{(k)}\ugdef\left(\bigoplus_{r\in\mz}\g^{([r])}\otimes
\mc t^r\right)\oplus\mc c\oplus\mc\partial,$$
where ${\rm{ad}}\partial= t\rapp{d}{dt}$.

It was proved by Kac (see \KM) that ${\hat{\g}}^{(k)}$ is a Kac-Moody
algebra of affine type depending just on $\g$ and $k$ and not on
the chosen Dynkin diagram automorphism (remark that
${\hat{\g}}^{(1)}$ is 
the non twisted
affine Kac-Moody algebra associated to $\g$; in case that $k>1$ the
algebra ${\hat{\g}}^{(k)}$ is said to be twisted, that is relative to a
non trivial automorphism of the Dynkin diagram); moreover all the affine
Kac-Moody algebras are isomorphic to ${\hat{\g}}^{(k)}$ for some finite
Kac-Moody algebra $\g$ and some $k\in\{1,2,3\}$.

We also have that the Dynkin diagram $\Gamma$ of ${\hat{\g}}^{(k)}$ is an
extension of the Dynkin diagram $\Gamma_0$ of 
$\g^{(0)}$ (which is called the finite diagram - or algebra - associated
to the affine one) 
by adding an extra vertex, which will be denoted by $0$;
remark that if $k=1$ then $\g^{(0)}=\g$.

\vskip .2 truecm

{\bf{\S \ntz. Some notations, structures and general properties.}}

\vskip .1 truecm

Let us fix an affine Kac-Moody algebra ${\hat{\g}}^{(k)}$. 
We shall now introduce some notations.

First of all let us denote by $\tn$ the number of vertices of the
Dynkin diagram of $\g$; this means that $\g$ is of type $X_{\tn}$
for some $X\in\{A,B,C,D,E,F,G\}$, and in this case we say
that ${\hat{\g}}^{(k)}$ is of type $X_{\tn}^{(k)}$. We denote by $n$ the
cardinality of $I_0$, where $I_0$ is the set of vertices of 
$\Gamma_0$
(it means
that the cardinality of
$I$, the set of vertices of $\Gamma$, is $n+1$) and note that 
$n=\tn\Leftrightarrow k=1$.
In particular it is possible to identify
$I$ with the set $\{0,1,...,n\}$ and $I_0$ 
with $\{1,...,n\}$ (so that $I=I_0\cup\{0\}$). One such identification
will be fixed in the next paragraph. 

Furthermore, $A=(a_{ij})_{ij\in I}$ will
denote the Cartan matrix of ${\hat{\g}}^{(k)}$ and $A_0$ will be
$(a_{ij})_{ij\in I_0}$ (the Cartan matrix of $\g^{(0)}$). Remember
that $(I,A)$ determines
${\hat{\g}}^{(k)}$ completely; remark also that $(I,A)$ can be extended to
$(I_{\infty},A_{\infty})$ where $I_{\infty}\ugdef I\cup\{\infty\}$ and
$A_{\infty}=(a_{ij})_{ij\in I_{\infty}}$ with $a_{i\infty}=a_{\infty
i}\ugdef\d_{i0}$ $\forall i\in I_{\infty}$.

We have the following structures associated to $(I,A)$:

{\bf{\rlatt}} Root lattice. 

The root lattice $Q$ and the extended
root lattice $Q_{\infty}$ of ${\hat{\g}}^{(k)}$, and the finite root
lattice
$Q_0$ of $\g^{(0)}$ are defined respectively by:
$$Q\ugdef\oplus_{i\in I}\mz\a_i,\ \ \
Q_{\infty}\ugdef\oplus_{i\in I_{\infty}}\mz\a_i,\
\
\ Q_0\ugdef\oplus_{i\in I_0}\mz\a_i.$$ Of course $Q_0\subseteq
Q=Q_0\oplus\mz\a_0\subseteq Q_{\infty}=Q\oplus\mz\a_{\infty}$. We also
define 
$Q^+$ and
$Q_0^+$ as follows:

$$Q^+\ugdef\sum_{i\in I}\mn\a_i=\left\{\sum_{i\in
I}r_i\a_i\in Q|r_i\geq 0\ \forall i\in I\right\},\ \ \
Q_0^+\ugdef Q_0\cap Q^+.$$ It can be useful to see $Q_{\infty}$, $Q$ and
$Q_0$ as lattices in a vector space over $\mr$ naturally defined as
$Q_{\infty,\mr}\ugdef\mr\otimes_{\mz}Q_{\infty}$. Set also $Q_{\mr}$
and $Q_{0,\mr}$ to be the
$\mr$-subspaces of $Q_{\infty,\mr}$ generated respectively by $Q$ and
$Q_0$.

{\bf{\sdfscpr}} Cartan matrix and symmetric bilinear form. 

The Cartan matrix $A$ is
symmetrizable, that is there exists a diagonal matrix $D={\rm{diag}}(d_i|i\in I)$
such that $DA$ is symmetric; the diagonal entries $d_i$'s can be chosen
to be coprime positive integers, and this condition determines
them uniquely. Remark that if we put $D_{\infty}\ugdef{\rm{diag}}(d_i|i\in
I_{\infty})$ with $d_{\infty}=d_0$ then $D_{\infty}A_{\infty}$ is
symmetric. 

Hence
$D_{\infty}A_{\infty}$ induces a
symmetric bilinear form $(\cdot|\cdot)$ on $Q_{\infty,\mr}$
($\mz$-valued on $Q_{\infty}$) defined by
$(\a_i|\a_j)\ugdef d_ia_{ij}$. We have that 
$(\cdot|\cdot)$ is non degenerate but it is not
positive definite, while 
$(\cdot|\cdot)\big|_{{}_{Q_0\times Q_0}}$ is positive definite
(that is $(\cdot|\cdot)$ induces a structure of Euclidean space on
$Q_{0,\mr}$) and $(\cdot|\cdot)\big|_{{}_{Q\times Q}}$ is degenerate
positive semidefinite: there exists a (unique) non zero element $\d\in Q^+$
such that $\d-\a_0\in Q_0^+$ and $(\d|\d)=0$; more precisely
$(\d|\a)=0$ $\forall\a\in Q$. Remark that $Q$ can be decomposed as
a direct sum also as $Q=Q_0\oplus\mz\d$, and that 
of course $(\a+r\d|\beta+s\d)=(\a|\beta)$ $\forall\a,\beta\in Q$
and $\forall r,s\in\mz$. 

{\bf{\homz}} Weight lattices. 

$(\cdot|\cdot)$ induces
an identification of Hom$(Q_0,\mz)$ with a subgroup $\Pc$ of
$Q_{0,\mr}$.
$\Pc$\ \ is, by the very definition, the lattice $\Pc=\oplus_{i\in
I_0}\mz\omc_i$ where the $\omc_i$'s, called fundamental weights, are
defined as the elements of $Q_{0,\mr}$ characterized by the property that
$(\omc_i|\a_j)\ugdef\d_{ij}$ $\forall j\in I_0$. $Q_0$ is
obviously a sublattice of $\Pc$, and it is worth introducing
another important sublattice $P$ of $\Pc$ as $P\ugdef\oplus_{i\in
I_0}\mz\omega_i$, where $\forall i\in I_0$ $\omega_i\ugdef
d_i\omc_i$; it is to be noticed that $Q_0\subseteq P\subseteq\Pc$.
The elements of $\Pc$ are called weights and the elements of
$$\Pc_+\ugdef\left\{\sum_{i\in I_0}r_i\omc_i\in\Pc|r_i\geq 0\
\forall i\in I_0\right\}$$ are called dominant weights; we also
denote by $P_+$ the set $P_+\ugdef\Pc_+\cap P$. Remark that
$[\Pc:P]=\prod_{i\in I_0}d_i$ and $[P:Q_0]={\rm{det}}(A_0)$, what
can be easily seen noticing that $\forall i\in I_0$ we have
$\a_i=\sum_{j\in I_0}a_{ji}\omega_j$.

It is also important to mention the realization
of $\Pc$ as a group of transformations of $Q$ as
follows: 
define $t:\Pc\rightarrow{\rm{Hom}}(Q,Q)$ by setting
$t_x(\a)\ugdef\a-(x|\a)\d$. Remark that $t$ is an injective
homomorphism of groups and that $\forall x\in\Pc$\ \ $t_x$ preserves 
$(\cdot|\cdot)\big|_{Q\times Q}$ and
fixes
$\d$; equivalently we can say that $t$ maps $\Pc$ into
$\autq$, where $\autq$ is the group of automorphisms of $Q$ preserving
$(\cdot|\cdot)\big|_{Q\times Q}$ and fixing
$\d$.

By abuse of notation we shall
usually denote by $x$ also its image $t_x$. 

{\bf{\wgrs}} Weyl and braid group.

The Weyl group $W$ is the subgroup of
$\autq$  generated by
$\{s_i|i\in I\}$ where
$s_i(\a_j)\ugdef\a_j-a_{ij}\a_i$. $W_0$, the Weyl group of 
$\g^{(0)}$, is the subgroup of $W$ generated by $\{s_i|i\in I_0\}$.

$Aut(\Gamma)$ is the group of automorphisms of the Dynkin
diagram of $\hat{\g}^{(k)}$ (which can be
naturally identified to a subgroup of $\autq$:
$\tau(\a_i)\ugdef\a_{\tau(i)}$). 

The subgroups of $\autq$ introduced above ($\Pc$, $W$, $W_0$,
$Aut(\Gamma)$) have the following properties: 

i) $W\cap Aut(\Gamma)=\{{\rm{id}}\}$, and $\tau s_i\tau^{-1}=s_{\tau(i)}$ 
$\forall i\in I$ $\forall\tau\in Aut(\Gamma)$; hence 
$$W\trianglelefteq\langle W,Aut(\Gamma)\rangle\ \ {\rm{and}}\ \  
\langle W,Aut(\Gamma)\rangle=W\rtimes Aut(\Gamma);$$ of course 
$\langle W,G\rangle=W\rtimes G$ $\forall G\leq Aut(\Gamma)$.

ii) $W_0\cap\Pc=\{{\rm{id}}\}$ and $wxw^{-1}=w(x)$ $\forall x\in\Pc$
$\forall w\in W$; hence $\Pc\trianglelefteq\langle W_0,\Pc\rangle$, so
that 
$\langle W_0,\Pc\rangle=\Pc\rtimes W_0$; as above 
$\langle W_0,L\rangle=L\rtimes W_0$ $\forall L\leq\Pc$ 
\ $W_0$-stable. 

iii) $W\leq\Pc\rtimes W_0$; more precisely 
$$W=\cases Q_0\!\!\check{}\,\,\rtimes W_0&{\rm{in\ the\ non\ twisted\
case}}\cr
\Pc\rtimes W_0&{\rm{in\ case\ }}A_{2n}^{(2)}\cr
Q_0\rtimes W_0 &{\rm{otherwise.}}\endcases$$
($Q_0\subseteq Q_0\!\!\check{}=\oplus_{i\in I_0}\mz\a_i\!\!\check{}
\subseteq\Pc$ \
where 
$\a_i\!\!\check{}={\a_i\over d_i}$).

iv) Let $\T\ugdef Aut(\Gamma)\cap(\Pc\rtimes W_0)$; then 
$$\tW\ugdef W\rtimes\T=\cases \Pc\rtimes W_0&{\rm{in\ the\ non\ twisted\
case\ and\ in\ case\ }}A_{2n}^{(2)}\cr
P\rtimes W_0&{\rm{otherwise.}}\endcases$$
For a precise description of $\T$ see \mapi. 

The length
function $l:\tW\rightarrow\mn$ is defined by 
$$l(\tw)\ugdef{\rm{min}}\{r\in\mn|\exists
i_1,...,i_r\in
I,\tau\in\T\ {\rm{s.t.}}\ \tw=s_{i_1}\cdot...\cdot s_{i_r}\tau\}.$$ 
Also, the braid group $\tB$ is the group
generated by $\{T_{\tw}|\tw\in\tW\}$ with
relations $T_{\tw}T_{\tw^{\prime}}=
T_{\tw\tw^{\prime}}$ whenever 
$l(\tw\tw^{\prime})=l(\tw)+l(\tw^{\prime})$ (see \matsu). Set 
$T_i\ugdef T_{s_i}$.

{\bf{\rsys}} Root system. 

The root system $\Phi\subseteq Q$ divides in two
components: 
$\Phi=\Phi^{{\rm{re}}}\cup\Phi^{{\rm{im}}}$ 
where 
$$\Phi^{{\rm{re}}}\ugdef W.\{\a_i|i\in I\}=
\tW.\{\a_i|i\in I\}=\{\a\in\Phi|(\a|\a)>0\}=\{{\rm{real\ roots}}\}$$ and
$$\Phi^{{\rm{im}}}\ugdef
\{m\d|m\in\mz\setminus\{0\}\}=\{\a\in\Phi|(\a|\a)=0\}=\{{\rm{imaginary\
roots}}\}.$$  

$\Phi^{{\rm{re}}}$ can be described in terms of 
$\Phi_0\ugdef W_0.\{\a_i|i\in I_0\}$ (the root
system of $\g^{(0)}$) as follows: 
$$\Phi^{{\rm{re}}}=\cases
\{m\d+\a|\a\in\Phi_0,m\in\mz\}&{\rm{in\ the\ non\ twisted\ case}}\cr
\{m\d+\a|\a\in\Phi_0,m\in\mz\}\cup\cr
{\phantom{a}}\cup\{(2m+1)\d+2\a|\a\in\Phi_0,(\a|\a)=2,m\in\mz\}
&{\rm{in\
case\ }}A_{2n}^{(2)}\cr
\{md_{\a}\d+\a|\a\in\Phi_0,m\in\mz\}
&{\rm{otherwise}},\endcases$$
where $d_{\a}\ugdef\rapp{(\a|\a)}{2}$ (or
equivalently $d_{\tw(\a_i)}=d_i$ $\forall i\in I,\ \tw\in\tW$); it is worth
remarking that all the subgroups of $\autq$ considered in \homz\ and
\wgrs\ leave $\Phi$ stable.

The multiplicity of each real root is 1; on the
other hand let us define, $\forall i\in I_0$,
$$\dk_i\ugdef\cases 1&{\rm{in\ the\ non\ twisted\ case\ or\ in\ case\ }}
A_{2n}^{(2)}\cr
d_i&{\rm{otherwise}}\ \endcases\ \ {\rm{and}}\ \
\tk\ugdef{\rm{max}}\{\dk_i|i\in I_0\},$$  
and set $I^m\ugdef\{i\in I_0|\dk_i|m\}$; 
then the multiplicity of
$m\d$ is 
$$\#
I^m=\cases n&{\rm{if}}\ k|m\cr{\tn-n\over k-1}&{\rm{if}}\
k\not|m,\endcases$$ which can be expressed also by saying that the
root system with multiplicities is 
$$\tPhi\ugdef\Phi^{{\rm{re}}}\cup
\tPhi^{{\rm{im}}}\ \ \ {\rm{where}}\ 
\tPhi^{{\rm{im}}}\ugdef
\{(m\d,i)|m\in\mz\setminus\{0\},\dk_i|m\}$$
(Denote by $p:\tPhi\rightarrow\Phi$ the natural map).

Remark that for $m\neq 0$ we have
$$(m\d,i)\in\tPhi^{{\rm{im}}}\Leftrightarrow
\dk_i|m\Leftrightarrow
m\d+\a_i\in\Phi
\Leftrightarrow m\omc_i\in\tW.$$ 

There is also another decomposition of $\Phi$: $\Phi=\Phi_+\cup\Phi_-$,
where 
$\Phi_+\ugdef\Phi\cap Q^+=$ $=\{{\rm{positive\ roots}}\}$ and
$\Phi_-\ugdef-\Phi_+=\{{\rm{negative\ roots}}\}$. It induces an analogous
decomposition of $\tPhi$ and of its real and imaginary parts, with obvious
notations.

Recall also the relation between the length
function and the root system: 
$l(\tw)=$ $=\#\Phi_+(\tw)$ where 
$$\Phi_+(\tw)\ugdef
\{\a\in\Phi_+|\tw^{-1}(\a)<0\}:$$
recall that if $s_{i_1}\cdot...\cdot s_{i_r}$ is a reduced
expression then 
$$\Phi_+(s_{i_1}\cdot...\cdot s_{i_r})=
\{s_{i_1}\cdot...\cdot s_{i_{k-1}}(\a_{i_k})|1\leq k\leq
r\}.$$
Finally,
$\l(\tw\tw^{\prime})=l(\tw)+l(\tw^{\prime})
\Leftrightarrow\Phi_+(\tw)\subseteq
\Phi_+(\tw\tw^{\prime})$. 

\vskip .2 truecm
{\bf{\S \classdgr. Dynkin diagrams and classification.}}
\vskip .1 truecm

In the following we list the affine Dynkin diagrams. The
labels under the vertices fix an identification between $I$
and $\{0,1,...,n\}$ such that $I_0$ corresponds to
$\{1,...,n\}$. For each type we also recall the coefficients
$r_i$ (for $i\in I_0$: recall that $r_0$ is always 1) in the
expression
$\d=\sum_{i\in I}r_i\a_i$. 

$$\sm(\Gamma,I)\vm\ \ \ 
(r_1,...,r_n)\ \leqno{X_{\tn=\tn(n)}^{(k)}\
\ \ \  n}$$
${\underline{{\phantom{\hskip 13.5 truecm}}}}$
$$\sm\vx{0}\dedg\vx{1}\om\ \ \ \ \ \ 
\nsm(1)\nstm
\leqno{\ \ A_1^{(1)}\ \ \ \ \ \ \ \ 1}
$$
$$\tm\vx{1}\,^{^{\diagup\!\!^{^{\diagup\!\!^{^{\diagup
}}}}}}
\!\!\!\!\!\!\!\!\!\!\!\edg\vx{2}\,.\,.\!
\!^{^{
{\phantom{\Big|}}
\!\!\!^{^{^{^{^{\vx{0}\!\!\!^{}}}}}}}}\!\!
\,\,.\vl{n-1}\edg\redg
\!\!\!^{^{\diagdown\!\!\!\!\!\!\!\!^{^{\diagdown
\!\!\!\!\!\!\!\!^{^{\diagdown}}}}}}
\,\,\,\vn
\cm\ \ \ \ 
\ndcm(1,...,1)\leqno{\ \
A_n^{(1)}\ \ \ \
\ >1}
$$
$$\tm\vx{1}\!<\!\!\!\dedg\vx{2}\edg\vx{3}\,...\vl{n-2}\edg\vl{n-1}
\vedgv{0}{}\edg\vn\qm\ \ \ \ \ 
\ndm\!(2,...,2,1)\ncm
\leqno{\ \ B_n^{(1)}\ \ \ \ \ 
>2}
$$
$$\dm\vx{1}\dedg\!\!>\!\!\vx{2}\edg\vx{3}\,...\vl{n-1}\edg
\vn\!<\!\!\!\dedg\vx{0}\qm\ \ \ \ \ \ 
(1,2,...,2)\ncm
\leqno{\ \ C_n^{(1)}\ \ \ \ \ 
>1}
$$
$$\dm\vx{2}\edg\vx{3}\vedgv{1}{}
\edg\vx{4}\,...\vl{n-2}
\edg\vl{n-1}\vedgv{0}{}\edg\vn\tm\ \ \ \ \ \ \ \ 
(1,1,2,...,2,1)\num
\leqno{\ \ D_n^{(1)}\ \ \ \ \ 
>3}
$$
$$\tm\vx{2}\edg\vx{3}\edg\vx{4}
\vedgv{1}{\vedgv{0}{}}\!\edg\vx{5}\edg\vx{6}\tm\
\ \ \ \ \ \ntm(2,1,2,3,2,1)\num
\leqno{\ \ E_6^{(1)}\ \ \ \ \ \ \ \ 
6}
$$
$$\um\vx{0}\edg\vx{2}\edg\vx{3}\edg\vx{4}\vedgv{1}{}
\edg\vx{5}\edg\vx{6}\edg\vx{7}\um\
\ \ \ \ \ \num(2,2,3,4,3,2,1)
\leqno{\ \ E_7^{(1)}\ \ \ \ \ \ \ \ 
7}
$$
$$\vx{2}\edg\vx{3}\edg\vx{4}\vedgv{1}{}\edg\vx{5}\edg\vx{6}
\edg\vx{7}\edg\vx{8}\edg\vx{0}\ \
\ \ \ \ (3,2,4,6,5,4,3,2)
\leqno{\ \ E_8^{(1)}\ \ \ \ \ \ \ \ 
8}
$$
$$\um\vx{1}\edg\vx{2}\!<\!\!\!\dedg
\vx{3}\edg\vx{4}\edg\vx{0}\qm\ \ \ \ 
\ncm(2,4,3,2)\ndm
\leqno{\ \ F_4^{(1)}\ \ \ \ \ \ \ \ 
4}
$$
$$\tm\vx{1}\!<\!\!\!\tedg\vx{2}\edg\vx{0}\stm\
\ \ \ \nsm(3,2)\nsm
\leqno{\ \ G_2^{(1)}\ \ \ \ \ \ \ \ 
2}
$$
$$\qm
\vx{1}\!<\!\!\!\qedg\vx{0}\om\ \ \ 
\nstm(2)\nstm
\leqno{\ \ A_2^{(2)}\ \ \ \ \ \ \ \ 
1}
$$
$$\um\vx{1}\!<\!\!\!\dedg\vx{2}\edg\vx{3}\,...
\vl{n-1}
\edg\vn\!<\!\!\!\dedg\vx{0}\qm\ \ \ \ 
\ndm(2,...,2)\nstm
\leqno{\ \ A_{2n}^{(2)}\ \ \ \ \ 
>1}
$$
$$\dm\vx{1}\dedg\!\!>\!\!\vx{2}\edg\vx{3}\,...\vl{n-2}\edg\vl{n-1}
\vedgv{0}{}\edg\vn\tm\ \ \ \ \ 
\ntm(1,2,...,2)\nnm
\leqno{\ \ A_{2n-1}^{(2)}\ \ \ 
>2}
$$
$$\um\vx{1}\!<\!\!\!
\dedg\vx{2}\edg\vx{3}\,...\vl{n-1}
\edg\vn\dedg\!\!>\!\!\vx{0}\tm\ \ \ \ 
\ncm(1,...,1)\nstm
\leqno{\ \ D_{n+1}^{(2)}\ \ \ \ 
>1}
$$
$$\dm\vx{0}\edg\vx{1}\edg
\vx{2}\!<\!\!\!\dedg\vx{3}\edg\vx{4}\qm\
\ \ \ 
\nqm(2,3,2,1)\ncm
\leqno{\ \ E_6^{(2)}\ \ \ \ \ \ \ \ 
4}
$$
$$\tm\vx{0}\edg\vx{1}\!<\!\!\!\tedg\vx{2}\stm\
\ \ \ \nsm(2,1)\nsm
\leqno{\ \ D_4^{(3)}\ \ \ \ \ \ \ \ 
2}
$$

\vskip .2 truecm
{\bf{\S\defstr. The quantum algebra.}}
\vskip .2 truecm
In this paragraph we  recall the definition of the quantum group
$\uq
$ 
associated to the affine
Kac-Moody algebra
$\hat\g^{(k)}$ of type $X_{\tn}^{(k)}$, and the main structures that
$\uq$ can be endowed with. 
\ddefi{\QNTALG}Let $\hat\g^{(k)}$ be the affine Kac-Moody algebra of
type $X_{\tn}^{(k)}$. 

We denote 
by $\uq=\uq(\hat\g^{(k)})=\uq(X_{\tn}^{(k)})$ 
the quantum algebra of $\hat\g^{(k)}$, that is the associative
$\mc(q)$-algebra generated by
$$\{E_i,F_i,K_i^{\pm 1},D^{\pm 1}=K_{\infty}^{\pm 1}|i\in I\}$$ 
with relations: 
$$K_{\lambda}K_{\mu}=K_{\lambda+\mu}\ \ 
\forall\lambda,\mu\in Q_{\infty},$$ 
$$K_{\lambda}E_i=q^{(\lambda|\alpha_i)}E_iK_{\lambda}\ \ 
\forall\lambda\in Q_{\infty},\forall i\in I,$$ 
$$K_{\lambda}F_i=q^{-(\lambda|\alpha_i)}F_iK_{\lambda}\ \ 
\forall\lambda\in Q_{\infty},\forall i\in I,$$ 
$$[E_i,F_j]=\delta_{ij}\rapp{K_i-K_i^{-1}}{q_i-q_i^{-1}}\ \ 
\forall i,j\in I$$ 
$$\sum_{r=0}^{1-a_{ij}}{1-a_{ij}\brack r}_{q_i} 
E_i^rE_jE_i^{1-a_{ij}-r}=0=\sum_{r=0}^{1-a_{ij}}{1-a_{ij}\brack r}_{q_i} 
F_i^rF_jF_i^{1-a_{ij}-r}\ \ \forall i\neq j\in I,$$ 
where we set:

i) $K_{\lambda}\ugdef \prod_{i\in I_{\infty}}K_i^{m_i}$
if $\lambda=\sum_{i\in I_{\infty}}m_i\alpha_i\in Q_{\infty}$; 

ii) $q_{\alpha}\ugdef q^{d_{\a}}$ 
$\forall\alpha\in\Phi^{{\rm{re}}}$, 
$q_i\ugdef q_{\alpha_i}$ $\forall i\in I$, 
$q_{(m\delta,i)}\ugdef q_i$ $\forall (m,i)\in\tPhi_+^{{\rm{im}}}$; 

iii) $[m]_{q^{r}}\ugdef\rapp{q^{rm}-q^{-rm}}{q^r-q^{-r}}$ 
and
$[m]_{q^{r}}!\ugdef 
\prod_{s=1}^{m}[s]_{q^{r}}$ $\forall m,r\in\mz$;

moreover 
${r\brack s}_{q^m}={[r]_{q^m}!\over[s]_{q^m}![r-s]_{q^m}!}$
$\forall m\in\mz,\forall s\leq r\in\mn$; 

iv) for further use let us also define $\forall\alpha\in\tilde\Phi_+$,
$\forall x\in\uq$ 
$${\rm{exp}}_{\alpha}(x)\ugdef 
\sum_{m\geq 0}\rapp{x^m}{(m)_{\alpha}!},$$ 
where $\forall\alpha\in\tilde\Phi_+, \forall
m\in\mn$ we have put
$$(m)_{\alpha}\ugdef\cases 
\rapp{q_{\alpha}^{2m}-1}{q_{\alpha}^2-1}&{\rm{if}}\  
\alpha\in\Phi_+^{{\rm{re}}}\cr 
m&{\rm{if}}\  
\alpha\in\tPhi_+^{{\rm{im}}}\endcases 
\ \ \ {\rm{and}}\ \ \ (m)_{\alpha}!\ugdef 
\prod_{s=1}^{m}(s)_{\alpha}.$$ 
\rem{\TTTUQ}
Recall that on 
$\uq$ we have the following structures: 

1) the $Q$-gradation $\uq=\oplus_{\eta\in Q}\u_{q,\eta}$ determined by the
conditions that
$E_i\in\u_{q,\a_i}$ and $F_i\in\u_{q,-\a_i}$ $\forall i\in I$, $K_i^{\pm
1}\in\u_{q,0}$ $\forall i\in I_{\infty}$ and 
$\u_{q,\a}\u_{q,\beta}\subseteq\u_{q,\a+\beta}$ $\forall\a,\beta\in Q$; 

2) the triangular decomposition: 
$\uq\cong\u_q^-\otimes\u_q^0\otimes\u_q^+\cong 
\u_q^+\otimes\u_q^0\otimes\u_q^-$, where 
$\u_q^-$, $\u_q^0$ and $\u_q^+$ are the subalgebras of $\uq$ 
generated respectively by 
$\{E_i|i\in I\}$, $\{K_i^{\pm 1}|i\in I_{\infty}\}$ and 
$\{F_i|i\in I\}$; define also $\u_q^{\geq 0}\ugdef(\u_q^+,\u_q^0)$,
$\u_q^{\leq 0}\ugdef(\u_q^-,\u_q^0)$;

3) the $\mc$-anti-linear anti-involution 
$\Omega:\uq\rightarrow\uq$ defined by 
$$\Omega(E_i)\ugdef F_i,\ \  
\Omega(F_i)\ugdef E_i\ \forall i\in I;\ \ 
\Omega(K_i)\ugdef K_i^{-1}\ \forall i\in I_{\infty};\ \  
\Omega(q)\ugdef q^{-1};$$ 

4) the $\mc(q)$-linear anti-involution 
$\Xi:\u_q\rightarrow\u_q$ defined
by
$$\Xi(E_i)\ugdef E_i,\ \  
\Xi(F_i)\ugdef F_i\ \forall i\in I;\ \ 
\Xi(K_i)\ugdef K_i^{-1}\ \forall i\in I_{\infty};$$ 

5) a braid group action commuting with $\Omega$, verifying 
$\Xi T_i=T_i^{-1}\Xi$ $\forall i\in I$, and preserving the
gradation (that is $T_w(\u_{q,\eta})=\u_{q,w(\eta)}$); 
more precisely 
$$T_w(E_i)\in\u_{q,w(\a_i)}^+\ {\rm{if}}\ w\in W\ {\rm{and}}
\ i\in I\ {\rm{are\ such\ that}}\ w(\a_i)\in Q^+\ ({\rm{i.e.}}\
l(ws_i)>l(w)),$$ and
$$T_i(E_i)=-F_iK_i\ \ \forall i\in I,\ \ \ 
T_i(E_j)=\sum_{r=0}^{-a_{ij}}(-1)^{r-a_{ij}}q_i^{-r}
E_i^{(-a_{ij}-r)}E_jE_i^{(r)}\ \ \forall i\neq j\in I$$
where $\forall m\in\mn$ $E_i^{(m)}\ugdef{E_i^m\over[m]_{q_i}!};$ 

6) a Hopf-structure ($\Delta$, $\epsilon$, $S$), 
whose coproduct $\Delta:\uq\rightarrow\uq\otimes\uq$ 
is given by 
$$\Delta(E_i)=E_i\otimes 1+K_i\otimes E_i,$$  
$$\Delta(F_i)=1\otimes F_i+F_i\otimes K_i^{-1},$$  
$$\Delta(K_i)=K_i\otimes K_i;$$ 
Let us remark that $\Delta\compo\Omega= 
\sigma\compo(\Omega\otimes\Omega)\compo\Delta$ 
where $\sigma:\uq\otimes\uq\rightarrow\uq\otimes\uq$ is defined by
$\sigma(x\otimes y)\ugdef y\otimes x$; 

7) the Killing form, that is the unique bilinear form 
$(\cdot,\cdot): 
\u_q^{\geq 0}\otimes\u_q^{\leq 0}\rightarrow\mc(q)$ such that 

i) $(x,y_1y_2)=(\Delta(x),y_1\otimes y_2)$ 
$\forall x\in\u_q^{\geq 0}$, $\forall y_1,y_2\in\u_q^{\leq 0}$; 

ii) $(x_1x_2,y)=(x_2\otimes x_1,\Delta(y))$ 
$\forall x_1,x_2\in\u_q^{\geq 0}$, $\forall y\in\u_q^{\leq 0}$; 

iii) $(\cdot,\cdot)\big|_{\u_{q,\eta}^{\geq 0}\times 
\u_{q,-\tilde\eta}^{\leq 0}}=0$ if $\eta\neq\tilde\eta$; 

iv) $(K_{\lambda},K_{\mu})=q^{-(\lambda,\mu)}$ 
$\forall\lambda,\mu\in Q_{\infty}$; 

v) $(E_i,F_i)=\rapp{1}{q_i^{-1}-q_i}$ $\forall i\in I$; 

$(\cdot,\cdot)$ has the following important properties: 

a) $(\cdot,\cdot)\big|_{\u_{q,\eta}^+\times 
\u_{q,-\eta}^-}$ is non degenerate $\forall\eta\in Q_+$; 

b) $(xK_{\lambda},y)=(x,y)=(x,yK_{\lambda})$ 
$\forall x\in\u_q^+$, $\forall y\in\u_q^-$, 
$\forall\lambda\in Q_{\infty}$. 

\vskip .5 truecm

{\bf{\becc. COPIES OF $A_1^{(1)}$ IN
${X}_{\tn}^{(k)}$. 
}}
\vskip .2 truecm

{\bf{\S \drv. Definition of the root vectors.}}

\vskip .1 truecm

In this section we shall shortly recall Beck's
results
for non-twisted affine
quantum algebras (see \beck\ and \bebr) and prove that his method
applies word by word to the twisted case, that is
we can find a copy of $\u_q(A_1^{(1)})$ in 
$\u_q({X}_{\tn}^{(k)})$ $\forall i\in I_0$,
provided that $({X}_{\tn}^{(k)},i)\neq
(A_{2n}^{(2)},1)$. 

To this aim let us first define, $\forall i\in I_0$,
an element 
$\lambda_i\in\Pc$ as follows: 
$$\lambda_i\ugdef \dk_i\omc_i$$ 
where we recall that $\dk_i={\rm{min}}\{
m>0|m\omc_i\in\tW\}$ (see \rsys). 

The next proposition and remark is what one
needs in order to define the positive root
vectors. 
\prop{\lro}
$\forall m\in\mn$ we have
$l((\lambda_1\cdot...\cdot \lambda_n)^m)=
m\sum_{i\in I_0}l(\lambda_i)$.
\dim 
See \bourb.\finedim
\rem{\fpfm}
$$\Phi_+^{{\rm{re}}}=\{m\d-\a\in\Phi_+|\a\in
Q_0^+,m>0\}\cup\{m\d+\a\in\Phi_+|\a\in
Q_0^+,m\in\mn\}.$$
Moreover
$$\{m\d-\a\in\Phi_+|\a\in
Q_0^+,m>0\}=
\bigcup_{m>0}\Phi_+((\lambda_1\cdot...\cdot \lambda_n)^m),$$
$$\{m\d+\a\in\Phi_+|\a\in
Q_0^+,m\in\mn\}=
\bigcup_{m>0}\Phi_+((\lambda_1\cdot...\cdot \lambda_n)^{-m}).
$$
\dim 
See \bebr.
\finedim
The preceding results suggest to define the following
sequence, by means of which it will be easy to 
give a suitable definition of the positive real root 
vectors. 
\ddefi{\seqs}
Let $N\ugdef l(\lambda_1\cdot...\cdot \lambda_n)$ and define
$\iota:\mz\rightarrow I_0$ to be such that 
$$\cases
\lambda_1\cdot...\cdot \lambda_n=s_{\iota_1}\cdot...\cdot s_{\iota_N}\tau\cr
\iota_{r+N}=\tau(\iota_r)\ \ \ \forall r\in\mz,
\endcases$$ where $\tau$ is a suitable element of
$\T$ (uniquely determined by the above condition); for further
determinations of $\iota$ see remark \rutai\ and notation \iot.

Notice that $\iota$ induces functions 
$$\mz\ni r\mapsto w_r\in W\ \ \ {\rm{and}}\ \ \ 
\mz\ni r\mapsto\beta_r\in\Phi_+^{{\rm{re}}}$$
by 
$$w_r\ugdef\cases s_{\iota_1}\cdot...\cdot s_{\iota_{r-1}}
{\rm{if}}\ r\geq 1\cr
s_{\iota_0}\cdot...\cdot s_{\iota_{r+1}}
{\rm{if}}\ r\leq 0,\endcases\ \ \ \ \
{\rm{and}}\ \ \ \ \ \beta_r\ugdef w_r(\a_{\iota_r}).$$
It is well known (see \bebr) that $\beta$ is a bijection. 

Then we arrive at the definition of the root
vectors: 
\ddefi{\rve}
The positive real root vectors
$E_{\a}$ (with $\a\in\Phi_+^{{\rm{re}}}$) are defined
by 
$$E_{\a}\ugdef\cases
T_{w_r}(E_{\iota_r})&{\rm{if}}\ r\geq 1\cr
T_{w_r^{-1}}^{-1}(E_{\iota_r})=
\Xi T_{w_r}(E_{\iota_r})&{\rm{if}}\ r\leq
0.\endcases$$
Also, positive imaginary root vectors 
$\te_{(m\dk_i\d,i)}$ (with $m>0$, $i\in I_0$) are
defined by 
$$\te_{(m\dk_i\d,i)}\ugdef-E_{m\dk_i\d-\a_i}E_i+
q_i^{-2}E_iE_{m\dk_i\d-\a_i}.$$
Similarly, the negative root vectors are defined
by
$$F_{\a}\ugdef\Omega(E_{\a})\ \ \ {\rm{if}}\
\a\in\Phi_+^{{\rm{re}}},\ \ \
\tf_{\a}\ugdef\Omega(\te_{\a})\ \ \ {\rm{if}}\
\a\in\tPhi_+^{{\rm{im}}}.$$
\rem{\rutai} 
It is worth remarking, for later use, that 
$\iota$ can be chosen (and it will be chosen!) so that 
$\exists N_0=0<N_1<...<N_n=N$ and $\exists\tau_1,...,\tau_n\in\T$
such that $\forall i=1,...,n$  
$\lambda_i=s_{\iota_{N_{i-1}+1}}\cdot...\cdot s_{\iota_{N_i}}\tau_i$
(notice that $\tau_1\cdot...\cdot\tau_n=\tau$). 

If this is the case then 
$$E_{m\dk_i\d+\a_i}=T_{\lambda_i}^{-m}(E_i)\ \ \ 
\forall m\in\mn,$$
$$E_{m\dk_i\d-\a_i}=T_{\lambda_i}^mT_i^{-1}(E_i)
\ \ \ 
\forall m>0,$$
which depends on the fact that 

1) $T_{\lambda_j}T_{\lambda_i}=
T_{\lambda_i}T_{\lambda_j}$ $\forall i,j\in I_0$;

2) $T_{\lambda_j}(E_i)=E_i$ and 
$T_{\lambda_j}T_i=T_iT_{\lambda_j}$ 
$\forall i,j\in I_0$ such that $j\neq i$.
\dim 
See \ld\ and \rmens.\finedim
\ddefi{\uqi}
If $i\in I_0$ we define $\u_q^{(i)}$ to be the $\mc(q^{d_i})$-subalgebra
of
$\u_q=\u_q({X}_{\tn}^{(k)})$ generated by
$\{E_i,E_{\dk_i\d-\a_i},F_i,F_{\dk_i\d-\a_i},K_i^{\pm
1},K_{\dk_i\d}^{\pm 1}\}$. 

Obviously $\u_q^{(i)}$ is $\Omega$-stable and pointwise fixed by
$T_{\lambda_j}$ for $j\in I_0\setminus\{i\}$. Moreover  it will become
clear that
$\u_q^{(i)}$ is the smallest 
$\mc(q^{d_i})$-subalgebra
of
$\u_q$ containing $E_i$ and $K_i$ and stable under $T_i^{\pm 1}, 
T_{\lambda_i}^{\pm 1}$; in particular we shall see that it contains 
$E_{m\dk_i\d\pm\a_i}$,
$F_{m\dk_i\d\pm\a_i}$, $\te_{(m\dk_i\d,i)}$ and 
$\tf_{(m\dk_i\d,i)}$ $\forall m>0$.
\vskip .2 truecm

{\bf{\S \dvba. The
homomorphisms $\varphi_i$.}}

\vskip .1 truecm
We are now ready to state the announced properties of
the root vectors just defined and of the braid 
group action, which will allow us to construct
homomorphisms 
$\varphi_i:\u_q(A_1^{(1)})\rightarrow
\u_q^{(i)}\subseteq\u_q( X_{\tn}^{(k)})$ for each $i\in I_0$, under the
only condition that 
$i\neq 1$ if $ X_{\tn}^{(k)}=A_{2n}^{(2)}$,
condition that we shall assume in the remaining part of
section \becc. These homomorphisms $\varphi_i$ behave
``well'' on the root vectors (that is, they
transform root vectors in root vectors, see corollary \firv), and
this fact immediately implies some important
consequences: in particular we shall translate the
commutation relations from 
$\u_q(A_1^{(1)})$ to $\u_q( X_{\tn}^{(k)})$. We
shall also underline the obstruction that we meet if
we try to apply the same argument when
$X_{\tn}^{(k)}=A_{2n}^{(2)}$ and 
$i=1$: this last situation will be studied in sections 
\adt\ and \2na.
\lem{\basbec}
Assume, as required, that $X_{\tn}^{(k)}\neq
A_{2n}^{(2)}$  or $i\neq 1$; then we have: 

1) $l(\lambda_is_i\lambda_i)=2l(\lambda_i)-1$; 

2) $(T_{\lambda_i}T_i^{-1})^2
\in\langle
T_{\lambda_j}|j\neq i\rangle$; 

3) $[E_{\dk_i\d-\a_i},F_i]=0$ 
and $E_iE_{\dk_i\d+\a_i}=q_i^{-2}E_{\dk_i\d+\a_i}E_i$.

4) $\forall j\in I_0\setminus\{i\}$ we have
$(T_{\lambda_i}T_i^{-1})^2T_{\lambda_j}^{{\dk_ia_{ij}\over \dk_j}}
\in\langle T_{\lambda_r}|r\in I_0
\setminus\{i,j\}\rangle$; in particular 
$T_{\lambda_i}T_i^{-1}(E_j)=T_{\lambda_j}^{-{\dk_ia_{ij}\over
\dk_j}}T_i(E_j)$.
\dim
For the proof see \ld\ and \beck. Here we
want just to point out that all the statements are
based on the fact that 
there is no root $\a$ of the form 
$\a=m\d+\varepsilon\a_i$ 
such that 
$\dk_i\d+\a_i\prec\a\prec\a_i$.

Remark that this is no longer true in case
$A_{2n}^{(2)}$ for $i=1$: indeed in this case we have 
$\d+\a_i\prec\d+2\a_i\prec\a_i$. Hence the existence 
of roots of the form $m\d+2\a$ with 
$\a\in\Phi_+\cap Q_0$ (which never happens if we are not
in case $A_{2n}^{(2)}$)
is what makes the
difference between this case and the others.\finedim
Since Beck's construction of the homomorphisms 
$\varphi_i:\u_q(A_1^{(1)})\rightarrow\u_q^{(i)}\subseteq\u_q(
X_n^{(1)})$ (and their properties)
is based on the statements of lemma \basbec\, the
argument used in \beck\ applies exactly to the present
situation, so that we have: 
\prop{\eccofi}
If $ X_{\tn}^{(k)}\neq A_{2n}^{(2)}$ 
or $i\neq 1$, then there exists
a $\mc$-algebra homomorphism 
$$\varphi_i:\u_q(A_1^{(1)})\rightarrow
\u_q^{(i)}\subseteq\u_q(X_{\tn}^{(k)})$$ verifying the following
properties:

1) $\varphi_i(K_1)=K_i$,
$\varphi_i(K_1K_0)=K_{\dk_i\d}$ and
$\varphi_i(q)=q_i=q^{d_i}$;  

2) $\varphi_i(E_1)=E_i$ and 
$\varphi_i(E_0)=E_{\dk_i\d-\a_i}$;

3) $\varphi_i\Omega=\Omega\varphi_i$; 

4) $\varphi_iT_1=T_i\varphi_i$ and 
$\varphi_iT_{\tau}=
T_{\lambda_i}T_i^{-1}\varphi_i$, where $\tau$ is the
non trivial Dynkin diagram automorphism of 
$A_1^{(1)}$.\finedim
\cor{\firv} 
Suppose again that 
$ X_{\tn}^{(k)}\neq A_{2n}^{(2)}$ 
or $i\neq 1$ and 
let $\nu_i$ be the group homomorphism
$\nu_i:Q(A_1^{(1)})\rightarrow Q(X_{\tn}^{(k)})$ 
defined in the obvious way
by 
$\nu_i(\a_1)\ugdef\a_i$, $\nu_i(\d)\ugdef \dk_i\d$;
furthermore let
the injection
$\tnu_i:\Phi_+(A_1^{(1)})\rightarrow
\tPhi_+( X_{\tn}^{(k)})$
be given by: 
$$\tnu_i(\a)\ugdef\cases\nu_i(\a)&{\rm{if\
}}\a{\rm{\ is\ real}}\cr(\nu_i(\a),i)&{\rm{if\
}}\a{\rm{\ is\ imaginary;}}\endcases$$
Then $\varphi_i$ has the property that 
$$\varphi_i(K_{\a})=K_{\nu_i(\a)}\ \
\forall\a\in Q(A_1^{(1)}),$$  
$$\varphi_i(E_{\a})=E_{\tnu_i(\a)}\ \
\forall\a\in\Phi_+(A_1^{(1)})$$  and
$$\varphi_i(\te_{\a})=\te_{\tnu_i(\a)}\ \ \
\forall\a\in\Phi_+^{{\rm{im}}}(A_1^{(1)}),$$
where the vectors $E_{(m\dk_i\d,i)}$ are defined by the
relation
$$1-(q_i-q_i^{-1})\sum_{m>0}\te_{(m\dk_i\d,i)}u^m=
{\rm{exp}}\left((q_i-q_i^{-1})
\sum_{m>0}E_{(m\dk_i\d,i)}u^m\right).$$\finedim
As an immediate consequence of corollary \firv, we
have the following proposition.
\prop{\commi}
Suppose as before that 
$(X_{\tn}^{(k)},i)\neq(A_{2n}^{(2)},1)$; then we have
the following relations in 
$\u_q(X_{\tn}^{(k)})$: 

1) $[\te_{(r\dk_i\d,i)},\te_{(s\dk_i\d,i)}]=0=
[E_{(r\dk_i\d,i)},E_{(s\dk_i\d,i)}]$ 
$\forall r,s>0$; 

2) $T_{\lambda_j}(\te_{(r\dk_i\d,i)})=\te_{(r\dk_i\d,i)}$
and $T_{\lambda_j}(E_{(r\dk_i\d,i)})=E_{(r\dk_i\d,i)}$ 
$\forall r>0, j\in I_0$; 

3) $[E_{(r\dk_i\d,i)},F_{(s\dk_i\d,i)}]=\d_{rs}
\rapp{[2r]_{q_i}}{r}
\rapp{K_{r\dk_i\d}-K_{r\dk_i\d}^{-1}}{q_i-q_i^{-1}}$ 
$\forall r,s>0$;

4) $[E_{(r\dk_i\d,i)},E_{s\dk_i\d+\a_i}]= 
\rapp{[2r]_{q_i}}{r}
E_{(r+s)\dk_i\d+\a_i}$ and

$[E_{(r\dk_i\d,i)},E_{(s+1)\dk_i\d-\a_i}]=
-\rapp{[2r]_{q_i}}{r} E_{(r+s+1)\dk_i\d-\a_i}$ 
$\forall r>0,s\geq 0$.
\dim 
The claim follows applying $\varphi_i$ to the analogous
relations holding in $\u_q(A_1^{(1)})$ (see \dasl2h\
and \beck).\finedim

\vskip .5 truecm

{\bf{\adt. THE CASE OF $A_2^{(2)}$.}}
\vskip .2 truecm 
In this section we study the case of $\sl_3^{(2)}=A_2^{(2)}$ (and $\u_q$ 
will indicate $\u_q(\sl_3^{(2)}))$. 
In particular we want to understand some important properties of the
imaginary root
vectors, their commutation rules, and their behaviour under the action of
the braid group. 
\vskip .2 truecm
{\bf{\S \rrv. Real root vectors.}}
\vskip .1 truecm
Recall that the Dynkin diagram of $A_2^{(2)}$ is 
$$
A_2^{(2)}\
\ \ \compo\big<\,^{_=}_{^=}\!\!^{_=}_{^=}\!\!^{_=}_{^=}\compo$$
\vskip -1.1 truecm 
$${\phantom{
A_2^{(2)}\
\ \ }}\ \ _1\ \ \ \ \ \
\ \! _0{\phantom{a}}$$
and that we can describe its roots as follows: $\d=\a_0+2\a_1$ and
$$\Phi_+^{{\rm{re}}}=\{m\d+\a_1, (m+1)\d-\a_1,
(2m+1)\d\pm 2\a_1|m\in\mn\},$$ 
$$\Phi_+^{{\rm{im}}}=\{(m+1)\d|m\in\mn\};$$
moreover
$\omega_1=\omc_1=s_0s_1$: indeed $d_1=1$ and
$$s_0s_1(\a_1)=-s_0(\a_1)=-(\a_1+\a_0)=-\d+\a_1=\omega_1(\a_1).$$ 

It follows that the real root vectors are given by: 
$$E_{m\d+\a_1}=(T_0T_1)^{-m}(E_1),$$ 
$$E_{(m+1)\d-\a_1}=(T_0T_1)^{m}T_0(E_1),$$
$$E_{(2m+1)\d+2\a_1}=(T_0T_1)^{-m}T_1^{-1}(E_0),$$
$$E_{(2m+1)\d-2\a_1}=(T_0T_1)^{m}(E_0).$$

We shall now give the following proposition which will semplify many
calculations. 
\prop{\tgammat} 
$\forall\a\in\Phi_+^{{\rm{re}}}\setminus\{\a_0\}\ \ \ 
T_0\Xi(E_{\a})=E_{s_0(\a)}$.

$\forall\a\in\Phi_+^{{\rm{re}}}\setminus\{\a_1\}\ \ \
\Xi T_1(E_{\a})=E_{s_1(\a)}$.
\dim
$\Xi(E_{\a})\in\u_{q,\a}^+$, hence
$T_0\Xi(E_{\a})\in\u_{q,s_0(\a)}$; it is then enough to prove that
$T_0\Xi(E_{\a})$ is a root vector, which is easily seen: 
$T_0\Xi(T_0T_1)^{-m}(E_1)=(T_0T_1)^{m}T_0(E_1)$ and similarly for the
other root vectors, so that the first assertion is proved. 
From this result and from the fact that $(T_0\Xi)(\Xi T_1)=T_0T_1$
the second assertion of the proposition follows immediately.\finedim 
We shall now give some simple commutation 
relations between the real root vectors. 
\lem{\ai2a}
$[E_{\d-\a_1},F_1]=-[4]_qK_1E_0$ and $E_1E_{\d+\a_1}-q^{-2}E_{\d+\a_1}E_1= 
q^{-2}[4]_qE_{\d+2\a_1}$.
\dim
The first identity is a straightforward computation, using that
$$E_{\d-\a_1}=T_0(E_1)=-E_0E_1+q^{-4}E_1E_0;$$
the second relation is $\Xi T_1$ applied to the first one. \finedim
As a first application of lemma \ai2a we state the following proposition:
\prop{\ngen}
$\{E_1,E_{\d-\a_1},F_1,F_{\d-\a_1},K_1^{\pm 1},K_0^{\pm 1}\}$ is a system of
generators of $\u_q$ as a $\mc(q)$-algebra. 
\dim
Let ${\Cal{A}}$ be the $\mc(q)$-subalgebra of $\u_q$ generated by 
$$\{E_1,E_{\d-\a_1},F_1,F_{\d-\a_1},K_1^{\pm 1},K_0^{\pm 1}\}.$$
It is then enough to prove that $E_0,F_0\in{\Cal{A}}$. Since ${\Cal{A}}$ is
$\Omega$-stable it is enough to prove that $E_0\in{\Cal{A}}$. But 
$E_0=-\rapp{1}{[4]_q}K_1^{-1}[E_{\d-\a_1},F_1]\in{\Cal{A}}$.\finedim 
\cor{\stngen}
The least subalgebra of $\uq$ containing $\u_q^0$ and $E_1$ and stable
under the action of $T_0T_1$ and $T_0\Xi$ is $\uq$.\finedim

\vskip .2 truecm

{\bf{\S \edelta. The imaginary root vector $\te_{\d}$.}}

\vskip .1 truecm
Here we want to study some properties of $\te_{\d}$, where we recall that 
$$\te_{\d}=-E_{\d-\a_1}E_1+q^{-2}E_1E_{\d-\a_1}.$$ 
\prop{\tgd1}
$T_0\Xi(\te_{\d})=\te_{\d}$. 
\dim
The claim follows immediately from the definitions and thanks to
proposition
\tgammat. 
\finedim
\prop{\comcon}
The following commutation rules hold: 
$$[\te_{\d},E_1]=-[3]_q!E_{\d+\a_1},\ \ \ 
[\te_{\d},F_1]=-[3]_q!K_1E_{\d-\a_1}.$$
\dim
The proof consists in simple computations using the following ingredients: 
$\te_{\d}=E_0E_1^2-q^{-3}[2]_qE_1E_0E_1+q^{-6}E_1^2E_0$ (by the definition
and the fact that $E_{\d-\a_1}=-E_0E_1+q^{-4}E_1E_0$),
$$E_{\d+\a_1}=T_1^{-1}T_0^{-1}(E_1)=T_1^{-1}(-E_1E_0+q^{-4}E_0E_1)=
-K_1^{-1}[T_1^{-1}(E_0),F_1]=$$
$$=-K_1^{-1}\left[\sum_{r=0}^{4}(-1)^rq^{-r}E_1^{(r)}E_0E_1^{(4-r)},
F_1\right]=-\sum_{r=0}^{3}(-1)^rq^{-2r}E_1^{(r)}E_0E_1^{(3-r)}$$ 
and 
$$[E_1^m,F_1]=[m]_q\rapp{q^{1-m}K_1-q^{m-1}K_1^{-1}}{q-q^{-1}}E_1^{m-1}.$$
\finedim
\prop{\tod1}
$\Xi T_1(\te_{\d})=\te_{\d}$ (hence $T_0T_1(\te_{\d})=\te_{\d}$, thanks to
proposition \tgd1). 
\dim
Since $\Xi T_1(E_1)=-K_1^{-1}F_1$ we have that
$$\Xi
T_1(\te_{\d})=K_1^{-1}F_1E_{\d+\a_1}-q^{-2}E_{\d+\a_1}K_1^{-1}F_1=
-K_1^{-1}[E_{\d+\a_1},F_1]=$$
$$=\rapp{1}{[3]_q!}K_1^{-1}[[\te_{\d},E_1],F_1]=-\rapp{1}{[3]_q!}K_1^{-1}
[E_1,[[\te_{\d},F_1]]=K_1^{-1}[E_1,K_1E_{\d-\a_1}]=\te_{\d}.$$\finedim
\cor{\1coma}
Applying a power of $T_0T_1$ to proposition \comcon\ we get that $\forall
m\in\mn$ the following commutation relations hold: 
$$[\te_{\d},E_{m\d+\a_1}]=-[3]_q!E_{(m+1)\d+\a_1},\ \ \ 
[\te_{\d},E_{(m+1)\d-\a_1}]=[3]_q!E_{(m+2)\d-\a_1}.$$\finedim
We conclude this paragraph with the commutation rules between $\te_{\d}$ and the simple root
vectors $E_0$ and $F_0$. 
\lem{\dcoma} 
$[\te_{\d},F_0]=-(q^2-1)[3]_qE_1^2K_0^{-1}$ and 
$[\te_{\d},E_0]=-(q^2-1)[3]_qE_{\d-\a_1}^2.$
\dim 
The first relation is found by straightforward calculations; 
the second relation is found by applying $T_0\Xi$ to the first
one.\finedim
\vskip .2 truecm
{\bf{\S \emdelta. The imaginary root vectors $\te_{m\d}$.}}
\vskip .1 truecm
Here we generalize to the vectors $\te_{m\d}=-E_{m\d-\a_1}E_1+
q^{-2}E_1E_{m\d-\a_1}$ what we
know for
$\te_{\d}$. 

We start by proving the equivalence between some commutation rules among the
imaginary root vectors and their behaviour under the action of $T_0T_1$. 

This is completely similar to what happens in the case of $A_1^{(1)}$. 
\lem{\eqct}
$\forall m>0$ the following implications hold:
$$[\te_{\d},\te_{m\d}]=0\Leftrightarrow T_0T_1(\te_{(m+1)\d})=
\te_{(m+1)\d}$$ and 
$$T_0T_1(\te_{(m+1)\d})=
\te_{(m+1)\d}\Leftrightarrow T_0\Xi(\te_{(m+1)\d})=\te_{(m+1)\d}=\Xi
T_1(\te_{(m+1)\d}).$$
\dim
See \dasl2h.\finedim
The following proposition is a fundamental tool in the study of the
imaginary root vectors. 
\prop{\comdind} 
Let $P_m$ be the following statement: 

$i)_m$ $T_0T_1(\te_{m\d})=\te_{m\d}=T_0\Xi(\te_{m\d})=\Xi T_1(\te_{m\d})$; 

$ii)_m$ $E_0E_{(m-1)\d+\a_1}-q^{-4}E_{(m-1)\d+\a_1}E_0=$
$$=\rapp{1}{[4]_q}\big((q^4-q^2-q^{-2})E_{\d-\a_1}\te_{(m-1)\d}+
(q^2+q^{-2}-q^{-4})\te_{(m-1)\d}E_{\d-\a_1}+$$
$$+q^2E_{2\d-\a_1}\te_{(m-2)\d}-
q^{-2}\te_{(m-2)\d}E_{2\d-\a_1}\big);$$ 

$iii)_m$ $[\te_{m\d},E_1]=(q^4-q^{-2})\te_{(m-1)\d}E_{\d+\a_1}+ 
(q^2-q^{-4})E_{\d+\a_1}\te_{(m-1)\d}+$
$$+q^2\te_{(m-2)\d}E_{2\d+\a_1}-
q^{-2}E_{2\d+\a_1}\te_{(m-2)\d};$$ 

$iv)_m$
$[\te_{m\d},E_{\d-\a_1}]=-(q^4-q^{-2})E_{2\d-\a_1}\te_{(m-1)\d}-
(q^2-q^{-4})\te_{(m-1)\d}E_{2\d-\a_1}+$
$$-q^2E_{3\d-\a_1}\te_{(m-2)\d}+
q^{-2}\te_{(m-2)\d}E_{3\d-\a_1};$$

$v)_m$ $[\te_{r\d},\te_{s\d}]=0\ \ \ \forall r,s\leq m$. 

Then, if we set $\te_{-\d}=0$ and $\te_{0\d}=-\rapp{1}{q-q^{-1}}$,
$P_m$ is true for every $m>0$. 
\dim
The proof is an induction on $m$. 

If $m=1$ remark that: $i)_1$ is proposition \tgd1 and proposition \tod1; 
$ii)_1$ is the expression of $E_{\d-\a_1}$ in terms of $E_0$ and $E_1$; 
$iii)_1$ is proposition \comcon; $iv)_1$ is $T_0\Xi$ applied to
$iii)_1$; $v)_1$ is obvious. 

Let $m$ be bigger than 1. Then: 

$i)_m$ follows immediately from $v)_{m-1}$ (which in particular implies 
$[\te_{\d},\te_{(m-1)\d}]=0$) and from lemma \eqct. 

$ii)_m$ is obtained by the following steps, where we shall indicate among the lines the
results that we need to pass from one side to the other of the identities:
$$E_0E_{(m-1)\d+\a_1}-q^{-4}E_{(m-1)\d+\a_1}E_0=$$
{\centerline{\{corollary \1coma\}}}
$$=-\rapp{1}{[3]_q!}(E_0\te_{\d}E_{(m-2)\d+\a_1}-
E_0E_{(m-2)\d+\a_1}\te_{\d}+$$
$$-q^{-4}\te_{\d}E_{(m-2)\d+\a_1}E_0+
q^{-4}E_{(m-2)\d+\a_1}\te_{\d}E_0)=$$
{\centerline{\{lemma \dcoma\}}}
$$=-\rapp{1}{[3]_q!}\big([\te_{\d},E_0E_{(m-2)\d+\a_1}-
q^{-4}E_{(m-2)\d+\a_1}E_0]+$$
$$+(q^2-1)[3]_q(E_{\d-\a_1}^2E_{(m-2)\d+\a_1}-q^{-4}
E_{(m-2)\d+\a_1}E_{\d-\a_1}^2)\big)=$$
{\centerline{\{$i)_{m-1}$, $ii)_{m-1}$ and definition of
$\te_{(m-1)\d}$\}}}
$$=-\rapp{1}{[3]_q!}\Big(\rapp{1}{[4]_q}[\te_{\d},
(q^4-q^2-q^{-2})E_{\d-\a_1}\te_{(m-2)\d}+
(q^2+q^{-2}-q^{-4})\te_{(m-2)\d}E_{\d-\a_1}+$$
$$+q^2E_{2\d-\a_1}\te_{(m-3)\d}-
q^{-2}\te_{(m-3)\d}E_{2\d-\a_1}]+$$
$$-(q^2-1)[3]_qE_{\d-\a_1}\te_{(m-1)\d}-
(1-q^{-2})[3]_q\te_{(m-1)\d}E_{\d-\a_1}\Big)=$$
{\centerline{\{corollary \1coma and $v)_{m-1}$\}}}
$$=-\rapp{1}{[4]_q}((q^4-q^2-q^{-2})E_{2\d-\a_1}\te_{(m-2)\d}+
(q^2+q^{-2}-q^{-4})\te_{(m-2)\d}E_{2\d-\a_1}+$$
$$+q^2E_{3\d-\a_1}\te_{(m-3)\d}-
q^{-2}\te_{(m-3)\d}E_{3\d-\a_1}+$$
$$-(q^2-1)(q^2+q^{-2})E_{\d-\a_1}\te_{(m-1)\d}-
(1-q^{-2})(q^2+q^{-2})\te_{(m-1)\d}E_{\d-\a_1})=$$
{\centerline{\{reordering the summands\}}}
$$=\rapp{1}{[4]_q}((q^4-q^2-q^{-2})E_{\d-\a_1}\te_{(m-1)\d}+
(q^2+q^{-2}-q^{-4})\te_{(m-1)\d}E_{\d-\a_1}+$$
$$+q^2E_{2\d-\a_1}\te_{(m-2)\d}-
q^{-2}\te_{(m-2)\d}E_{2\d-\a_1})+$$
$$-\rapp{1}{[4]_q}([\te_{(m-1)\d},E_{\d-\a_1}]+
(q^4-q^{-2})E_{2\d-\a_1}\te_{(m-2)\d}+(q^2-q^{-4})\te_{(m-2)\d}E_{2\d-\a_1}+$$
$$+q^2E_{3\d-\a_1}\te_{(m-3)\d}-
q^{-2}\te_{(m-3)\d}E_{3\d-\a_1})=$$
{\centerline{\{$iv)_{m-1}$\}}}
$$=\rapp{1}{[4]_q}((q^4-q^2-q^{-2})E_{\d-\a_1}\te_{(m-1)\d}+
(q^2+q^{-2}-q^{-4})\te_{(m-1)\d}E_{\d-\a_1}+$$
$$+q^2E_{2\d-\a_1}\te_{(m-2)\d}-
q^{-2}\te_{(m-2)\d}E_{2\d-\a_1}),$$
which is the claim. 

$iii)_m$ is the result of the following computations, where the first identity is true thanks 
to $i)_m$:
$$[\te_{m\d},E_1]=
[-E_{(m-1)\d-\a_1}E_{\d+\a_1}+q^{-2}E_{\d+\a_1}E_{(m-1)\d-\a_1},E_1]=$$
$$=-E_{(m-1)\d-\a_1}E_{\d+\a_1}E_1+q^{-2}E_{\d+\a_1}E_{(m-1)\d-\a_1}E_1+$$
$$+E_1E_{(m-1)\d-\a_1}E_{\d+\a_1}-q^{-2}E_1E_{\d+\a_1}E_{(m-1)\d-\a_1}=$$ 
{\centerline{\{lemma \ai2a\ and definition of $\te_{(m-1)\d}$\}}}
$$=-q^2E_{(m-1)\d-\a_1}E_1E_{\d+\a_1}+[4]_qE_{(m-1)\d-\a_1}E_{\d+2\a_1}+
q^{-4}E_{\d+\a_1}E_1E_{(m-1)\d-\a_1}+$$
$$-q^{-2}E_{\d+\a_1}\te_{(m-1)\d}+
E_1E_{(m-1)\d-\a_1}E_{\d+\a_1}-q^{-2}E_1E_{\d+\a_1}E_{(m-1)\d-\a_1}=$$
{\centerline{\{definition of $\te_{(m-1)\d}$ and lemma \ai2a\}}}
$$=q^2\te_{(m-1)\d}E_{\d+\a_1}-q^{-2}E_{\d+\a_1}\te_{(m-1)\d}+$$
$$+[4]_q(E_{(m-1)\d-\a_1}E_{\d+2\a_1}-q^{-4}E_{\d+2\a_1}E_{(m-1)\d-\a_1})=$$
{\centerline{\{$\Xi T_1(ii)_m)$ and $i)_{m-1,m-2}$\}}}
$$=q^2\te_{(m-1)\d}E_{\d+\a_1}-q^{-2}E_{\d+\a_1}\te_{(m-1)\d}+
(q^4-q^2-q^{-2})\te_{(m-1)\d}E_{\d+\a_1}+$$
$$+(q^2+q^{-2}-q^{-4})E_{\d+\a_1}\te_{(m-1)\d}+q^2\te_{(m-2)\d}E_{2\d+\a_1}-
q^{-2}E_{2\d+\a_1}\te_{(m-2)\d}=$$
$$=(q^4-q^{-2})\te_{(m-1)\d}E_{\d+\a_1}+
(q^2-q^{-4})E_{\d+\a_1}\te_{(m-1)\d}+$$
$$+q^2\te_{(m-2)\d}E_{2\d+\a_1}-
q^{-2}E_{2\d+\a_1}\te_{(m-2)\d},$$
which is $iii)_m$.

$iv)_m$ is found by applying $T_0\Xi$ to $iii)_m$. 

$v)_m$: we have to prove that $[\te_{m\d},\te_{r\d}]=0$
$\forall r<m$. Since $[\te_{m\d},\te_{r\d}]\in\u_{q,(m+r)\d}^+$ it is enough
to prove that it commutes with $F_0$ and $F_1$ (see \liq): in particular,
proving that it belongs to the center of 
$\u_q$ will solve the problem; but the set of the elements
of $\uq$ commuting with $[\te_{m\d},\te_{r\d}]$ is a subalgebra of $\uq$
obviously containing $\u_q^0$ and stable by the action of $T_0T_1$ and
$T_0\Xi$ (because $[\te_{m\d},\te_{r\d}]$ is fixed by $T_0T_1$ and
$T_0\Xi$). Hence the problem reduces (see corollary \stngen) to show that
$[[\te_{m\d},\te_{r\d}],E_1]=0$, that is 
$[\te_{m\d},[\te_{r\d},E_1]]=[\te_{r\d},[\te_{m\d},E_1]]$: this will be
done by induction on $r$: 

$r=1$: proposition \comcon, $i)_m$, $(T_0T_1)^s(iii)_m)$ and $v)_{m-1}$
imply that  $$[\te_{m\d},[\te_{\d},E_1]]=-[3]_q![\te_{m\d},E_{\d+\a_1}]=$$
$$=-[3]_q!
\big((q^2-1)[3]_q\te_{(m-1)\d}E_{2\d+\a_1}+(1-q^{-2})[3]_qE_{2\d+\a_1}\te_{(m-1)\d}+$$
$$+q^2\te_{(m-2)\d}E_{3\d+\a_1}-q^{-2}E_{3\d+\a_1}\te_{(m-2)\d}\big)=$$
$$=[\te_{\d},(q^2-1)[3]_q\te_{(m-1)\d}E_{\d+\a_1}+(1-q^{-2})[3]_qE_{\d+\a_1}\te_{(m-1)\d}+$$
$$+q^2\te_{(m-2)\d}E_{2\d+\a_1}-q^{-2}E_{2\d+\a_1}\te_{(m-2)\d}]=$$
$$=[\te_{\d},[\te_{m\d},E_1]],$$ so that $[\te_{m\d},\te_{\d}]=0$. 

$r>1$: notice that if we put $x_1\ugdef q^4-q^{-2}$, $y_1\ugdef
q^2-q^{-4}$, $x_2\ugdef q^2$, $y_2\ugdef -q^{-2}$, 
we have that 
$\forall u\leq m$, $\forall v\in\mn$ 
(applying $(T_0T_1)^{-v}$ to $iii)_u$)
$$[\te_{u\d},E_{v\d+\a_1}]=
\sum_{s=1}^2(x_s\te_{(u-s)\d}E_{(v+s)\d+\a_1}+
y_sE_{(v+s)\d+\a_1}\te_{(u-s)\d});$$  hence, using the inductive
hypothesis and $v)_{m-1}$, we get that 
$$[\te_{m\d},[\te_{r\d},E_1]]=
\left[\te_{m\d},\sum_{s=1}^2(x_s\te_{(r-s)\d}E_{s\d+\a_1}+
y_sE_{s\d+\a_1}\te_{(r-s)\d})\right]=$$
$$=\sum_{s,t=1}^2(x_sx_t\te_{(r-s)\d}
\te_{(m-t)\d}E_{(s+t)\d+\a_1}+
x_sy_t\te_{(r-s)\d}E_{(s+t)\d+\a_1}\te_{(m-t)\d}+$$
$$+y_sx_t\te_{(m-t)\d}E_{(s+t)\d+\a_1}
\te_{(r-s)\d}+y_sy_tE_{(s+t)\d+\a_1}\te_{(m-t)\d}\te_{(r-s)\d})=$$
$$=\left[\te_{r\d},\sum_{t=1}^2x_t\te_{(m-t)\d}E_{t\d+\a_1}+
y_tE_{t\d+\a_1}\te_{(m-t)\d}\right]=[\te_{r\d},[\te_{m\d},E_1]],$$  which is the claim.\finedim
\vskip .2 truecm
{\bf{\S \cref. Commutation between positive and negative real
root vectors.}}
\vskip .1 truecm
The preceding proposition is first of all useful in proving other
commutation relations, now involving the real root vectors. 
\lem{\comefa}
$\forall r,s\in\mn$ we have that: 

$a)_{r,s}$ 
$[E_{r\d+\a_1},F_{s\d+\a_1}]=\cases
\rapp{K_{r\d+\a_1}-K_{r\d+\a_1}^{-1}}{q-q^{-1}}&{\rm{if}}\ r=s\cr 
-K_{s\d+\a_1}\te_{(r-s)\d}&{\rm{if}}\ r>s\cr
-\tf_{(s-r)\d}K_{r\d+\a_1}^{-1}&{\rm{if}}\ r<s;\endcases$

$b)_{r,s}$ 
$[E_{(r+1)\d-\a_1},F_{(s+1)\d-\a_1}]=\cases
\rapp{K_{(r+1)\d-\a_1}-K_{(r+1)\d-\a_1}^{-1}}{q-q^{-1}}&{\rm{if}}\ r=s\cr 
K_{(s+1)\d-\a_1}^{-1}\te_{(r-s)\d}&{\rm{if}}\ r>s\cr
\tf_{(s-r)\d}K_{(r+1)\d-\a_1}&{\rm{if}}\ r<s.\endcases$
\dim
$a)_{r,r}$ is proved by applying $(T_0T_1)^{-r}$ to the identity 
$[E_1,F_1]=\rapp{K_1-K_1^{-1}}{q-q^{-1}}$; if $r>s$ $a)_{r,s}$ follows
from the identity
$-E_{(r-s)\d-\a_1}E_1+q^{-2}E_1E_{(r-s)\d-\a_1}=\te_{(r-s)\d}$ again
applying to it $(T_0T_1)^{-r}$; if $r>s$ $a)_{r,s}$ is nothing but
$\Omega$ applied to $a)_{s,r}$. 

As for assertion $b)$, we have that $b)_{r,s}=\Xi
T_1(a)_{r+1,s+1})$.\finedim
\vskip .2 truecm
{\bf{\S \comed. The imaginary root vectors $E_{m\d}$.}}
\vskip .1 truecm
Recall that the subalgebra of $\u_q$ generated by the $\te_{m\d}$'s is a
commutative subalgebra pointwise fixed by $T_0T_1$ and $T_0\Xi$ (see
proposition
\comdind). Hence it makes sense to introduce new imaginary root vectors by
the following equation:
\ddefi{\nuoed}
Let us define some new elements in the subalgebra of $\u_q$ generated by
the $\te_{m\d}$'s by the following identity: 
$$1-(q-q^{-1})\sum_{m>0}\te_{m\d}u^m=
{\rm{exp}}\left((q-q^{-1})\sum_{m>0}E_{m\d}u^m\right).$$
Remark that of course $\forall m>0$ $E_{m\d}\in\u_{q,m\d}^+$
and 
$E_{m\d}+\te_{m\d}$ belongs to the subalgebra of $\u_q$ generated by 
$\{\te_{r\d}|r<m\}$. 
\prop{\rnorcomd}
$\forall m>0$ we have the following commutation relations: 
$$[E_{m\d},E_1]=
\rapp{[2m]_q}{m}(q^{2m}+(-1)^{m-1}+q^{-2m})E_{m\d+\a_1}.$$
\dim
The idea of this proof is completely similar to the proof of the analogous
result for $A_1^{(1)}$ (see \beck).

Let $\ts_m^{\pm},S_m^{\pm}:\u_q\rightarrow\u_q$ be the operators of left
(${}^+$) and right (${}^-$) multiplication respectively by $\te_{m\d}$ and
$E_{m\d}$ (where we keep on indicating by $\te_{0\d}$ and $\te_{-\d}$
respectively $-\rapp{1}{q-q^{-1}}$ and 0) and let us set 
$$\ts^{\pm}(u)\ugdef\sum_{m\geq 0}\ts_m^{\pm}u^m,\ \ \ 
S^{\pm}(u)\ugdef\sum_{m>0}S_m^{\pm}u^m;$$ moreover let us indicate by 
$T$ the automorphism $(T_0T_1)^{-1}$: remark that $\forall m$ 
$$\ts_m^{\pm}T=T\ts_m^{\pm}\ {\rm{and}}\ S_m^{\pm}T=TS_m^{\pm}.$$ 
Now observe that point $iii)$ of proposition \comdind\ can be written by
saying that $\forall m$ 
$$(\ts_m^+-\ts_m^-)(E_1)=$$
$$=\big((q^4-q^{-2})\ts_{m-1}^+T+ 
(q^2-q^{-4})\ts_{m-1}^-T+q^2\ts_{m-2}^+T^2-
q^{-2}\ts_{m-2}^-T^2\big)(E_1)$$ 
or, also, (multiplying by $u^m$ and summing over $m$)
$$\ts^+(u)(1+q^{-2}Tu)(1-q^4Tu)(E_1)=$$
$$=\ts^-(u)(1+q^2Tu)(1-q^{-4}Tu)(E_1);$$
multiplying both sides of this identity by $-(q-q^{-1})$ and considering
that $${\rm{log}}(-(q-q^{-1})\ts^{\pm}(u))=(q-q^{-1})S^{\pm}(u)$$ we get
that  $$(q-q^{-1})(S^+(u)-S^-(u))(E_1)=$$
$$=
({\rm{log}}(1+q^2Tu)+{\rm{log}}(1-q^{-4}Tu)-
{\rm{log}}(1+q^{-2}Tu)-{\rm{log}}(1-q^4Tu))(E_1),$$
that is 
$$[E_{m\d},E_1]=
\rapp{1}{m(q-q^{-1})}((-1)^{m-1}q^{2m}-q^{-4m}+(-1)^mq^{-2m}+q^{4m})
E_{m\d+\a_1}=$$
$$=\rapp{[2m]_q}{m}(q^{2m}+(-1)^{m-1}+q^{-2m})E_{m\d+\a_1}.$$\finedim
The last proposition generalizes immediately to the following more general
relations: 
\cor{\comda}
Let $m$ be bigger than 0; then $\forall r\in \mn$ we have:

$a)_r$ $[E_{m\d},E_{r\d+\a_1}]=
\rapp{[2m]_q}{m}(q^{2m}+(-1)^{m-1}+q^{-2m})E_{(m+r)\d+\a_1}$; 

$b)_r$ $[E_{m\d},E_{(r+1)\d-\a_1}]=
-\rapp{[2m]_q}{m}(q^{2m}+(-1)^{m-1}+q^{-2m})E_{(m+r+1)\d-\a_1}$; 

$c)_r$ $[E_{m\d},F_{r\d+\a_1}]=$
$$=\cases
\rapp{[2m]_q}{m}(q^{2m}+(-1)^{m-1}+q^{-2m})K_{r\d+\a_1}E_{(m-r)\d-\a_1}&
{\rm{if}}\ r<m\cr 
-\rapp{[2m]_q}{m}(q^{2m}+(-1)^{m-1}+q^{-2m})K_{m\d}F_{(r-m)\d+\a_1}&
{\rm{if}}\ r\geq m;\endcases$$ 

$d)_r$ $[E_{m\d},F_{(r+1)\d-\a_1}]=$
$$=\cases
-\rapp{[2m]_q}{m}(q^{2m}+(-1)^{m-1}+
q^{-2m})E_{(m-r-1)\d+\a_1}K_{(r+1)\d-\a_1}^{-1}&
{\rm{if}}\ r<m\cr 
\rapp{[2m]_q}{m}(q^{2m}+(-1)^{m-1}+q^{-2m})K_{m\d}^{-1}F_{(r-m+1)\d-\a_1}&
{\rm{if}}\ r\geq m.\endcases$$ 
\dim
$a)_r$ is $(T_0T_1)^{-r}$ applied to 
proposition \rnorcomd;
$b)_r$ is $\Xi T_1$ applied to $a)_{r+1}$; $c)_r$ is $(T_0T_1)^{-r}\Xi
T_1$ applied to $a)_0$; finally $d)_r$ is $\Xi T_1$ applied to 
$c)_{r+1}$.\finedim 
\cor{\comdd}
The imaginary root vectors form a ``Heisenberg'' algebra, with commutation
relations given by: 
$$[E_{r\d},F_{s\d}]=\d_{r,s}\rapp{[2r]_q}{r}(q^{2r}+(-1)^{r-1}+
q^{-2r})\rapp{K_{r\d}-K_{r\d}^{-1}}{q-q^{-1}}.$$
\dim
The claim is a straightforward consequence of proposition \comdind\ and
corollary
\comda; see also \cp.\finedim

\vskip .5 truecm

{\bf{\2na. $A_{2n}^{(2)}$ AND A COPY OF $A_2^{(2)}$.}}
\vskip .2 truecm

{\bf{\S \vf. Definition of $\varphi_1$.}}

\vskip .1 truecm

The aim of this section is to prove that the particular situation
of $\sl_3^{(2)}=A_2^{(2)}$ studied in the preceding section plays,
in some cases, the same role as $\sl_2^{(1)}$ does in general.
More precisely we know that if $i$ is a vertex of the finite
Dynkin diagram associated to an affine algebra $\hat{\g}^{(k)}$,
then, provided that $\hat{\g}^{(k)}$ is not of type $A_{2n}^{(2)}$
or that $i$ is not 1, there exists a homomorphism of
$\mc$-algebras
$\varphi_i:\u_q(\sl_2^{(1)})\rightarrow\u_q(\hat{\g}^{(k)})$ such
that $\varphi_i(E_1)=E_i$, $\varphi_i(E_0)=E_{\dk_i\d-\a_i}$,
$\varphi_i(q)=q_i$ and $\Omega\varphi_i=\varphi_i\Omega$; moreover
$\varphi_i$ is such that $T_i\varphi_i=\varphi_i T_1$ and
$T_{\lambda_i}\varphi_i=\varphi_i T_{\omc_1}(=\varphi_i T_0T_{\tau})$
(that is $T_1$ corresponds to $T_i$, $T_{\tau}=T_0T_{\tau}T_1^{-1}$ to
$T_{\lambda_i}T_i^{-1}$ and $T_0=T_0T_{\tau}^2$ to 
$T_{\lambda_i}^2T_i^{-1}$) (see proposition \eccofi).

Here we want to study the remaining cases, namely the situation
where $\hat{\g}^{(k)}$ is of type $A_{2n}^{(2)}$ and $i=1$: in this case
we construct a homomorphism of $\mc$-algebras
$\varphi_1:\u_q(\sl_3^{(2)})\rightarrow\u_q^{(1)}\subseteq\u_q(\hat{\g}^{(k)})$
such that $\varphi_1(q)=q_1(=q)$ ($\varphi_1$ is indeed a homomorphism
of $\mc(q)$-algebras), $\varphi_1(E_1)=E_1$ and
$\Omega\varphi_1=\varphi_1\Omega$; moreover $\varphi_1$ is such
that $T_1\varphi_1=\varphi_1 T_1$ and
$T_{\omega_1}\varphi_1=\varphi_1 T_{\omega_1}(=\varphi_1 T_0T_1)$
(that is $T_1$ corresponds to $T_1$ and $T_0$ to
$T_{\omega_1}T_1^{-1}$). Of course now the image of $E_0$ can't be
$E_{\d-\a_1}$ (since $\d-\a_1$ has the same length as $\a_1$ while
$\a_0$ is longer): in this case we will have that
$\varphi_1(E_0)=E_{\d-2\a_1}$.
\ddefi{\homoinj}Let
$\varphi_1:\u_q(A_2^{(2)})\rightarrow\u_q^{(1)}\subseteq\u_q(A_{2n}^{(2)})$
be the homomorphism of $\mc(q)$-algebras defined on the generators as
follows: $$\varphi_1(E_1)\ugdef E_1,\ \ \ \varphi_1(F_1)\ugdef
F_1,\ \ \ \varphi_1(K_1^{\pm 1})\ugdef K_1^{\pm 1},$$
$$\varphi_1(E_0)\ugdef E_{\d-2\a_1},\ \ \ \varphi_1(F_0)\ugdef
F_{\d-2\a_1},\ \ \ \varphi_1(K_0^{\pm 1})\ugdef K_{\d-2\a_1}^{\pm
1}.$$ Our program is to prove first of all that $\varphi_1$ is
well defined (of course if this is the case then it commutes with
$\Omega$), point to which we devote paragraph \bendef, after \S \combrw,
whose aim is recalling and describing some (mainly combinatorial)
properties of the root system and Weyl group of type $A_{2n}^{(2)}$; then,
in \S \tvft, we shall study the relations between the braid
group action and $\varphi_1$. These are the main tools which will allow us
to construct a PBW-basis.

\vskip .2 truecm

{\bf{\S \combrw. Root system and Weyl group: some properties.}}

\vskip .1 truecm

In this paragraph we shall study the properties of the root system
and the Weyl group of $A_{2n}^{(2)}$ which will enable us to
attack the central problem of this section, that is the research
of a copy of $\u_q(A_2^{(2)})$ in $\u_q(A_{2n}^{(2)})$. What we need
in particular is a more explicit description of the root vectors
(in particular of some root vectors) and their behaviour under the action
of
the braid group, goal
which can be achieved via a closer analysis of the Weyl group. 
More precisely our
attention will concentrate on the study of a reduced
expression of $\omega_1$ (remark that $\omega_1=\omc_1$, since $d_1=1$) and
of some of its properties, 
what will help us in manipulating the root 
vectors.
\rem{\lungho}
$l(\omega_1)=n(n+1)$. 
\dim
It is enough to compute the cardinality of $\Phi_+(\omega_1)$: 
$$\Phi_+(\omega_1)=\{\a\in\Phi_+|\omega_1^{-1}(\a)<0\}=
\left\{\d-\varepsilon\sum_{r=1}^k\a_r\big|1\leq
k\leq n,\varepsilon=1,2\right\}\cup$$
$$\cup
\left\{\varepsilon\d-2\sum_{r=1}^k\a_r-\sum_{r=k+1}^l\a_r\big|1\leq 
k<l\leq n,\varepsilon=1,2\right\},$$ 
so that 
$$\#\Phi_+(\omega_1)=2\left(n+\rapp{n(n-1)}{2}\right)=n(n+1),$$ 
which is the claim.\finedim 
The element of the Weyl group which we are now going to introduce
will play an important role in what follows. 
\ddefi{\vudexx}
In the following $w$ will denote the element of $W$ given by:
$$w\ugdef s_0s_ns_{n-1}\cdot...\cdot s_1.$$
\lem{\acvudexx}
The action of $w$ on $Q$ can be described as follows: 
$$w(\a_i)=\cases\d-2\sum_{r=1}^{n-1}\a_r&{\rm{if}}\
i=0\cr
-\d+\a_1+...+\a_n&{\rm{if}}\
i=1\cr
\d-(\a_2+...+\a_n)&{\rm{if}}\ i=2\cr
\a_{i-1}&{\rm{if}}\
3\leq i\leq n.
\endcases$$ 
Furthermore $\forall k\in\mz$
$$w^k(\a_i)=\cases
\left(1+2\left[\rapp{k}{n}\right]\right)
\d-2\sum_{r=1}^{\overline{n-k}}\a_r&
{\rm{if}}\ 
i=0\cr
\left[-\rapp{k}{n}\right]\d+\sum_{r=1}^{n+1-\overline{k}}\a_r
&{\rm{if}}\
i=1\cr
\a_{\overline{i-k}}&{\rm{if}}\
2\leq i\leq n\ {\rm{and}}\ \overline{i-k}\neq 1\cr
\d-(\a_2+...+\a_n)&{\rm{if}}\
2\leq i\leq n\ {\rm{and}}\ \overline{i-k}=1,
\endcases$$
where 
$\overline{}:\mz\rightarrow\{1,...,n\}$ is the application such that 
$\overline{k}\equiv k$ (mod $n$) $\forall k\in\mz$.
\dim 
It is clear that if $3\leq i\leq n$ then 
$$w(\a_i)=s_0s_n\cdot...\cdot s_1(\a_i)=s_0s_n\cdot...\cdot s_{i-1}(\a_i)=
s_0s_n\cdot...\cdot s_{i+1}(\a_{i-1})=\a_{i-1};$$ 
also, it is straightforward to see that 
$w(\a_1)=-\d+(a_1+...+\a_n)$,
$w(\a_2)=$ $=\d-(\a_2+...+\a_n)$
and $w(\a_0)=\a_0+2\a_n= 
\d-2\sum_{r=1}^{n-1}\a_r$, which proves the first assertion. 

Remark that $w^2(\a_2)=\a_n$, so that the second assertion immediately
follows for $2\leq i\leq n$. 

The fact that $w(\a_1+\a_2)=\a_1$ implies that $w(\a_1+...+\a_r)=
\a_1+...+\a_{r-1}$  $\forall r>1$, which easily leads to the remaining
assertions for $i=0,1$.\finedim
\cor{\tomwn}
$\omega_1=w^n$ and 
$(s_0s_n\cdot...\cdot s_1)^n$ is a reduced expression of 
$\omega_1$; in particular we have that
$T_{\omega_1}=T_w^n=(T_0T_n\cdot...\cdot T_1)^n$.
\dim
That $\omega_1=w^n$ follows from the fact that
$w^n(\a_i)=\a_i$ when $i\in I_0\setminus\{1\}$ and $w^n(\a_1)=-\d+\a_1$ 
(see lemma \acvudexx).  The remaining assertion is an immediate consequence
of lemma
\lungho. 
\finedim
\nota{\iot}
We fix $\iota$ so that $\iota_r+r\equiv 1$ (mod $n+1$) whenever
$1\leq r\leq  n(n+1)$ (this can be done thanks to corollary \tomwn; see
also remark \rutai). 
\cor{\eoeo}
$w^{n-1}(\a_0)=\d-2\a_1$; in particular 
$E_{\d-2\a_1}=T_{w^{n-1}}(E_0)$. 
\dim
That $w^{n-1}(\a_0)=\d-2\a_1$ is a direct consequence of lemma \acvudexx. 
Since $T_{w^{n-1}}(E_0)$ is a root vector, it must then be 
$E_{\d-2\a_1}$.\finedim
The lemma that we are now going to propose is devoted to the study of the
reduced expressions of some suitable elements of the Weyl group. Here lies
perhaps the main difference between this case and the general case that
has been recalled in section \becc: indeed in that case
$l(\omega_is_i\omega_i)=2l(\omega_i)-1$, and this implied that 
$(T_{\omega_i}T_i^{-1})^2$ acted as the identity on the
subalgebra of $\u_q$ generated by $\{E_i,E_{\d-\a_i}\}$ (see lemma
\basbec\ and remark \rutai); this is no more true here. However we can
still say something. 
\lem{\redx}
$l(s_0ws_1\omega_1)=l(s_0ws_1)+l(\omega_1)$, so that 
$T_{s_0ws_1w^k}=T_{s_0ws_1}T_{w}^k$ $\forall k\in\{0,...,n\}$; 

$l(s_0ws_1\omega_1)=l(ws_1)+l(s_1w^{-1}s_0ws_1\omega_1)$, that is 
$T_{s_0ws_1\omega_1}=T_{ws_1}T_{s_1w^{-1}s_0ws_1\omega_1}$; 
moreover $s_0ws_1w^{n-1}(\a_0)=\a_0$ and 
$s_1w^{-1}s_0ws_1\omega_1s_1(\a_1)=\a_1$. 
\dim
The first statement means that 
$\Phi_+(s_0ws_1)\subseteq\Phi_+(s_0ws_1\omega_1)$. 

Now $(s_0ws_1)^{-1}(\Phi_+(s_0ws_1))\subseteq\sum_{i=2}^n\mz\a_i$ which
is pointwise fixed by $\omega_1$, so that 
$$\a\in\Phi_+(s_0ws_1)\Rightarrow(s_0ws_1\omega_1)^{-1}(\a)=
(s_0ws_1)^{-1}(\a)<0\Rightarrow\a\in\Phi_+(s_0ws_1\omega_1).$$

The second part of the lemma  is equivalent to saying that 
$\Phi_+(ws_1)\subseteq\Phi_+(s_0ws_1\omega_1)$. Of course
$$\Phi_+(ws_1)=\{s_0s_n\cdot...\cdot s_{r+1}(\a_r)|2\leq r\leq
n+1\}=\{\a_0+\a_n+...+\a_r|2\leq r\leq n+1\}$$ so that, if $2\leq r\leq
n+1$, 
$(s_0ws_1\omega_1)^{-1}(\a_0+\a_n+...+\a_r)=$
$$=\cases\omega_1^{-1}(\a_0+2\a_n+...+2\a_{r+1}+\a_r+...+\a_2)
&{\rm{if}}\ r\leq n\cr
\omega_1^{-1}(\a_0+2\a_n+...+2\a_2)&{\rm{if}}\ r=n+1\endcases$$
$$=\cases-\d-(\a_r+...+\a_2+2\a_1)<0
&{\rm{if}}\ r\leq n\cr
-\d-2\a_1<0&{\rm{if}}\ r=n+1.\endcases$$
Now we have that 
$$s_0ws_1w^{n-1}(\a_0)=s_0ws_1(\d-2\a_1)=s_0w(\d+2\a_1)=
s_0(2\d-(w(\d-2\a_1)))=$$
$$=s_0(2\d-\omega_1(\a_0))=
s_0\omega_1(\d+2(\a_1+...+\a_n))=
s_0(-\a_0)=\a_0.$$
Finally we get that 
$$s_1w^{-1}s_0ws_1\omega_1s_1(\a_1)=s_1w^{-1}s_0w(\d+\a_1)=$$
$$=s_1w^{-1}(\a_0+\a_1+...+\a_n)=s_1(-\a_1)=\a_1,$$
which concludes the proof.\finedim
\cor{\othusg}
$E_{\d-2\a_1}=T_{s_0ws_1}^{-1}(E_0)$ and 
$E_{\d-\a_1}=T_{s_0ws_1}^{-1}T_{ws_1}(E_1)$. 
\dim 
Straightforward consequence of lemma \redx.\finedim

\vskip .2 truecm

{\bf{\S \bendef. $\varphi_1$ is well defined.}}

\vskip .1 truecm

Here we want to prove that $\varphi_1$ is well defined, which means
that the relations between the ``canonical'' generators of
$\u_q(A_2^{(2)})$  (that is $E_i$,
$F_i$,
$K_i^{-1}$ for $i=0,1$) are preserved by $\varphi_1$ (so that
$\varphi_1:\uq(A_2^{(2)})\rightarrow\uq(A_{2n}^{(2)})$ is well defined),
and that $E_{\d-2\a_1}\in\u_q^{(1)}$ (so that $\varphi_1$ effectively maps
$\uq(A_2^{(2)})$ in $\u_q^{(1)}$).

For some of these relations this result is immediate as indicated
in proposition \banale. The others will require some more work.

Let us recall which are the relations that we are speaking about:
\ddefi{\relsl32} $\u_q(A_2^{(2)})$ is the $\mc(q)$-algebra
generated by $\{E_i,F_i,K_i^{\pm 1}|i=0,1\}$ with the following
relations:

$$K_1K_0=K_0K_1;\leqno[KK]$$

$$K_1E_1=q^2E_1K_1,\ K_1E_0=q^{-4}E_0K_1,\
K_0E_1=q^{-4}E_1K_0,\ K_0E_0=q^8E_0K_0;\leqno{[KE]} $$

$$K_1F_1=q^{-2}F_1K_1,\ K_1F_0=q^{4}F_0K_1,\ 
K_0F_1=q^{4}F_1K_0,\ K_0F_0=q^{-8}F_0K_0;\leqno{[KF]} $$

$$[E_1,F_1]=\rapp{K_1-K_1^{-1}}{q-q^{-1}},\ \ \
[E_0,F_0]=\rapp{K_0-K_0^{-1}}{q^4-q^{-4}},\ \ \
[E_1,F_0]=0=[E_0,F_1];\leqno{[EF]} $$

$$\sum_{r=0}^2(-1)^rE_0^{(r)}E_1E_0^{(2-r)}=0=
\sum_{r=0}^5(-1)^rE_1^{(r)}E_0E_1^{(5-r)};\leqno{[EE]} $$

$$\sum_{r=0}^2(-1)^rF_0^{(r)}F_1F_0^{(2-r)}=0=
\sum_{r=0}^5(-1)^rF_1^{(r)}F_0F_1^{(5-r)}.\leqno{[FF]} $$

\prop{\banale} The relations [KK], [KE], [KF] and [EF] are
obviously preserved by $\varphi_1$. 
\dim 
The subalgebra of $\u_q(A_{2n}^{(2)})$ generated by
$\{E_i|i\neq 1\}$ (to which $E_{\d-2\a_1}$ belongs) commutes with
$F_1$, so that $[E_{\d-2\a_1},F_1]=0$ (and of course
$[E_1,F_{\d-2\a_1}]=$
$=\Omega([E_{\d-2\a_1},F_1])=0$). The other relations
are trivial.\finedim

The point is now to prove the remaining relations.

\rem{\stabom} Once we have proved the compatibility of relations
[EE] with the definition of $\varphi_1$, it is enough to apply
$\Omega$ in order to find the same result for [FF]: indeed both the left
and the right hand sides of [FF] are found by applying $\Omega$ to the
corresponding expressions of [EE].

\lem{\m1m21} In $\u_q(A_{2n}^{(2)})$ we have that
$E_{\d-\a_1}E_{\d-2\a_1}=q^{-4}E_{\d-2\a_1}E_{\d-\a_1}$.

\dim 
Since $E_{\d-2\a_1}=T_{w^{n-1}}(E_0)$ and $
E_{\d-\a_1}=T_{\omega_1}T_1^{-1}(E_1)=T_{w^n}T_1^{-1}(E_1)$
we can apply Levendorskii-Soibelman formula (see \lswg): the claim follows
from the fact that 
$(\d-2\a_1|\d-\a_1)=4$, and that the only possibility for a sum of
positive roots of the form $m\d-\a$ (with $\a$ a positive multiple of a
root of the finite system associated to $A_{2n}^{(2)}$) to be equal to
$2\d-3\a_1$ is that these positive roots are $\d-2\a_1$ and
$\d-\a_1$.\finedim

\lem{\tconcor}$-E_0T_nT_{n-1}\cdot...\cdot
T_3(E_2)+q^{-4}T_nT_{n-1}\cdot...\cdot T_3(E_2)E_0=
T_0T_nT_{n-1}\cdot...\cdot T_3(E_2)$.

\dim Since $T_nT_{n-1}\cdot...\cdot T_3(E_2)\in\u_{q,\a_n+...+\a_2}^+$ we
can write it as $
\sum xE_ny$ where
$x,y\in(E_2,...,E_{n-1})$; in particular $x$ and $y$ commute with
$E_0$ and are fixed points for $T_0$. Then we have
$$-E_0T_nT_{n-1}\cdot...\cdot
T_3(E_2)+q^{-4}T_nT_{n-1}\cdot...\cdot T_3(E_2)E_0= \sum
x(-E_0E_n+q^{-4}E_nE_0)y=$$

$$=\sum xT_0(E_n)y=T_0\left(\sum xE_ny\right)=
T_0T_nT_{n-1}\cdot...\cdot T_3(E_2).$$\finedim

\cor{\am21} In $\u_q(A_{2n}^{(2)})$ we have that
$-E_{\d-2\a_1}E_1+q^{-4}E_1E_{\d-2\a_1}=E_{\d-\a_1}$.

\dim The proof consists of a computation in which lemma \tconcor\ and
corollary \othusg\ play a fundamental role: if $n\geq 2$
$$-E_{\d-2\a_1}E_1+q^{-4}E_1E_{\d-2\a_1}=
-T_{s_0ws_1}^{-1}(E_0)E_1+q^{-4}E_1T_{s_0ws_1}^{-1}(E_0)=$$
$$=T_{s_0ws_1}^{-1}(-E_0T_{s_0ws_1}(E_1)+q^{-4}T_{s_0ws_1}(E_1)E_0)=$$
$$=T_{s_0ws_1}^{-1}(-E_0T_{s_0ws_1s_2}(-E_2E_1+q^{-2}E_1E_2)+
q^{-4}T_{s_0ws_1s_2}(-E_2E_1+q^{-2}E_1E_2)E_0)=$$
$$=-T_{s_0ws_1}^{-1}((-E_0T_n\cdot...\cdot
T_3(E_2)+q^{-4}T_n\cdot...\cdot T_3(E_2)E_0)E_1+$$ 
$$-
q^{-2}E_1(-E_0T_n\cdot...\cdot
T_3(E_2)+q^{-4}T_n\cdot...\cdot T_3(E_2)E_0))=$$
$$=-T_{s_0ws_1}^{-1}(T_0T_nT_{n-1}\cdot...\cdot T_3(E_2)E_1-
q^{-2}E_1T_0T_nT_{n-1}\cdot...\cdot T_3(E_2))=$$
$$=T_{s_0ws_1}^{-1}T_0T_nT_{n-1}\cdot...\cdot T_3(-E_2E_1+q^{-2}E_1E_2)=
T_{s_0ws_1}^{-1}T_{ws_1}(E_1)=E_{\d-\a_1},$$
where we have used that $[E_1,E_0]=0$.\finedim

\prop{\releed}The following relation holds in $\u_q(A_{2n}^{(2)})$:
$$\sum_{r=0}^2(-1)^r{2\brack
r}_{q^4}E_{\d-2\a_1}^rE_1E_{\d-2\a_1}^{2-r}=0.$$

\dim The proof is a simple application of lemma \m1m21 and
corollary \am21. 
\finedim
So we are now left to prove just the last relation [EE]. To this
aim we start by recalling a very simple and basic combinatorial
lemma.

\lem{\combkkp1} Let $a,b$ be elements of a $\mc(q)$-algebra
satisfying the relation

$$\sum_{r=0}^k(-1)^rq^{pr}{k\brack r}_qa^rba^{k-r}=0$$ for some
natural number $k$ and some integer $p$. Then $a,b$ also satisfy
the relation $$\sum_{r=0}^{k+1}(-1)^rq^{(p+1)r}{k+1\brack
r}_qa^rba^{k+1-r}=0.$$

\dim The proof is a straightforward computation based on the
identity $$q^{(p+1)r}{k+1\brack r}_q= q^{pr}{k\brack 
r}_q+q^{k+p+1}q^{p(r-1)}{k\brack r-1}_q.$$\finedim

Before going to the point, we study some other simple
relations.

\lem{\comem2f} The following is true in $\u_q(A_{2n}^{(2)})$:
$$[E_{\d-2\a_1},F_i]=0\ \ \ \forall i\in\{3,4,...,n\};$$
$$T_0^{-1}(E_{\d-2\a_1})\in\u_q^+;$$
$$[T_{w^{n-2}}(E_0),E_1]=0=[T_{w^{n-2}}(E_0),F_2].$$

\dim The proof consists in some easy remarks:

1) if $i\in\{3,4...,n\}$ then $T_{s_0ws_1}(F_i)=T_nT_{n-1}\cdot...\cdot
T_{i-1}(F_i)=F_{i-1}$, which commutes with
$E_0$, so that $$[E_{\d-2\a_1},F_i]=T_{s_0ws_1}^{-1}([E_0,F_{i-1}])=0;$$ 
2)
$T_0^{-1}(E_{\d-2\a_1})=T_0^{-1}T_{w^{n-1}}(E_0)=T_{s_0w^{n-1}}(E_0)$
which is obviously in $\u_q^+$;

3) $T_{w^{n-2}}(E_0)\in(E_0,E_3,E_4,...,E_n)$, so that it
commutes with $E_1$ and $F_2$.\finedim
The following proposition is what we still need for our goal; it
concerns the last relation [EE], and the proof is adapted from the
analogous situation in case $A_2^{(1)}$ (see \beck).

\prop{\relee5} In $\u_q(A_{2n}^{(2)})$ the following relation is
verified:
$$\sum_{r=0}^5(-1)^rE_1^{(r)}E_{\d-2\a_1}E_1^{(5-r)}=0.$$

\dim Let
$\E\ugdef\sum_{r=0}^5(-1)^rE_1^{(r)}E_{\d-2\a_1}E_1^{(5-r)}$; we
want to prove that $\E=0$. The strategy to achieve this goal is to
prove that $T_i^{-1}(\E)\in\u_q^+$ $\forall i=0,... n$ (see \liq). Recall
that
$[\E,F_i]=0$ implies $T_i^{-1}(\E)\in\u_q^+$.

Lemma \comem2f evidently implies that $[\E,F_i]=0$ when
$i=3,...,n$ and $T_0^{-1}(\E)\in\u_q^+$ (since
$T_0^{-1}(E_1)\in\u_q^+$).

Moreover, thanks to proposition \banale, $[\E,F_1]=$
$$=\sum_{r=1}^5(-1)^r\rapp{q^{1-r}K_1-q^{r-1}K_1^{-1}}
{q-q^{-1}}E_1^{(r-1)}E_{\d-2\a_1}E_1^{(5-r)}+$$
$$+\sum_{r=0}^4(-1)^r
E_1^{(r)}E_{\d-2\a_1}\rapp{q^{r-4}K_1-q^{4-r}K_1^{-1}}
{q-q^{-1}}E_1^{(4-r)}=$$
$$=\sum_{r=0}^4(-1)^r\rapp{-q^{-r}K_1+q^rK_1^{-1}+
q^{-r}K_1-q^rK_1^{-1}}{q-q^{-1}}E_1^{(r)}E_{\d-2\a_1}E_1^{(4-r)}=0.$$ So
now we have just to show that $[\E,F_2]=0$. To this aim we rewrite $\E$
(when $n>1$, which is the case we are interested in)  noticing that
$$E_{\d-2\a_1}=T_{w^{n-1}}(E_0)=T_{w^{n-2}}T_0T_n(E_0)=
T_{w^{n-2}}T_n^{-1}(E_0)=$$
$$=T_{w^{n-2}}(E_0E_n^{(2)}-q^{-2}E_nE_0E_n+q^{-4}E_n^{(2)}E_0)=$$
$$=T_{w^{n-2}}(E_0)E_2^{(2)}-q^{-2}E_2T_{w^{n-2}}(E_0)E_2+
q^{-4}E_2^{(2)}T_{w^{n-2}}(E_0),$$ so that 
$\E$ is given by the expression
$$\sum_{r=0}^5(-1)^rE_1^{(r)}
(T_{w^{n-2}}(E_0)E_2^{(2)}- q^{-2}E_2T_{w^{n-2}}(E_0)E_2+
q^{-4}E_2^{(2)}T_{w^{n-2}}(E_0))E_1^{(5-r)}.$$ Then we are ready
to conclude thanks to the last part of lemma \comem2f:
$$[\E,F_2]=$$
$$=\sum_{r=0}^5(-1)^rE_1^{(r)}\Big(T_{w^{n-2}}(E_0)
\rapp{q^{-2}K_2-q^2K_2^{-1}}{q^2-q^{-2}}E_2-q^{-2}\rapp{K_2-K_2^{-1}}{q^2-q^{-2}}
T_{w^{n-2}}(E_0)E_2+$$
$$-q^{-2}E_2T_{w^{n-2}}(E_0)\rapp{K_2-K_2^{-1}}{q^2-q^{-2}}+
q^{-4}\rapp{q^{-2}K_2-q^2K_2^{-1}}{q^2-q^{-2}}
E_2T_{w^{n-2}}(E_0)\Big)E_1^{(5-r)}=$$
$$=T_{w^{n-2}}(E_0)\sum_{r=0}^5(-1)^r
X E_1^{(r)}E_2E_1^{(5-r)}-\sum_{r=0}^5(-1)^r
E_1^{(r)}E_2E_1^{(5-r)}YT_{w^{n-2}}(E_0)$$
where 
$$X={q^{2r-2}K_2-q^{2-2r}K_2^{-1}-q^{2r-6}K_2+q^{2-2r}K_2^{-1}\over
q^2-q^{-2}}=q^{2r-4}K_2$$
and 
$$Y=\rapp{q^{2r-8}K_2-q^{4-2r}K_2^{-1}-q^{2r-12}K_2+
q^{4-2r}K_2^{-1}}{q^2-q^{-2}}=q^{2r-10}K_2,$$
so that
$$[\E,F_2]=q^{-4}T_{w^{n-2}}(E_0)K_2\left(\sum_{r=0}^5(-1)^rq^{2r}
E_1^{(r)}E_2E_1^{(5-r)}\right)+$$
$$-q^{-10}\left(\sum_{r=0}^5(-1)^rq^{2r}
E_1^{(r)}E_2E_1^{(5-r)}\right)K_2T_{w^{n-2}}(E_0)$$ which is zero thanks to
lemma \combkkp1 and to the relation, valid in
$\u_q(A_{2n}^{(2)})$ (if $n>1$),
$$\sum_{r=0}^3(-1)^rE_1^{(r)}E_2E_1^{(3-r)}=0.$$ This completely proves the
claim, that is $\E=0$.\finedim

The paragraph can now be closed with the
next proposition, which states the achievement of our goal:
\prop{\welldef} $\varphi_1$ is well defined.
\dim
The only thing to prove is that $E_{\d-2\a_1}\in\u_q^{(1)}$. This follows
immediately from proposition \banale\ and corollary \am21: indeed by
straightforward computations 
$$[E_{\d-\a_1},F_1]=-[4]_qK_1E_{\d-2\a_1}.$$\finedim

\vskip .2 truecm

{\bf{\S \tvft. Relations between $\varphi_1$ and the braid group
action.}}

\vskip .1 truecm

The results of paragraph \bendef\ are the first step in our attempt
to prove that the relations between the root vectors in
$\u_q(A_2^{(2)})$ found in section \adt\ are still valid in
$\u_q(A_{2n}^{(2)})$. To conclude this piece of program we need to
analyze how the braid group acts with respect to $\varphi_1$. As
we already announced, we want to prove that
$T_1\varphi_1=\varphi_1T_1$ and
$T_{\omega_1}T_1^{-1}\varphi_1=\varphi_1 T_0$. We start from this
second statement.

\prop{\fvto} The action of $T_0$ on $\u_q(A_2^{(2)})$ corresponds,
under $\varphi_1$, to the action of $T_{\omega_1}T_1^{-1}$, that
is $T_{\omega_1}T_1^{-1}\varphi_1=\varphi_1 T_0$.

\dim Of course it is enough to prove that
$T_{\omega_1}T_1^{-1}\varphi_1=\varphi_1 T_0$ on the generators of
$\u_q(A_2^{(2)})$, that is on $E_1$, $E_0$, $F_1$, $F_0$, $K_1$,
$K_0$. For the elements of $\u_q^0(A_2^{(2)})$ the claim is
trivial. 

Remark that $T_{\omega_1}T_1^{-1}$, $\varphi_1$ and $T_0$
all commute with $\Omega$, so that we are reduced to prove that
$T_{\omega_1}T_1^{-1}\varphi_1(E_1)=\varphi_1 T_0(E_1)$ and
$T_{\omega_1}T_1^{-1}\varphi_1(E_0)=\varphi_1 T_0(E_0)$; this
means that we want to prove on one hand that
$$\varphi_1(-E_0E_1+q^{-4}E_1E_0)=E_{\d-\a_1}$$ (which is nothing
but corollary \am21), and on the other hand that
$$-F_{\d-2\a_1}K_{\d-2\a_1}=T_{\omega_1}T_1^{-1}(E_{\d-2\a_1}),$$
or equivalently (see corollaries \tomwn, \eoeo\ and \othusg) that
$-F_0K_0=T_{ws_1}T_{s_0ws_1}^{-1}(E_0)$; but of course
$T_{ws_1}T_{s_0ws_1}^{-1}(E_0)=T_0(E_0)$, which proves the
claim.\finedim

\rem{\thpnt} We concentrate now on the question of proving that
$T_1\varphi_1=\varphi_1T_1$, or equivalently that
$T_1^{-1}\varphi_1=\varphi_1T_1^{-1}$. The same considerations as
before tell us that this problem reduces to the question whether
$T_1^{-1}\varphi_1(E_1)=\varphi_1T_1^{-1}(E_1)$ and
$T_1^{-1}\varphi_1(E_0)=\varphi_1T_1^{-1}(E_0)$: notice that the
first relation is obviously true while the second one can be
translated into
$$T_1^{-1}(E_{\d-2\a_1})=
\sum_{r=0}^4(-1)^rq^{-r}E_1^{(r)}E_{\d-2\a_1}E_1^{(4-r)}.$$ 
To face this problem let us set the following definition: 
\ddefi{\access}
We introduce the elements $X_0$,...,$X_n$ of $\u_q(A_{2n}^{(2)})$
as follows: $$X_n\ugdef E_0,\ \ \ X_{s-1}\ugdef
T_s^{-1}(X_s)\ \ \forall s=1,...,n.$$
Remark that $X_1=E_{\d-2\a_1}$ (corollary \othusg) and
$X_0=T_1^{-1}(E_{\d-2\a_1})$. 

Our strategy will divide in two steps: 

1) proving that $X_{s-1}=
\sum_{r=0}^2(-1)^rq^{-2r}E_s^{(r)}X_sE_s^{(2-r)}$ $\forall s=2,...,n$; 

2) computing $X_0$ thanks to the description of $X_1$ given in
1): this will be the claim. 
\lem{\sfzt}
Let $s$ be in $\{2,...,n-1\}$ and let $X\in\u_q$ be an element
of $\u_q(A_{2n}^{(2)})$ commuting with $E_s$ and fixed by $T_s$. 
Then 
$$\sum_{r=0}^2(-1)^rq^{-2r}E_s^{(r)}E_{s+1}XE_{s+1}E_s^{(2-r)}=
T_s^{-1}(E_{s+1}XE_{s+1}).$$
\dim 
Using the Serre relation between $E_s$ and $E_{s+1}$ we have that
$$\sum_{r=0}^2(-1)^rq^{-2r}E_s^{(r)}E_{s+1}XE_{s+1}E_s^{(2-r)}=$$
$$=E_{s+1}E_sXE_{s+1}E_s-\rapp{1}{q^2+q^{-2}}E_{s+1}E_sXE_sE_{s+1}-
q^{-2}E_sE_{s+1}XE_{s+1}E_s+$$
$$+q^{-4}E_sE_{s+1}XE_sE_{s+1}-
\rapp{q^{-4}}{q^2+q^{-2}}E_{s+1}E_sXE_sE_{s+1}=$$
$$=(-E_{s+1}E_s+q^{-2}E_sE_{s+1})X(-E_{s+1}E_s+q^{-2}E_sE_{s+1})=$$
$$=
T_s^{-1}(E_{s+1})XT_s^{-1}(E_{s+1})=T_s^{-1}(E_{s+1}XE_{s+1}).$$\finedim
\cor{\susf}Let $s$ and $X$ be as in lemma \sfzt, and let 
$Y\ugdef\sum_{r=0}^2(-1)^rq^{-2r}E_{s+1}^{(r)}XE_{s+1}^{(2-r)}$; then 
$$\sum_{r=0}^2(-1)^rq^{-2r}E_s^{(r)}YE_s^{(2-r)}=
T_s^{-1}(Y).$$
\dim The claim is an immediate application of lemma \sfzt.\finedim
\cor{\xsind}
Let $s$ be in $\{2,...,n\}$; then 
$X_{s-1}=\sum_{r=0}^2(-1)^rq^{-2r}E_s^{(r)}X_sE_s^{(2-r)}$. 

In particular
$E_{\d-2\a_1}=\sum_{r=0}^2(-1)^rq^{-2r}E_2^{(r)}X_2E_2^{(2-r)}$. 
\dim 
The statement is proved by induction on $s$: if $s=n$ the 
claim is obvious (it is the definition of $T_n^{-1}(E_0)$). 

If $2\leq s<n$ define $X$ by $X\ugdef X_{s+1}$; then
$X\in(E_0,E_n,...,E_{s+2})$, hence it commutes with
$E_s$  and is fixed by $T_s$. Moreover, by the
inductive hypothesis, the element
$Y$ defined in corollary \susf\ is nothing but $X_s$. Hence we have that 
$$X_{s-1}=T_s^{-1}(X_s)=T_s^{-1}(Y)=
\sum_{r=0}^2(-1)^rq^{-2r}E_s^{(r)}YE_s^{(2-r)}=$$
$$=\sum_{r=0}^2(-1)^rq^{-2r}E_s^{(r)}X_sE_s^{(2-r)}.$$\finedim
We have thus solved the first step. For the second one let us develop some
computations. 
\lem{\sfilza} Let $X$ be an element of $\u_q$ commuting with $E_1$.
Then we have the following identities:
$$E_2XE_2E_1^{(4)}=$$
$$=E_2E_1^{(2)}XE_2E_1^{(2)}-
\rapp{[2]_q}{[3]_q}E_2E_1^{(2)}XE_1E_2E_1+
\rapp{1}{q^2+q^{-2}}E_2E_1^{(2)}XE_1^{(2)}E_2;$$
$$E_1E_2XE_2E_1^{(3)}=$$
$$=E_1E_2E_1XE_2E_1^{(2)}-
\rapp{1}{[2]_q}E_1E_2E_1XE_1E_2E_1+
\rapp{1}{[3]_q}E_1E_2E_1XE_1^{(2)}E_2;$$
$$E_1^{(3)}E_2XE_2E_1=$$
$$=E_1^{(2)}E_2XE_1E_2E_1-
\rapp{1}{[2]_q}E_1E_2E_1XE_1E_2E_1+
\rapp{1}{[3]_q}E_2E_1^{(2)}XE_1E_2E_1;$$
$$E_1^{(4)}E_2XE_2=$$
$$=E_1^{(2)}E_2XE_1^{(2)}E_2-
\rapp{[2]_q}{[3]_q}E_1E_2E_1XE_1^{(2)}E_2+
\rapp{1}{q^2+q^{-2}}E_2E_1^{(2)}XE_1^{(2)}E_2.$$
\dim It is a matter of manipulation of the left hand sides using the
relation
$$\sum_{r=0}^3(-1)^rE_1^{(r)}E_2E_1^{(3-r)}=0.$$
\cor{\concsflz} Let $X$ be an element of $\u_q$ commuting with $E_1$ and
fixed by $T_1$. Then
$$\sum_{r=0}^4(-1)^rq^{-r}E_1^{(r)}E_2XE_2E_1^{(4-r)}=
T_1^{-1}(E_2XE_2).$$
\dim The claim is a consequence of lemma \sfilza\ once
that one remarks that
$$T_1^{-1}(E_2)=E_2E_1^{(2)}-q^{-1}E_1E_2E_1+q^{-2}E_1^{(2)}E_2;$$
indeed
$$\sum_{r=0}^4(-1)^rq^{-r}E_1^{(r)}E_2XE_2E_1^{(4-r)}=$$
$$=E_2E_1^{(2)}XE_2E_1^{(2)}-
\rapp{[2]_q}{[3]_q}E_2E_1^{(2)}XE_1E_2E_1+
\rapp{1}{q^2+q^{-2}}E_2E_1^{(2)}XE_1^{(2)}E_2+$$
$$-q^{-1}E_1E_2E_1XE_2E_1^{(2)}+
\rapp{q^{-1}}{[2]_q}E_1E_2E_1XE_1E_2E_1-
\rapp{q^{-1}}{[3]_q}E_1E_2E_1XE_1^{(2)}E_2+$$
$$+q^{-2}E_1^{(2)}E_2XE_2E_1^{(2)}+$$
$$-q^{-3}E_1^{(2)}E_2XE_1E_2E_1+
\rapp{q^{-3}}{[2]_q}E_1E_2E_1XE_1E_2E_1-
\rapp{q^{-3}}{[3]_q}E_2E_1^{(2)}XE_1E_2E_1+$$
$$+q^{-4}E_1^{(2)}E_2XE_1^{(2)}E_2-
\rapp{q^{-4}[2]_q}{[3]_q}E_1E_2E_1XE_1^{(2)}E_2+
\rapp{q^{-4}}{q^2+q^{-2}}E_2E_1^{(2)}XE_1^{(2)}E_2=$$
$$=E_2E_1^{(2)}XE_2E_1^{(2)}-q^{-1}E_2E_1^{(2)}XE_1E_2E_1+
q^{-2}E_2E_1^{(2)}XE_1^{(2)}E_2+$$
$$-q^{-1}E_1E_2E_1XE_2E_1^{(2)}+
q^{-2}E_1E_2E_1XE_1E_2E_1-q^{-3}E_1E_2E_1XE_1^{(2)}E_2+$$
$$+q^{-2}E_1^{(2)}E_2XE_2E_1^{(2)}
-q^{-3}E_1^{(2)}E_2XE_1E_2E_1
+q^{-4}E_1^{(2)}E_2XE_1^{(2)}E_2=$$
$$=(E_2E_1^{(2)}-q^{-1}E_1E_2E_1+
q^{-2}E_1^{(2)}E_2)X(E_2E_1^{(2)}-q^{-1}E_1E_2E_1+
q^{-2}E_1^{(2)}E_2)=$$
$$=T_1^{-1}(E_2XE_2).$$\finedim
We can now conclude: 
\prop{\finpar} $\varphi_1$ commutes with $T_1$.
\dim
As already remarked, we want to prove that 
$$T_1^{-1}(E_{\d-2\a_1})=
\sum_{r=0}^4(-1)^rq^{-r}E_1^{(r)}E_{\d-2\a_1}E_1^{(4-r)};$$ 
following the notations of definition \access\ and the result of corollary
\xsind\
$$E_{\d-2\a_1}=\sum_{r=0}^2(-1)^rq^{-2r}E_2^{(r)}X_2E_2^{(2-r)}\ \
{\rm{with}}\ \ 
[X_2,E_1]=0\ \
{\rm{and}}\ \ T_1(X_2)=X_2;$$ then corollary \concsflz\
implies the claim.\finedim  
To sum up, we can state the
following proposition, which makes the importance of the present section
explicit: indeed proposition \rotug\ allows us to read in $A_{2n}^{(2)}$
all the relations found in section \adt\ for $A_2^{(2)}$; at the same
time it allows to define $E_{(m\d,1)}$ $\forall m>1$. 
\prop{\rotug} Let $\nu$ be the group homomorphism
$\nu:Q(A_2^{(2)})\rightarrow Q(A_{(2n)}^{(2)})$ defined in the obvious way
by 
$\nu(\a_1)\ugdef\a_1$, $\nu(\d)\ugdef\d$ and let the injection
$\tnu:\Phi_+(A_2^{(2)})\rightarrow\tPhi_+(A_{(2n)}^{(2)})$ be given by: 
$$\tnu(\a)\ugdef\cases\nu(\a)&{\rm{if\
}}\a{\rm{\ is\ real}}\cr(\nu(\a),1)&{\rm{if\ }}\a{\rm{\ is\
imaginary;}}\endcases$$
Then $\varphi_1$ has the property that $$\varphi_1(E_{\a})=E_{\tnu(\a)}\ \
\forall\a\in\Phi_+(A_2^{(2)}).$$  Moreover
$\varphi_1(\te_{\a})=\te_{\tnu(\a)}$
$\ \ \forall\a\in\Phi_+^{{\rm{im}}}(A_2^{(2)})$.\finedim 
\dim
The claim for $\te_{\a}$ when $\a\in\Phi_+^{{\rm{im}}}(A_2^{(2)})$
follows from the assertion on the real roots. 
Then $E_{(m\d,1)}\in\uq(A_{2n}^{(2)})$ can be defined $\forall m>0$ (see
the analogous definition, corollary \firv) and the claim for $E_{m\d}$
is also abvious. 

For the (real) roots of the form $m\d-\a_1$ the claim follows from
propositions \fvto\ and \finpar, and from remark \rutai; 
on the other hand, since 
$E_{\d-2\a_1}\in\u_q^{(1)}$, $T_{\omc_i}(E_{\d-2\a_1})=E_{\d-2\a_1}$ 
$\forall i\in I_0\setminus\{1\}$; thus
$$
E_{(2m+1)\d-2\a_1}=
(T_{\omc_1}\cdot...\cdot T_{\omc_n})^m T_w^{n-1}(E_0)=$$
$$=(T_{\omc_1}\cdot...\cdot T_{\omc_n})^m(E_{\d-2\a_1})=
T_{\omc_1}^m(E_{\d-2\a_1})=\varphi_1(T_0T_1)^m(E_0)$$
and analogously
$$
E_{(2m+1)\d+2\a_1}=
(T_{\omc_1}\cdot...\cdot T_{\omc_n})^{-m-1}
T_w^{n-1}T_0(E_0)=$$
$$=T_{\omc_1}^{-m-1}(-F_{\d-2\a_1}K_{\d-2\a_1})=
\varphi_1(T_0T_1)^{-m-1}T_0(E_0),$$
which concludes the proof.\finedim

\vskip .5 truecm 

{\bf{\ij. RELATIONS BETWEEN $\u_q^{(i)}$ AND $\u_q^{(j)}$ ($i\neq j$).}}
\vskip .2 truecm 
In this section we shall complete the description of the
commutation relations among the root vectors. The basis for this
section is, as for section \becc, an application of Beck's work (with
something more to check in the twisted case: see \S \prtclr); the argument
works also in the case $A_{2n}^{(2)}$.

In particular the aim of this part of the work is to understand the
commutation rules between the vectors from $\u_q^{(i)}$ and
$\u_q^{(j)}$ when $i\neq j$.
\rem{\trivcomm}
$\forall i,j\in I_0$ such that $a_{ij}=0$ we have 
$[\u_{q}^{(i)},\u_{q}^{(j)}]=0$.
\dim
$[E_i,\Omega^r(E_{\dk_j\d-\a_j})]=
T_{\lambda_j}([E_i,\Omega^rT_j^{-1}(E_j)])=0$.\finedim
Hence we are reduced to the case when $a_{ij}<0$; recall also that, if this
is the case, then either $a_{ij}=-1$ or $a_{ji}=-1$.
\vskip .2 truecm
{\bf{\S \gnrl. First level of commutation.}}
\vskip .1 truecm
In this paragraph we sketch Beck's recursive approach to the
problem. 
\lem{\bancomm}
Let $i\neq j\in I_0$. Then $\forall r,s\in\mn$ the following relations
hold:
$$[E_{r\dk_i\d+\a_i},F_{s\dk_j\d+\a_j}]=0\ \ \ {\rm{and}}\ \ \  
E_{r\dk_i\d+\a_i}E_{(s+1)\dk_j\d-\a_j}= 
q_i^{a_{ij}}E_{(s+1)\dk_j\d-\a_j}E_{r\dk_i\d+\a_i}.$$
\dim
It is enough to apply $T_{\lambda_i}^{-r}T_{\lambda_j}^{-s}$
and $T_{\lambda_i}^{-r}T_{\lambda_j}^{s+1}$
to the relation $[E_i,F_j]=0$. \finedim
\lem{\fststp}
Let $i,j\in I_0$ be such that $a_{ij}=-1$ (then ${\dk_i\over
\dk_j}\in\mn$).  If $r\geq 1$, we have: 
$$[\te_{(r\dk_i\d,i)},E_j]=\cases
-E_{\dk_i\d+\a_j}&{\rm{if}}\ \ r=1\cr
T_{\lambda_j}^{-{\dk_i\over \dk_j}}
(q_iE_j\te_{((r-1)\dk_i\d,i)}-q_i^{-1}\te_{((r-1)\dk_i\d,i)}E_j)
&{\rm{if}}\ \ r>1 \endcases$$
and 
$$[\te_{(r\dk_i\d,i)},F_j]=\cases
-K_jE_{\dk_i\d-\a_j}&{\rm{if}}\ \ r=1\cr
K_{\dk_i\d}T_{\lambda_j}^{{\dk_i\over \dk_j}}
(q_i^{-1}F_j\te_{((r-1)\dk_i\d,i)}-q_i\te_{((r-1)\dk_i\d,i)}F_j)
&{\rm{if}}\ \ r>1;\endcases$$
similarly, if $r\geq{\dk_i\over \dk_j}$, we get: 
$$[\te_{(r\dk_j\d,j)},E_i]=\cases
[a_{ji}]_{q_j}E_{\dk_i\d+\a_i}&{\rm{if}}\ \ r={\dk_i\over \dk_j}\cr
T_{\lambda_i}^{-1}
(q_j^{-a_{ji}}E_i\te_{((r\dk_j-\dk_i)\d,j)}-q_j^{a_{ji}}
\te_{((r\dk_j-\dk_i)\d,j)}E_i)
&{\rm{if}}\ \ r>{\dk_i\over \dk_j}\endcases$$
and 
$$[\te_{(r\dk_j\d,j)},F_i]=\cases
[a_{ji}]_{q_j}K_iE_{\dk_i\d-\a_i}&{\rm{if}}\ \ r={\dk_i\over \dk_j}\cr
K_{\dk_i\d}T_{\lambda_i}
(q_j^{a_{ji}}F_i\te_{((r\dk_j-\dk_i)\d,j)}\!-q_j^{-a_{ji}}
\te_{((r\dk_j-\dk_i)\d,j)}F_i)
&{\rm{if}}\ \ r>{\dk_i\over \dk_j}.\endcases$$
\dim
The result on which this lemma is based on is that 
$T_{\lambda_i}T_i^{-1}(E_j)=T_{\lambda_j}^{{\dk_i\over \dk_j}}T_i(E_j)$
and 
$T_{\lambda_i}T_i^{-1}(F_j)=T_{\lambda_j}^{{\dk_i\over \dk_j}}T_i(F_j)$
(see \beck\ and point 4) of lemma \basbec): indeed notice that if
$a_{ij}=-1$ then $(X_{\tn}^{(k)},i)\neq (A_{2n}^{(2)},1)$. 
The fact that ${\dk_i\over \dk_j}\in\mn$ 
is a consequence of the following
simple remark: we have either $\dk_i=1$ $\forall i\in I_0$ or $\dk_i=d_i$
$\forall i\in I_0$, so that ${\dk_i\over \dk_j}$ is either $1$ or
${d_i\over d_j}={a_{ji}\over a_{ij}}=-a_{ji}$.\finedim
\vskip .2 truecm
{\bf{\S \prtclr. Some particular computations: low twisted cases.}}
\vskip .1 truecm
What happens of $[\te_{(r\dk_j\d,j)},E_i]$ and $[\te_{(r\dk_j\d,j)},F_i]$
when
$1\leq r<{\dk_i\over \dk_j}$ (and $i,j\in I_0$ are such that $a_{ij}=-1$)
is the matter of the present paragraph; of course this situation occurs
neither for the non twisted algebras, nor for
$A_{2n}^{(2)}$, since in these cases we have $\tk=1$. 
Hence we have now to figure out what happens in cases
$A_{2n-1}^{(2)}$, $D_{n+1}^{(2)}$, $E_6^{(2)}$ and 
$D_{4}^{(3)}$ (when $1\leq r<{\dk_i\over\dk_j}$). 

The aim of this paragraph is to prove that the 
brackets above are zero. 

The problem will be first reduced to compute the commutation of $F_j$ with
a suitable element lying in $\u_q^+$; this simple remark will
immediately solve the case
$D_{n+1}^{(2)}$ (see proposition \dn2).

In the other cases the question will be further modified in looking for an
element of the Weyl group with ``good'' properties (see remark \sndstp\
and lemma \linw); the last part of the paragraph will be devoted to a
closer and more detailed investigation of $A_{2n-1}^{(2)}$, 
$D_{4}^{(3)}$ and $E_6^{(2)}$. 
\rem{\wtlf}
Let $i,j\in I_0$ be such that $a_{ij}=-1\neq a_{ji}$ and let $r$ be such
that
$1\leq r<{\dk_i\over \dk_j}$ (notice that $\dk_j=d_j=1$, $\dk_i=d_i=k>1$
and 
$\lambda_j=\omega_j=\omc_j$); we
want to prove that
$[\te_{(r\d,j)},E_i]=0$ and
$[\te_{(r\d,j)},F_i]=0$. 

Since, by a simple manipulation depending on $a_{ij}=-1$ (which
implies that $-E_jE_i+q_i^{-1}E_iE_j=T_i^{-1}(E_j)$), and on the
definition of $\te_{(r\d,i)}$, we see that 
$$[\te_{(r\d,j)},E_i]=
-T_i^{-1}T_j^{-1}([T_jT_iT_{\omega_j}^rT_j^{-1}(E_j),F_j]K_j)$$ 
and 
$$[\te_{(r\d,j)},F_i]=
-q_iK_{\a_i-\a_j}T_j^{-1}
([T_jT_iT_{\omega_j}^rT_j^{-1}(E_j),F_j]),$$ the
problem reduces to prove that
$[T_jT_iT_{\omega_j}^rT_j^{-1}(E_j),F_j]=0$. 

Remark that from 
$l(s_js_i\omega_j^r)=l(\omega_j^r)+2$ we get  
$T_jT_iT_{\omega_j}^rT_j^{-1}(E_j)\in\u^+_{q,r\d-s_j(\a_i+\a_j)}$.
\prop{\dn2}
Let $X_{\tn}^{(k)}=D_{n+1}^{(2)}$ and let $i,j\in I_0$ be such that
$a_{ij}=-1\neq a_{ji}$ (this means that $i=2$, $j=1$).

Then we have that $[\te_{(\d,j)},E_i]=0$ and
$[\te_{(\d,j)},F_i]=0$. 
\dim
Thanks to remark \wtlf\ we need to show that
$[T_1T_2T_{\omega_1}T_1^{-1}(E_1),F_1]=0$, where 
$T_1T_2T_{\omega_1}T_1^{-1}(E_1)\in\u^+_{q,\d-s_1(\a_1+\a_2)}$.
But, since $\d=\sum_{r=0}^n\a_r$ and $a_{12}=-2$,
$$\d-s_1(\a_1+\a_2)=\d-(\a_1+\a_2)=\a_0+\sum_{r=3}^n\a_r,$$ so that 
$T_1T_2T_{\omega_1}T_1^{-1}(E_1)$ belongs to the subalgebra of
$\u_q^+$ generated by $\{E_r|r\neq 1\}$, which immediately implies
that
it commutes with $F_1$ (since $E_r$ does for $r\neq 1$).\finedim
\rem{\sndstp}
Remark that for an element $x\in\u_q^+$ the property 
$[x,F_j]=0$ is equivalent to the condition $T_j^{\pm 1}(x)\in\u_q^+$ (see
\liq); since for any $i,j\in I_0$ and $r>0$ we obviously have
that 
$T_j^{-1}T_jT_iT_{\omega_j}^rT_j^{-1}(E_j)=
T_iT_{\omega_j}^rT_j^{-1}(E_j)\in\u_q^+$, remark \wtlf\ allows us to
translate the problem of understanding the behaviour of
$[\te_{(r\d,j)},E_i]$ and
$[\te_{(r\d,j)},F_i]$ (when $1\leq r<\dk_i$ and $i,j\in I_0$ are such that 
$a_{ij}=-1\neq a_{ji}$) into the following
question: is
$T_j^2T_iT_{\omega_j}^rT_j^{-1}(E_j)$ an element of $\u_q^+$?

The strategy that we shall use to prove that the answer is ``yes'' is to
find an element $w$ of the Weyl group $\tW$ and an element $t\in I$ with
the property that $l(s_jws_t)=l(w)+2$ and that 
$T_j^2T_iT_{\omega_j}^rT_j^{-1}(E_j)=T_{s_jw}(E_t)$. 

The next lemma illustrates more precisely what we are looking for. 
\lem{\linw}
Fix $i\neq j\in I_0$ and $r>0$, and
let $w\in\tW$ satisfy the following
conditions: 

1) $l(w^{-1}s_js_i\omega_j^rs_j)=l(s_js_i\omega_j^rs_j)-l(w)$;

2) $w^{-1}s_js_i\omega_j^rs_j(\a_j)$ is a simple root;

3) $w^{-1}(\a_j)>0$.

Then $T_j^2T_iT_{\omega_j}^rT_j^{-1}(E_j)\in\u_q^+$.
\dim
The first hypothesis means that
$T_jT_iT_{\omega_j}^rT_j^{-1}=T_wT_{w^{-1}s_js_i\omega_j^rs_j}$ and the
third means that $T_jT_w=T_{s_jw}$. 

Then we have:
$T_j^2T_iT_{\omega_j}^rT_j^{-1}(E_j)=
T_{s_jw}T_{w^{-1}s_js_i\omega_j^rs_j}(E_j)=
T_{s_jw}(E_t)$, where $t$ is such that 
$w^{-1}s_js_i\omega_j^rs_j(\a_j)=\a_t$;
hence the only thing to prove is that $s_jw(\a_t)>0$, which is
obvious since
$s_jw(\a_t)=s_i\omega_j^rs_j(\a_j)=r\d-(\a_i+\a_j)$.\finedim
The aim of the next lemmas is now to find an element $w$ with the
properties stated in lemma \linw: in order to exhibit 
such an element we need to understand something about $\omega_j$.
\ddefi{\dadnp1}
Consider the case of $A_{2n-1}^{(2)}$. Let $\tau$ be the
Dynkin diagram automorphism defined by
$$\tau(t)\ugdef\cases t&{\rm{if}}\ 1\leq t\leq n-1\cr
0&{\rm{if}}\ t=n\cr
n&{\rm{if}}\ t=0\endcases$$
and define elements $y_r\in\tW$ (for $r=1,...,n$) as follows:
$y_r\ugdef\tau s_ns_{n-1}\cdot...\cdot s_r$.
\lem{\ladnp1}
In $A_{2n-1}^{(2)}$ the following assertions hold: 

a) if $1\leq s\leq r\leq n-1$ and $r<i\leq n$
or
$i=0$, then 
$$y_{s+1}\cdot...\cdot y_r(\a_i)=\cases\a_{i-r+s}&{\rm{if}}\ i\neq 0\cr 
\a_0+\sum_{t=1}^{r-s}(\a_{n-t}+\a_{n-t+1})&{\rm{if}}\ i=0;\endcases$$

b) $\forall r=2,...,n$ $\Phi_+(y_r)=\{\tau(\a_i+...+\a_n)|r\leq i\leq
n\}$; 

c) if $1\leq r<i\leq n$ we have
$$y_2\cdot...\cdot y_r(\tau(\a_i+...+\a_n))=
\d-(\a_1+2(\a_2+...+\a_{i-r})+\a_{i-r+1}+...+\a_{n-r+1});$$

d) if $1\leq r<i\leq n$ we have
$$s_2\omega_2^{-1}s_1s_2y_2...y_r(\tau(\a_i+...+\a_n))=$$
$$=\cases
\a_2&{\rm{if}}\ i-r=1=n-r\cr
-\d-(\a_1+...+\a_{n-r+1})&{\rm{if}}\ i-r=1<n-r\cr
-(\a_3+...+\a_{n-r+1})&{\rm{if}}\ i-r=2\cr
-\d-(\a_1+2(\a_2+...+\a_{i-r})+\a_{i-r+1}+...+\a_{n-r+1})&{\rm{if}}\
i-r>2.\endcases$$
\dim 
Statement a) is easily proved by induction on $r-s$. 

For b) one 
immediately shows that $y_r^{-1}(\tau(\a_i+...+\a_n))=s_r\cdot...\cdot
s_n(\a_i+...+\a_n)<0$ if
$r\leq i\leq n$, which means that $\{\tau(\a_i+...+\a_n)|r\leq i\leq
n\}\subseteq\Phi_+(y_r)$; but $\#\Phi_+(y_r)=l(y_r)\leq n-r+1$, which
proves the claim. 

Assertion c) follows from a) remarking that 
$\tau(\a_i+...+\a_n)=\a_i+...+\a_{n-1}+\a_0$ and 
$y_2\cdot...\cdot y_r(\a_0)=\d-(\a_1+2(\a_1+...+\a_{n-r})+\a_{n-r+1})$. 

d) is just 
$s_2\omega_2^{-1}s_1s_2$ applied to c).\finedim
\prop{\adnp1}
Let $X_{\tn}^{(k)}=A_{2n-1}^{(2)}$ and let $i,j\in I_0$ be such that
$a_{ij}=-1\neq a_{ji}$ (this means that $i=1$, $j=2$). 

Then we have: 

1) $l((y_2\cdot...\cdot y_{n-1})^{-1}s_2s_1\omega_2s_2)=
l(s_2s_1\omega_2s_2)-l(y_2\cdot...\cdot y_{n-1})$; 

2) $(y_2\cdot...\cdot y_{n-1})^{-1}s_2s_1\omega_2s_2(\a_2)=\a_0$;

3) $(y_2\cdot...\cdot y_{n-1})^{-1}(\a_2)>0$. 

In particular $T_2^2T_1T_{\omega_2}T_2^{-1}(E_2)\in\u_q^+$, so that 
$[\te_{(\d,2)},E_1]=0$ and
$[\te_{(\d,2)},F_1]=0$.
\dim
To prove 1) it is enough to show that
$\Phi_+(y_r)\subseteq\Phi_+((y_2\cdot...\cdot
y_{r-1})^{-1}s_2s_1\omega_2s_2)$
$\forall r=2,...,n-1$, and this is a straightforward consequence of points
b) and d) of lemma \ladnp1.

From point a) of lemma \ladnp1\ we see also that 
$y_2\cdot...\cdot y_{n-1}(\a_0)=\d-(\a_1+\a_2)=$
$=s_2s_1\omega_2s_2(\a_2)$,
which is 2), and $y_2\cdot...\cdot y_{n-1}(\a_n)=\a_2$, which implies
3).\finedim
\prop{\dqt}
Consider the case $D_4^{(3)}$ and let $i,j\in I_0$ be such that
$a_{ij}=-1\neq a_{ji}$: this means that $i=2$ and $j=1$. Then 
$[\te_{(r\d,1)},E_2]=0$ and
$[\te_{(r\d,1)},F_2]=0$ for $r=1,2$.
\dim
First of all remark that 
$\d=\a_0+2\a_1+\a_2$, hence
$s_1s_2\omega_1s_1(\a_1)=\a_0$, so that 
lemma \linw\ applies with $w=$id; it follows from remark \sndstp\ that 
$[\te_{(\d,1)},E_2]=0$ and
$[\te_{(\d,1)},F_2]=0$.

Consider now the element $s_0s_1s_2$; then 

1) $l((s_0s_1s_2)^{-1}s_1s_2\omega_1^2s_1)=l(s_1s_2\omega_1^2s_1)-3$:
indeed
$$s_1\omega_1^{-2}s_2s_1(\a_0),\ \
s_1\omega_1^{-2}s_2s_1(\a_0+\a_1)\
\ {\rm{and}}\ \ s_1\omega_1^{-2}s_2s_1(3\a_0+3\a_1+\a_2)$$ are negative
roots (they are respectively $-\d+\a_1$, $-3\d-(\a_1+\a_2)$ and
$-3\d-\a_2$); 

2) $(s_0s_1s_2)^{-1}s_1s_2\omega_1^2s_1(\a_1)=\a_1$ which is a simple
root; 

3) $(s_0s_1s_2)^{-1}(\a_1)=\a_0>0$. 

As above, remark \sndstp\ and lemma \linw\ apply to the present
situation and imply that $[\te_{(2\d,1)},E_2]=0$ and
$[\te_{(2\d,1)},F_2]=0$.\finedim
\prop{\esix}
Consider the case $E_6^{(2)}$ and let $i,j\in I_0$ be such that
$a_{ij}=-1\neq a_{ji}$: this means that $i=3$, $j=2$. 
Then 
$[\te_{(\d,2)},E_3]=0$ and
$[\te_{(\d,2)},F_3]=0$.
\dim 
Consider the element $w\ugdef
s_0s_1s_2s_3s_4s_2s_3s_2s_1s_0s_2s_3s_2s_1s_2s_3$; then, through
direct computations very similar (even though longer!) to those used for
$D_4^{(3)}$ (see the proof of proposition \dqt, 1), 2), 3)), one
sees that:

1) $l(w^{-1}s_2s_3\omega_2s_2)=l(s_2s_3\omega_2s_2)-l(w)$; 

2) $s_2s_3\omega_2s_2(\a_2)=w(\a_2)$;

3) $w^{-1}(\a_2)(=\a_0)>0$. 

The claim follows.\finedim
\vskip .2 truecm
{\bf{\S \comconcl. General commutation formulas.}}
\vskip .1 truecm
We are now ready to discuss the commutation relations, as announced at the
beginning of this section. Before giving them explicitely, let us
introduce some definitions and notations. 
\ddefi{\simlac}
In the twisted case different from $A_{2n}^{(2)}$ let us define a matrix
$A^s=(a_{ij}^s)_{ij\in I_0}$ (called the simply laced matrix associated to
$(I,A)$) as follows:
$$a_{ij}^s\ugdef{\rm{max}}\{a_{ij},-1\};$$
this can be said equivalently by:
$$\cases a_{ii}^s=2&{}\cr
a_{ij}^s=0&{\rm{if}}\ a_{ij}=0\cr
a_{ij}^s=-1&{\rm{if}}\ a_{ij}<0.\endcases$$
Remark that $(I_0,A^s)$ is a Kac-Moody datum of type $A_n$.

Define also $\cnt_{ij}\ugdef\cases 1&{\rm{if}}\ a_{ij}\geq-1\cr
d_j&{\rm{otherwise}}.\endcases$
\teo{\cascomm}
Let $i,j\in I_0$ and $r,s>0$ be such that $\dk_i|r$ and $\dk_j|s$; then we
have: 
$$[E_{(r\d,i)},E_{(s-\dk_j)\d+\a_j}]=
\cases x_{ijr}E_{(s-\dk_j+r)\d+\a_j}&{\rm{if}}\ \dk_j|r\cr
0&{\rm{otherwise}};\endcases\leqno (1)$$
$$[E_{(r\d,i)},E_{s\d-\a_j}]=
\cases -x_{ijr}E_{(s+r)\d-\a_j}&{\rm{if}}\ \dk_j|r\cr
0&{\rm{otherwise}};\endcases\leqno (2)$$
$$[E_{(r\d,i)},F_{(s-\dk_j)\d+\a_j}]=
\cases x_{ijr}K_{(s-\dk_j)\d+\a_j}E_{(r-s+\dk_j)\d-\a_j}&{\rm{if}}\
\dk_j|r\ {\rm{and}}\ r\geq s\cr 
-x_{ijr}K_{r\d}F_{(s-r-\dk_j)\d+\a_j}&{\rm{if}}\ \dk_j|r\
{\rm{and}}\ r<s\cr
0&{\rm{otherwise}};\endcases\leqno (3)$$
$$[E_{(r\d,i)},F_{s\d-\a_j}]=
\cases -x_{ijr}E_{(r-s)\d+\a_j}K_{\a_j-s\d}&{\rm{if}}\ \dk_j|r\
{\rm{and}}\  r\geq s\cr
x_{ijr}F_{(s-r)\d-\a_j}K_{r\d}^{-1}&{\rm{if}}\ \dk_j|r\ {\rm{and}}\ 
r<s\cr
0&{\rm{otherwise}};\endcases\leqno (4)$$
$$[E_{(r\d,i)},E_{(s\d,j)}]=0;\leqno (5)$$
$$[E_{(r\d,i)},F_{(s\d,j)}]=\d_{rs}x_{ijr}{K_{r\d}-K_{r\d}^{-1}\over
q_j-q_j^{-1}};\leqno
(6)$$ 
where $$x_{ijr}=
\cases(o(i)o(j))^r\rapp{[ra_{ij}]_{q_i}}{r}
&{\rm{in\ the\ non\ twisted\ case\ and
}}\cr &{\rm{\ 
\ in\ case\
}}A_{2n}^{(2)}\  {\rm{with}}\ (i,j)\neq(1,1)\cr
\rapp{[2r]_q}{r}(q^{2r}+(-1)^{r-1}+
q^{-2r})&{\rm{in\ case\
}}A_{2n}^{(2)}\  {\rm{with}}\ i=j=1\cr
(o(i)o(j))^{r\over
\dk_i\cnt_{ij}}\rapp{\dk_i\cnt_{ij}[ra_{ij}^s]_q}{r[d_i]_q}
&{\rm{otherwise}}\endcases$$ 
and $o:I_0\rightarrow\{\pm 1\}$ is such that $a_{ij}<0\Rightarrow
o(i)o(j)=-1$.
\dim
For $i=j$, the claim is nothing but proposition \commi\ and $\varphi_1$
applied to \comda\ and \comdd. 

For $a_{ij}=0$, the claim follows from remark \trivcomm.

For $a_{ij}<0$, (1) and (2) are found as an application of lemma \fststp\ 
and of the results of paragraph \prtclr\ by a standard argument: see
\beck\ and proposition \rnorcomd; (3) and (4) are found
applying a suitable power of $T_{\lambda_j}$ respectively to (2) and (1);
(5) follows immediately from (1) and (2) while (6) follows from (3)
and (4) (using proposition \commi,2)). 
\finedim

\vskip .5 truecm

{\bf{\secpbw. A PBW-BASIS.}} 
\vskip .2 truecm 

\vskip .1 truecm
An important and useful consequence of the commutation relations given
in the preceding section is the construction of a PBW-basis for the
general quantum algebra of affine type.  The aim of this section is thus to
exhibit such a basis using a convex ordering of the root system, to recall
its principal properties (such as the Levendorskii-Soibelman formula)
and to introduce some notations. 
\vskip .2 truecm

{\bf{\S \acrd. A convex ordering.}}

\vskip .1 truecm

\ddefi{\totord}
We denote by $\preceq$ the total ordering of $\tPhi_+$ defined by: 
$$\forall r,s\in\mz:\ \ \ \beta_r\preceq\beta_s\Leftrightarrow\cases
s\leq r\leq 0&{\rm{or}}\cr
r\leq 0<s&{\rm{or}}\cr
1\leq s\leq r,
\endcases$$
$$\forall r\in\mz,\forall\a\in\tPhi_+^{{\rm{im}}}:\ \ \ 
\beta_r\preceq\a\Leftrightarrow r\leq 0,$$
$$\forall (r\d,i),(s\d,j)\in\tPhi_+^{{\rm{im}}}:\ \ \ 
(r\d,i)\preceq(s\d,j)\Leftrightarrow\cases
r>s&{\rm{or}}\cr
r=s\ \ {\rm{and}}\ \ i\leq j.\endcases$$
\rem{\cov}
$\preceq$ is a convex ordering of $\tPhi_+$ (see \bebr). 
\nota{\prb}
$\forall\eta\in Q$ define $\Pp(\eta)$ as follows: 
$$\Pp(\eta)\ugdef\big\{\ug=(\gamma_1,...,\gamma_r)
\in\bigcup_{r\in\mn}\tPhi_+^r|
\gamma_1\preceq...\preceq\gamma_r, \sum_{u=1}^rp(\gamma_u)=\eta
\big\};$$ 
define also $\Pp$ to be $\Pp\ugdef\bigcup_{\eta\in Q}\Pp(\eta)$.

We denote by $\tp:\bigcup_{r\in\mn}\tPhi_+^r\rightarrow\Pp$ the natural
(reordering) map; moreover $\forall\a\in\tPhi_+$ 
$\mu^{\a}:\Pp\rightarrow\mn$ ($\ug\mapsto\mu_{\ug}^{\a}$) indicates the
occurrence of $\a$ in $\ug$, that is if $\ug=(\gamma_1,...,\gamma_r)$
then 
$\mu_{\ug}^{\a}\ugdef\#\{u=1,...,r|\gamma_u=\a\}$.
\ddefi{\abnord}
By abuse of notation we denote by $\preceq$ also the
induced lessicographical ordering of $\Pp$: if
$\ug=(\gamma_1,...,\gamma_r)$ and
$\ug^{\prime}=(\gamma_1^{\prime},...,\gamma_s^{\prime})$ we set  
$$\ug\preceq\ug^{\prime}\Leftrightarrow {\rm{if}}\ \ 
l={\rm{min}}\{u=1,...,r|\gamma_u\neq\gamma_u^{\prime}\}\ \ {\rm{then}}\ \ 
\gamma_l\prec\gamma_l^{\prime}.$$
\vskip .2 truecm

{\bf{\S \pbls. PBW-basis and L-S formula.}}

\vskip .1 truecm
\nota{\sqnlt}
Let $x:\tPhi_+\rightarrow\uq$ ($\a\mapsto x_{\a}$) be any function. 
Then $\forall\ug=(\gamma_1,...,\gamma_r)\in\Pp$ we set
$x(\ug)\ugdef x_{\gamma_1}\cdot...\cdot x_{\gamma_r}$.

We are now ready to state the main theorem of this section.
\teo{\pbwt}
The set
$\{E(\ug)|\ug\in\Pp\}$ is a $\mc(q)$-basis of $\u_q^+$; more precisely
$\forall\eta\in Q$ $\{E(\ug)|\ug\in\Pp(\eta)\}$ is a $\mc(q)$-basis of
$\u_{q,\eta}^+$, while more generally 
$\{E(\ug)K_{\lambda}F(\ug^{\prime})|\lambda\in Q,\ug,\ug^{\prime}\in\Pp\}$
is a
$\mc(q)$-basis of
$\u_q$.
\dim
The argument used in \lusz\ applies here word for word.\finedim
The following theorem is the Levendorskii-Soibelman formula:
\teo{\lsfr} 
1) $\forall\a,\beta\in\tPhi_+$, $\a\succ\beta\Rightarrow
E_{\a}E_{\beta}-q^{(p(\a)|p(\beta))}E_{\beta}E_{\a}=
\sum_{\ug\succ(\beta,\a)}a_{\ug}E(\ug);$

2) $\forall\ug,\ug^{\prime}\in\Pp\ $, $\ug\succ\ug^{\prime}\Rightarrow
E(\ug)E(\ug^{\prime})-q^{r(\ug;\ug^{\prime})}E(\ug^{\prime})E(\ug)=
\sum_{\tilde{\ug}\succ\tp(\ug,\ug^{\prime})}a_{\tilde{\ug}}E(\tilde{\ug})$,
where $(\ug,\ug^{\prime})\mapsto r(\ug;\ug^{\prime})$ is the only function
compatible with 1). 
\dim
See \lswg\ and \beck.\finedim
\vskip .2 truecm

{\bf{\S \sntt. Some subspaces of $\u_q^+$.}}

\vskip .1 truecm
In this paragraph we first list some notations and then give a
characterization of some subspaces of $\u_q^+$.
\nota{\notinf}
$\forall r\in\mz\cup\{\pm\infty,\pm\tilde\infty,{\rm{im}}\}$ we define a
a subset $\tPhi_+(r)$ of $\tPhi_+$:
$$\tPhi_+(r)=\cases
\{\beta_s|1\leq s<r\}&{\rm{if}}\ r\geq 1\cr
\{\beta_s|0\geq s>r\}&{\rm{if}}\ r\leq 0\cr
\tPhi_+^{{\rm{im}}}&{\rm{if}}\ r={\rm{im}}\cr
\bigcup_{s>0}\tPhi_+(\pm s)&{\rm{if}}\ r=\pm\infty\cr
\tPhi_+^{{\rm{im}}}\cup\tPhi_+(\pm\infty)&{\rm{if}}\ r=\pm\tilde\infty;
\endcases$$
moreover we define $\u_q^+(r)$ and $\u_q^{\geq 0}(r)$ by:
$$\u_q^+(r)\ugdef{\rm{span}}_{\mc(q)}\{E(\ug)|\ug=(\gamma_1,...,\gamma_s)\in\Pp,\gamma_u\in
\tPhi_+(r)\ \forall u=1,...,s\},$$ 
$\u_q^{\geq 0}(r)\ugdef\u_q^0\u_q^+(r)$.

In order to describe $\u_q^+(r)$ we need the following lemma.
\lem{\stbwg}
$\forall\a\in\tPhi_+^{{\rm{im}}}$ we have that 
$$T_{w_r}^{-1}(E_{\a})\in\u_q^+\ \ \forall r\geq 1,\ \ \ 
T_{w_r^{-1}}(E_{\a})\in\u_q^+\ \ \forall r\leq 0.$$
\dim
Following \bebr, we have to prove that $\forall j\in I_0$ if 
$\lambda_j=\tau_j s_{j_1}\cdot...\cdot s_{j_d}$ then 
$T_{j_r}\cdot...\cdot T_{j_d}(E_{\a})\in\u_q^+$ $\forall r=1,...,d$; let
$\a=(m\dk_i\d,i)$ with $i\in I_0$, $m>0$. 

If $i\neq j$ the claim is obvious; if $i=j$ and 
$(X_{\tn}^{(k)},i)\neq(A_{2n}^{(2)},1)$, see 
\bebr. 

Finally we need to study the case when 
$(X_{\tn}^{(k)},i,j)=(A_{2n}^{(2)},1,1)$: first of all remark that 
$j_d=1$ and that $\forall m>0$
$$\varphi_1\Xi(E_{m\d-\a_1})=\varphi_1\Xi(T_0T_1)^mT_1^{-1}(E_1)=
\varphi_1(T_0^{-1}T_1^{-1})^mT_1(E_1)=$$
$$=\varphi_1 T_1(T_0T_1)^{-m}(E_1)=T_1T_{\omega_1}^{-m}(E_1)
$$
which belongs to $\u_q^+$; then 
$$T_1(E_{(m\d,1)})=\varphi_1T_1(E_{m\d})=\varphi_1\Xi(E_{m\d})=$$ 
$$=-E_1T_1T_{\omega_1}^{-m}(E_1)+q^{-2}T_1T_{\omega_1}^{-m}(E_1)E_1\in
\u_q^+.$$
The claim now follows remarking that 
$l(s_{j_1}\cdot...\cdot s_{j_{d-1}}s_1)=d$, that $m>0$ and that 
$T_{j_r}\cdot...\cdot T_{j_{d-1}}T_1T_{\omega_1}^{-1}=
(T_{j_1}\cdot...\cdot T_{j_{r-1}})^{-1}$.
\finedim
\prop{\descrinf}
$\forall r\in\mz\cup\{\pm\infty,\pm\tilde\infty,{\rm{im}}\}$ $\u_q^+(r)$
is a subalgebra of $\u_q^+$. 

Moreover
$$\u_q^+(\pm\infty)=\{x\in\u_q^+|\exists m>0\ {\rm{s.t.}}\ 
T_{w_{\pm m}^{\pm 1}}^{\mp 1}(x)\in\u_q^{\leq 0}\}=$$ 
$$=\{x\in\u_q^+|
T_{w_{\pm m}^{\pm 1}}^{\mp 1}(x)\in\u_q^{\leq 0}\ \forall m>>0\},$$
$$\u_q^+(\pm\tilde\infty)=\{x\in\u_q^+|
T_{w_{\mp m}^{\mp 1}}^{\pm 1}(x)\in\u_q^+\ \forall m>0\},$$
$$\u_q^+({\rm{im}})=\u_q^+(\tilde\infty)\cap\u_q^+(-\tilde\infty)=
\{x\in\u_q^+|T_{\lambda_i}(x)=x\ \forall i\in I_0\}.$$
\dim
The first assertion follows from L-S formula and from the convexity of
$\preceq$; for the second it is enough to apply lemma \stbwg\ (see \bebr\
and \rmens), while the third follows immediately from the
preceding ones, remarking that $T_{\lambda_i}(E_{\a})=E_{\a}$
$\forall i\in I_0$ $\forall\a\in\tPhi_+^{{\rm{im}}}$.\finedim

\vskip .5 truecm

{\bf{\rmat. THE $R$-MATRIX.}}
\vskip .2 truecm
Here it will be shown, applying to the present situation 
the strategy  
and the results used 
for the non twisted
quantum algebras (see \lss\ and \rmens),
how 
it is possible to exhibit an explicit
multiplicative formula for the $R$-matrix (denoted by $R$) in the
twisted case. 

Hence most of this 
section will be devoted to follow the fundamental steps
which lead to the description of the $R$-matrix for the non twisted quantum
algebras, to recall the main properties used to establish them, and to
remark how, thanks to the results developed in the preceding parts of
this paper, the same properties (and consequently the same results) are
still valid in this more general setting: we shall 
refer to \lsrm\ and \rmens\ for the proofs. Also, some minor differences
occurring in the twisted case (mainly in case $A_{2n}^{(2)}$) will be
discussed, and the computations needed to complete the argument performed. 
\vskip .2 truecm
{\bf{\S \rmgn. General strategy for the construction of the $R$-matrix.}}
\vskip .1 truecm
First of all we recall a very general result connecting the $R$-matrix
to the Killing form (see
\tan). 
\teo{\rta}
Let $C_{\eta}\in\u_{q,\eta}^+\otimes\u_{q,-\eta}^-$ (with $\eta\in Q^+$)
be the canonical element of $(\cdot,\cdot)\big|_{\u_{q,\eta}^+\times 
\u_{q,-\eta}^-}$. 
Then $$R=\bigg(\sum_{\eta\in Q^+}C_{\eta}\bigg)q^{-t_{\infty}},$$
where $q^{-t_{\infty}}$ is the diagonal operator on $V\otimes\tilde V$
($V$ and $\tilde V$ $\uq$-modules) which acts as $q^{-(\eta|\tilde\eta)}$
on
$V_{\eta}\otimes\tilde V_{\tilde\eta}$.
\finedim
Then we recall the behaviour of the real root vectors with
respect to the coproduct and to the Killing form.
\prop{\dkrrv}
If $\a=\beta_r\in\Phi_+^{{\rm{re}}}$ then
$\Delta(E_{\a})-(E_{\a}\otimes 1+K_{\a}\otimes E_{\a})\in
\cases\u_q^{\geq 0}(r)\otimes\u_q^+&{\rm{if}}\ r\geq 1\cr
\u_q^{\geq 0}\otimes\u_q^+(r)&{\rm{if}}\ r\leq 0.\endcases$

Moreover let $\underline\gamma=(\gamma_1,...,\gamma_r)$,
$\underline\gamma^{\prime}=(\gamma_1^{\prime},...,\gamma_s^{\prime})\in\Pp$
be such that $\gamma_u,\gamma_v^{\prime}\in\tPhi_+(\infty)$ $\forall u,v$
or 
$\gamma_u,\gamma_v^{\prime}\in\tPhi_+(-\infty)$ $\forall u,v$; then 
$$(E(\underline\gamma),F(\underline\gamma^{\prime}))=
\d_{\underline\gamma,\underline\gamma^{\prime}}
\prod_{\a\in\tPhi_+}
\rapp{(\mu_{\underline\gamma}^{\a})_{\a}!}{(q_{\a}^{-1}-q_{\a})^{\mu_{\ug}^{\a}}}.$$
\dim
The proof given in \lsrm\ and \rmens\ depends only on the action of the
braid group on $\uq$ and on its connection with the coproduct
$\Delta$ (that is on the construction of ``partial
$R$-matrices'' and on their properties), hence it applies word for
word to the present situation.\finedim
Before going to the next proposition we shall introduce a notation and give
a simple lemma.
\nota{\dxijr}
We denote by det${}_rX_{\tn}^{(k)}$ the determinant of the matrix 
$$\Big({r[d_i]_q\over \dk_i[r]_q}x_{ijr}\Big)_{ij\in I^r}=
{r\over[r]_q}{\rm{diag}}\bigg({[d_i]_q\over \dk_i}\bigg|i\in I^r\bigg)
(x_{ijr})_{_{ij\in I^r}}$$ for the
quantum algebra of type
$X_{\tn}^{(k)}$.
\lem{\nondegxij}
$\forall r>0$ we have det${}_rX_{\tn}^{(k)}\neq 0$ (equivalently 
det$(x_{ijr})_{_{ij\in I^r}}\neq 0$).

More precisely
$${\rm{det}}_rX_{n}^{(1)}=\cases
[n+1]_{q^r}&{\rm{in\ case\ }}A_n^{(1)}\cr
[2]_{q^r}^{n-1}[2]_{q^{r(2n-1)}}&{\rm{in\ case\ }}B_n^{(1)}\cr
[2]_{q^r}[2]_{q^{r(n+1)}}&{\rm{in\ case\ }}C_n^{(1)}\cr
[2]_{q^r}[2]_{q^{r(n-1)}}&{\rm{in\ case\ }}D_n^{(1)}\cr
{[3]_{q^r}[2]_{q^{6r}}\over[2]_{q^{2r}}}&{\rm{in\ cases\ }}E_6^{(1)}\ 
{\rm{and}}\ G_2^{(1)}\cr
{[2]_{q^r}[2]_{q^{9r}}\over[2]_{q^{3r}}}&{\rm{in\ case\ }}E_7^{(1)}\cr
[2]_{q^r}[2]_{q^{7r}}-[3]_{q^r}={[2]_{q^r}[2]_{q^{15r}}\over
[2]_{q^{3r}}[2]_{q^{5r}}}&{\rm{in\
case\ }}E_8^{(1)}\cr 
{[2]_{q^r}^2[2]_{q^{9r}}\over[2]_{q^{3r}}}&{\rm{in\
case\ }}F_4^{(1)}\cr 
[2]_{q^r}^n[2n+1]_{q^r}&{\rm{in\ case}}\ A_{2n}^{(2)}\ \ {\rm{if}}\
2\!\not|r\cr 
[2]_{q^r}^{n-1}[2]_{q^{r(2n+1)}}&{\rm{in\ case}}\ A_{2n}^{(2)}\ \ {\rm{if}}\
2|r\cr
[{\tn-n\over k-1}+1]_{q^r}&{\rm{if}}\ \ \tk\not |r\cr 
[2]_{q^{rn}}&{\rm{in\ cases}}\ A_{2n-1}^{(2)}\ {\rm{and}}\
D_{n+1}^{(2)}\ {\rm{if}}\ 2|r\cr 
{[2]_{q^{6r}}\over[2]_{q^{2r}}}&{\rm{in\ case\ }}E_6^{(2)}\ {\rm{if}}\
2|r\cr
{[2]_{q^{3r}}\over[2]_{q^r}}&{\rm{in\ case\ }}D_4^{(3)}\ {\rm{if}}\
3|r.\endcases$$
\dim
The proof consists in simple computations: 
in the non twisted case we have that det${}_rX_{\tn}^{(k)}=
{\rm{det}}([d_ia_{ij}]_{q^r})_{ij\in I_0}$: see 
\rmcen\ and \rmens; the result when $\tk\not|r$
follows from the remark that 
${\rm{det}}_rX_{\tn}^{(k)}={\rm{det}}([a_{ij}^s]_{q^r})_{ij\in I^r}$ (and
from the result for 
$A_n^{(1)}$).
The remaining assertions follow considering that 
$${\rm{det}}_rA_{2n}^{(2)}=
(q^{2r}+(-1)^{r-1}+q^{-2r})[2]_{q^r}^n
{\rm{det}}_{2r}A_{n-1}^{(1)}-[2]_{q^r}^n{\rm{det}}_{2r}A_{n-2}^{(1)},$$
$$1<k|r,\ X_{\tn}^{(k)}=A_{2n-1}^{(2)}, D_{\tn}^{(k)}\Rightarrow
{\rm{det}}_rX_{\tn}^{(k)}=
[2]_{q^r}{\rm{det}}_rA_{n-1}^{(1)}-k{\rm{det}}_rA_{n-2}^{(1)},$$ 
$${\rm{det}}_rE_6^{(2)}=[2]_{q^r}{\rm{det}}_rD_3^{(2)}-
{\rm{det}}_rA_2^{(1)}.$$\finedim
\prop{\dkirv}
1) Let $\a\in\tPhi_+^{{\rm{im}}}$; then 
$\Delta(E_{\a})-(E_{\a}\otimes 1+K_{p(\a)}\otimes E_{\a})\in
\u_q^{\geq 0}(\infty)\otimes\u_q^+(-\infty)$;

2) Let $\a\in\tPhi_+$ and $\underline\gamma\in{\Pp}$ be such that 
$(E_{\a},F(\underline\gamma))\neq 0$ or $(E(\underline\gamma),F_{\a})\neq
0$; then $\underline\gamma=(\beta)$ with $p(\a)=p(\beta)$.

3) Let $\{\bar E_{\a}|\a\in\tPhi_+\}$ be elements of $\u_q^+$ satisfying
the following requirements:

i)  $\bar E_{\a}=E_{\a}$ $\forall\a\in\Phi_+^{{\rm{re}}}$; 

ii) $\forall r>0$ $\{\bar E_{(r\d,i)}|i\in I^r\}$ and 
$\{E_{(r\d,i)}|i\in I^r\}$ 
are bases of 
the same 
$\mc(q)$-vector subspace of $\uq$; 

iii) $(\bar E_{(r\d,i)},\bar F_{(r\d,j)})=0$ if $i\neq j$ (where $\bar
F_{\a}\ugdef\Omega(\bar E_{\a})$);

then $\forall\underline\gamma,\underline\gamma^{\prime}\in{\Pp}$
we have that 
$$(\bar E(\underline\gamma),\bar F(\underline\gamma^{\prime}))=
\d_{\underline\gamma,\underline\gamma^{\prime}}
\prod_{\a\in\Phi_+^{{\rm{re}}}}
\rapp{(\mu_{\underline\gamma}^{\a})_{\a}!}
{(q_{\a}^{-1}-q_{\a})^{\mu_{\underline\gamma}^{\a}}}
\prod_{\a\in\tPhi_+^{{\rm{im}}}}
\mu_{\underline\gamma}^{\a}!(\bar E_{\a},\bar
F_{\a})^{\mu_{\underline\gamma}^{\a}}.$$
\dim
For the proof see \rmens. Indeed the argument used in the non twisted case
is based on the following properties, which are still valid in the twisted
case:

a) $T_{\lambda_1\cdot...\cdot\lambda_n}(E_{\beta_r})=E_{\beta_{r+N}}$ if
$r>0$ or
$r\leq-N$ and $T_{\lambda_1\cdot...\cdot\lambda_n}(E_{\a})=E_{\a}$ if
$\a\in\tPhi_+^{{\rm{im}}}$; more precisely see proposition \descrinf;

b) the matrix $(x_{ijr})_{ij\in I^r}$ is non degenerate (see lemma
\nondegxij); 

c) $\prec$ is a convex ordering of $\tPhi_+$.\finedim
\cor{\rmaf}
Let $\bar E_{\a}$, $\bar F_{\a}$ ($\a\in\tPhi_+$) be as in
proposition \dkirv; then 
$$R=\bigg(\prod_{\a\in\tPhi_+}{\rm{exp}}_{\a}\Big(
\rapp{\bar E_{\a}\otimes\bar F_{\a}}{(\bar E_{\a},\bar
F_{\a})}
\Big)\bigg)q^{-t_{\infty}}=
\bigg(\prod_{\a\in\tPhi_+}{\rm{exp}}_{\a}\Big(
\rapp{q_{\a}^{-1}-q_{\a}}{c_{\a}}
\bar E_{\a}\otimes\bar F_{\a}
\Big)\bigg)
q^{-t_{\infty}}$$
where $c_{\a}=(q_{\a}^{-1}-q_{\a})(\bar E_{\a},\bar F_{\a})=\cases
1&{\rm{if}}\ \a\in\Phi_+^{{\rm{re}}}\cr 
(q_{\a}^{-1}-q_{\a})(\bar E_{\a},\bar F_{\a})&{\rm{if}}\
\a\in\tPhi_+^{{\rm{im}}}.\endcases$
\dim
See \lsrm\ and \rmens.\finedim
\rem{\rcan}
Remark that $\forall r>0$
$$\prod_{i\in I^{r}}{\rm{exp}}_{(r\d,i)}\rapp{\bar E_{(r\d,i)}\otimes\bar
F_{(r\d,i)}}{(\bar E_{(r\d,i)},\bar
F_{(r\d,i)})}={\rm{exp}}\bigg(\sum_{i\in I^{r}}\rapp{\bar
E_{(r\d,i)}\otimes\bar F_{(r\d,i)}}{(\bar E_{(r\d,i)},\bar
F_{(r\d,i)})}\bigg).$$ Hence corollary \rmaf\ can be rewritten as 
$$R\!\!=\!\!\!\!\prod_{m\leq 0}\!\!{\rm{exp}}_{\beta_m}
((q_{\beta_m}^{-1}-q_{\beta_m})E_{\beta_m}\otimes F_{\beta_m}\!)
\cdot\!\prod_{r>0}\!\!{\rm{exp}}\tc_r
\cdot\!\prod_{m\geq
1}\!\!{\rm{exp}}_{\beta_{m}}((q_{\beta_m}^{-1}-q_{\beta_m})E_{\beta_m}\otimes
F_{\beta_m})q^{-t_{\infty}}$$ 
where $\tc_r=\sum_{i\in
I^{r}}\rapp{\bar E_{(r\d,i)}\otimes\bar F_{(r\d,i)}}{(\bar
E_{(r\d,i)},\bar F_{(r\d,i)})}$ is the canonical element of the
restriction of
$(\cdot,\cdot)$ to 
$${\rm{span}}_{\mc(q)}\{E_{(r\d,i)}|i\in I^r\}\times{\rm{span}}_{\mc(q)}\{
F_{(r\d,i)}|i\in I^r\},$$
and both the first and the third products in the right hand side are
performed in decreasing order.\finedim
\vskip .2 truecm
{\bf{\S \kirv. The Killing form on the imaginary root vectors.}}
\vskip .1 truecm
We have seen how to reduce the explicit description of $R$ to the
construction of $\tc_m$ or equivalently to the construction of elements 
$\bar E_{\a}$ with the properties described in proposition \dkirv, and to
the computation of $(\bar E_{\a},\bar F_{\a})$ for
$\a\in\tPhi_+^{{\rm{im}}}$. Of course for both approaches the problem is
that of understanding $(\cdot,\cdot)$. Thus we start from the computation
of 
$(E_{(r\d,i)},F_{(r\d,j)})$, following the path indicated in
\rmens\ and remarking that in case $A_{2n}^{(2)}$ it has to be slightly
modified: the problem arises, here too (see lemma \basbec), from the fact
that in this case $\d-2\a_1$ is a root.
\ddefi{\prcff}
Let $i\in I_0$; we shall denote by $\pi_i$ the natural projection 
$\pi_i:\u_q\otimes\u_q^{\geq 0}\rightarrow\u_q\otimes\mc(q)E_i$
(with respect to the $Q$-gradation). 

If $r>0$ and $i,j\in I^r$ we define $c_r^{ij}$ by the
following property:
$$\pi_j(\Delta(E_{(r\d,i)}))=c_r^{ij}E_{r\d-\a_j}K_j\otimes E_j$$
($c_r^{ij}$ is well defined since $\pi_j(\Delta(E_{(r\d,i)}))$
is indeed a multiple of $E_{r\d-\a_j}K_j\otimes E_j$: see \rmens,
using proposition \dkirv, 1) and the convexity of $\prec$).
\lem{\erc}
Let $r>0$ and $i,j\in I^r$; then 
$(E_{(r\d,i)},F_{(r\d,j)})=-q_j^2{c_r^{ij}\over (q_j^{-1}-q_j)^2}$. 
\dim
The proof is based on proposition \dkirv, 2), on \TTTUQ, 7,i), and on
proposition \dkrrv\ (see \rmens).
\finedim
\lem{\cij}
Let $r>0$ and $i,j\in I^r$; then $c_r^{ij}=q_j^{-2}(q_j-q_j^{-1})x_{ijr}$.
\dim
The proof consists in showing that $c_r^{ij}={c_1^{jj}\over
x_{jj1}}x_{ijr}$ and in computing ${c_1^{jj}\over
x_{jj1}}
$. 

For $c_r^{ij}={c_1^{jj}\over
x_{jj1}}x_{ijr}$ see \rmens. 

For ${c_1^{jj}\over x_{jj1}}=q_j^{-2}(q_j-q_j^{-1})$ see
\rmens, 
when we are not in the case $A_{2n}^{(2)}$ with
$j=1$. 

So it remains to study $c_1^{11}$ just in case $A_{2n}^{(2)}$: remark that
$\pi_1(\Delta(E_{\d-2\a_1}))=0$ (because
$\d-3\a_1\not\in\tPhi$); 
then, using corollary \am21, we have
that 
$$\pi_1(\Delta(E_{\d-\a_1}))=-E_{\d-2\a_1}K_1\otimes E_1+
q^{-4}K_1E_{\d-2\a_1}\otimes E_1=-(1-q^{-8})E_{\d-2\a_1}K_1\otimes E_1;$$
it follows that 
$$\pi_1(\Delta(E_{(\d,1)}))=\pi_1(\Delta(E_{\d-\a_1}))(E_1\otimes
1)-q^{-2}(E_1\otimes 1)\pi_1(\Delta(E_{\d-\a_1}))+$$
$$+E_{\d-\a_1}K_1\otimes
E_1-q^{-2}K_1E_{\d-\a_1}\otimes
E_1=(q^2+1-q^{-4}-q^{-6})E_{\d-\a_1}K_1\otimes E_1,$$ which means that
$c_1^{11}=q^2+1-q^{-4}-q^{-6}=
q_1^{-2}(q_1-q_1^{-1})x_{111}$.\finedim
\cor{\kc}
Let $r>0$ and let $x$ (respectively $y$) belong to the linear span 
of $\{E_{(r\d,i)}|i\in I^r\}$ (respectively $\{F_{(r\d,i)}|i\in I^r\}$).

Then $[x,y]=-(x,y)(K_{r\d}-K_{r\d}^{-1})$
and $(\Omega(y),\Omega(x))=-\Omega((x,y))$.
\dim
The claim follows from lemmas \erc\ and \cij, from theorem \cascomm,
from the bilinearity of $(\cdot,\cdot)$ and 
$[\cdot,\cdot]$ and from the fact that $\Omega$ is an
antiautomorphism.\finedim
\vskip .2 truecm
{\bf{\S \rsnsd. A ``canonical'' form for the $R$-matrix.}}
\vskip .1 truecm 
Remark \rcan\ reduces the problem of describing $R$ to that of finding
$\tc_r$ $\forall r>0$. The goal of the present paragraph, which
can be achieved thanks to the results of \S \kirv, is to describe $\tc_r$ 
in terms of $\{E_{(r\d,i)},F_{(r\d,i)}|i\in I^r\}$.

The following is a simple lemma of linear algebra. 
\lem{\linalg}
Let $r>0$ and let $\tc_r=\sum_{i,j\in I^r}y_{ij}^rE_{(r\d,i)}\otimes
F_{(r\d,j)}$. Then if $y^r\ugdef(y_{ij}^r)_{ij\in I^r}$ and 
$H^r\ugdef(E_{(r\d,i)},
F_{(r\d,j)})_{ij\in I^r}$ we have $y^r=(^t\!H^r)^{-1}$.\finedim
\prop{\cmt}
In the notation of lemma \linalg, if we put
$$z_{ij}^{r}={(o(i)o(j))^{r\over
k_r^{\prime}}\dk_i[r]_q\over
r[d_i]_q[d_j]_q(q^{-1}-q)}y_{ij}^{r}\ \  {\rm{with}}\ \ 
k_r^{\prime}\ugdef\cases k&{\rm{if}}\ k|r,\ X=D,E\cr
1&{\rm{otherwise}},\endcases$$ we can describe $z_{ij}^r$ as follows: 
of course $z_{ij}^r\dk_j=z_{ji}^r\dk_i$ (because $H^r$ is symmetric); if
$i\leq j$ then 
$$\eqalign{
&A_n^{(1)}\hskip 1 truecm
z_{ij}^r={[i]_{q^r}
[n-j+1]_{q^r}\over[n+1]_{q^r}}\cr
&B_n^{(1)}\hskip 1 truecm
z_{ij}^r=\cases{[n-j+1]_{q^{2r}}
\over[2]_{q^{r(2n-1)}}}&{\rm{if}}\ \ i=1\cr
{[2]_{q^{r(2i-3)}}
[n-j+1]_{q^{2r}}
\over[2]_{q^r}[2]_{q^{r(2n-1)}}}&{\rm{otherwise}}\endcases\cr 
&C_n^{(1)}\hskip 1 truecm
z_{ij}^r=\cases{[n]_{q^r}\over[2]_{q^r}[2]_{q^{r(n+1)}}}
&{\rm{if}}\ \ i=j=1\cr
{[n-j+1]_{q^r}\over[2]_{q^{r(n+1)}}}
&{\rm{if}}\ \ 
j>1=i\cr
{[2]_{q^{ri}}
[n-j+1]_{q^r}\over[2]_{q^{r(n+1)}}}&
{\rm{otherwise}}\endcases\cr 
&D_n^{(1)}\hskip 1 truecm z_{ij}^r=\cases
{[n]_{q^r}\over[2]_{q^r}[2]_{q^{r(n-1)}}}&{\rm{if}}\ \ i=j\leq 2\cr
{[n-2]_{q^r}\over[2]_{q^r}[2]_{q^{r(n-1)}}}&{\rm{if}}\ \ i=1,j=2\cr
{[n-j+1]_{q^r}
\over[2]_{q^{r(n-1)}}}&{\rm{if}}\ \ i\leq 2<j
\cr 
{[2]_{q^{r(i-2)}}[n-j+1]_{q^r}\over[2]_{q^{r(n-1)}}}&{\rm{otherwise}}
\endcases\cr}$$
$$\eqalign{
&E_n^{(1)}\hskip 1 truecm {\rm{det}}_rE_n^{(1)}z_{ij}^r=\cases
[n]_{q^r}&{\rm{if}}\ \ i=j=1\cr
[j-1]_{q^r}[n-3]_{q^r}&{\rm{if}}\ \ i=1<j\leq 3\cr
[2]_{q^r}[2]_{q^{r(n-2)}}&{\rm{if}}\ \ i=j=2\cr
[i-1]_{q^r}[n-1]_{q^r}&{\rm{if}}\ \ 2\leq i\leq 3=j\cr
[3]_{q^r}[n-j+1]_{q^r}&{\rm{if}}\ \ i=1,j>3\cr
[2]_{q^r}[i-1]_{q^r}[n-j+1]_{q^r}&{\rm{if}}\ \ 2\leq i\leq 4\leq j\cr
[5]_{q^r}[n-j+1]_{q^r}&{\rm{if}}\ \ i=5\cr
[2]_{q^r}[2]_{q^{4r}}[n-j+1]_{q^r}&{\rm{if}}\ \ i=6\cr 
{[3]_{q^r}[2]_{q^{6r}}[n-j+1]_{q^r}\over[2]_{q^{2r}}}&{\rm{if}}\ \ i=7\cr 
{[2]_{q^r}[2]_{q^{9r}}\over[2]_{q^{3r}}}&{\rm{if}}\ \ i=j=8
\endcases\cr
&F_4^{(1)}\hskip 1 truecm z_{ij}^r=\cases
{[2]_{q^{3r}}[2]_{q^{5r}}\over[2]_{q^{9r}}}&{\rm{if}}\
\ i=j=1\cr
{[2]_{q^{3r}}[i]_{q^r}[5-j]_{q^{2r}}\over[2]_{q^{9r}}}&{\rm{if}}\
\ i\leq 2\leq j\cr {[2]_{q^{3r}}[3]_{q^r}[5-j]_{q^{2r}}
\over[2]_{q^{9r}}[2]_{q^r}}&{\rm{if}}\ \ i=3
\cr 
{[2]_{q^{3r}}[2]_{q^{4r}}\over[2]_{q^{9r}}[2]_{q^r}}&{\rm{if}}\ \ i=j=4
\endcases\cr
&G_2^{(1)}\hskip 1 truecm
(z_{ij}^r)_{ij\in I_0}={[2]_{q^{2r}}\over[2]_{q^{6r}}[3]_{q^r}}
\left(\matrix[6]_{q^r}&[3]_{q^r}\cr[3]_{q^r}&[2]_{q^r}\endmatrix\right)
\cr 
&A_{2n}^{(2)}\hskip 1 truecm z_{ij}^r=\cases
{[2i-1]_{q^r}
[n-j+1]_{q^{2r}}\over[2]_{q^r}[2n+1]_{q^r}}&{\rm{if}}\
2\not|r\cr
{[2]_{q^{r(2i-1)}}
[n-j+1]_{q^{2r}}\over[2]_{q^r}[2]_{q^{r(2n+1)}}}&{\rm{if}}\
2|r
\endcases\cr
&A_{2n-1}^{(2)}\hskip 1 truecm z_{ij}^r=\cases
{[i-1]_{q^r}
[n-j+1]_{q^r}\over[n]_{q^r}}&{\rm{if}}\ \
k\not|r\cr
{[n]_{q^{r}}\over[2]_{q^{rn}}}&{\rm{if}}\ k|r,i=j=1\cr
{(-1)^{{r\over 2}}2[n-j+1]_{q^r}\over [2]_{q^{rn}}}&{\rm{if}}\ \ k|r, 
i=1<j\cr
{[2]_{q^{r(i-1)}}[n-j+1]_{q^r}\over
[2]_{q^{rn}}}&{\rm{otherwise}}\endcases\cr  
&D_{n+1}^{(2)}\hskip 1 truecm
z_{ij}^r=\cases
{1\over[2]_{q^r}}&{\rm{if}}\ \ k\not|r\cr
{[n-j+1]_{q^r}\over
[2]_{q^{rn}}}&{\rm{if}}\ \ k|r, i=1\cr
{[2]_{q^{r({i-1)}}}[n-j+1]_{q^r}\over
[2]_{q^{rn}}}&{\rm{otherwise}}\endcases\cr 
&E_6^{(2)}\hskip 1 truecm
z_{ij}^r=\cases
{[i]_{q^r}[3-j]_{q^r}\over[3]_{q^r}}\ \ {\rm{if}}\ \
k\not|r\cr
{[2]_{q^{2r}}[2]_{q^{3r}}\over[2]_{q^{6r}}}&{\rm{if}}\
\ k|r, i=j\in\{1,4\}\cr 
{(-1)^{{r\over 2}}[2]_{q^{2r}}[5-j]_{q^r}\over[2]_{q^{6r}}}&{\rm{if}}\ \
k|r, i=1<j\cr 
{[2]_{q^{2r}}[i]_{q^r}[5-j]_{q^r}\over[2]_{q^{6r}}}&{\rm{otherwise}}
\endcases\cr
&D_4^{(3)}\hskip 1 truecm
z_{ij}^r=\cases
{1\over[2]_{q^r}}&{\rm{if}}\ \ k\not|r\cr
{[2]_{q^r}[i]_{q^r}[3-j]_{q^r}\over[2]_{q^{3r}}}&{\rm{if}}\ \
k|r\endcases\cr}. $$
\dim
Thanks to lemma \linalg, the proof consists in the simple
exercise of inverting the matrix
$(x_{ijr})_{ij\in I^r}$.\finedim
\teo{\rfln}
With the notations of proposition \cmt, we have that 
$Rq^{t_{\infty}}=$
$$=\!\!\!\prod_{r\leq
0}\!\!{\rm{exp}}_{\beta_r}\!(q_{\beta_r}^{-1}-q_{\beta_r}\!)E_{\beta_r}\otimes
F_{\beta_r}
\!\!\!\prod_{r>0}\!\!{\rm{exp}}(\!\!\sum_{ij\in
I_r}\!\!\!y_{ij}^r E_{(r\d,i)}\otimes F_{(r\d,j)})\!\!\prod_{r>0}\!\!
{\rm{exp}}_{\beta_r}\!(q_{\beta_r}^{-1}-q_{\beta_r}\!)E_{\beta_r}\!\otimes
F_{\beta_r}.$$
\dim
The claim is an immediate consequence of remark \rcan, of lemma \linalg\
and of proposition \cmt.\finedim
\vskip .2 truecm
{\bf{\S \teal. Linear triangular transformation: the $\bar E_{\a}$'s.}}
\vskip .1 truecm
It can be useful to determine explicitly the expression of $R$ in the form
given in corollary \rmaf, that is to describe $R$ in terms of ``dual''
bases. To this aim we shall look for elements $\bar E_{\a}$ with
$\a\in\tPhi_+$ satisfying the requirements of proposition \dkirv. 
\lem{\pvlf}
For $r>0$ let $\bar E_{(r\d,i)}$ ($i\in I^r$) be elements of $\u_{q,r\d}^+$
such that:

1) $\forall i\in I^r$ $\bar E_{(r\d,i)}$ is a linear combination of 
$\{E_{(r\d,j)}|j\in I^r, j\geq i\}$;

2) $\forall i<j\in I^r$ $(\bar E_{(r\d,i)},F_{(r\d,j)})=0$.

Then, if $\bar F_{(r\d,i)}\ugdef\Omega(\bar E_{(r\d,i)})$, we have that 
$(\bar E_{(r\d,i)},\bar F_{(r\d,j)})=0$ $\forall i\neq j\in I^r$. 
\dim
The assertion is clear for $i<j$ once that one notices that $\bar
F_{(r\d,j)}$ lies in the linear span of $\{F_{(r\d,l)}|l\in I^r, l\geq
j\}$; for $i>j$, it follows from corollary \kc.\finedim
\rem{\wvlf}
Lemma \pvlf\ allows us to reduce the construction of 
$\bar E_{(r\d,i)}$ to the solution of the following system of linear
equations:
$$\sum_{l\geq i}A_{il}^{(r)}x_{ljr}=0\ \ \ {\rm{for}}\
j\in I^r\ {\rm{with}}\ i<j\leqno{(*)_{r,i}}$$
under the condition $A_{ii}^{(r)}\neq 0$, 
where $\bar E_{(r\d,i)}=\sum_{l\geq
i}A_{il}^{(r)}E_{(r\d,l)}$.

Notice that in the non twisted case
$$x_{ljr}=(o(l)o(j))^r{[ra_{lj}]_{q_l}\over
r}=(o(l)o(j))^r{[r]_q[d_la_{lj}]_{q^r}\over r[d_l]_q}$$ so that 
$(*)_{r,i}$ can be written as 
$$\sum_{l\geq i}{A_{il}^{(r)}o(l)^r\over[d_l]_q}[d_la_{lj}]_{q^r}=0
\ \ \ {\rm{for}}\ j\in I_0\ {\rm{with}}\ i<j.\leqno{(*)_{r,i}^{(1)}}$$ 
This system has been solved in \rmcen, and those solutions will be used
here in order to simplify the computations in the general case. 
\prop{\sntw}
Consider the non twisted case ($k=1$) and let $r>0$, $i\in I_0$. 

A solution of $(*)_{r,i}^{(1)}$ is given by: 
$${A_{ij}^{(r)}o(j)^r\over[d_j]_q}=\cases
[n-j+1]_{q^{2r}}&{\rm{in\ case\ }}B_n^{(1)}{\rm{and}}\cr
&{\rm{if\ }}(i,j)\neq(1,1)\ {\rm{in\ case\ }}F_4^{(1)}\cr
[n-j+1]_{q^r}[2]_{q^r}&{\rm{if\ }}j>i=1\ {\rm{in\ case\
}}C_n^{(1)}{\rm{and}}\cr 
&{\rm{if\ }}j-1>i=1\ {\rm{in\ case\
}}D_n^{(1)}\cr [n-2]_{q^r}&{\rm{if\ }}i=1,j=2\ {\rm{in\ case\
}}D_n^{(1)}\cr [j-1]_{q^r}[n-3]_{q^r}&{\rm{if\ }}4>j>i=1\ {\rm{in\ case\
}}E_n^{(1)}\cr [n-j+1]_{q^r}[3]_{q^r}&{\rm{if\ }}j\geq 4,i=1\ {\rm{in\
case\ }}E_n^{(1)}\cr 
[2]_{q^{r(n+1)}}&{\rm{if\ }}i=j=1\ {\rm{in\ cases\
}}F_4^{(1)}\ {\rm{and}}\ G_2^{(1)}\cr
[n-j+1]_{q^r}&{\rm{otherwise,}}
\endcases$$
where $j\in I_0$ is such that $i\leq j$.
\dim
See \rmcen.\finedim
We can now pass to the study of the twisted case. We shall consider three
different situations: the case $A_{2n}^{(2)}$, the twisted case different
from $A_{2n}^{(2)}$ with $k\not|r$ and the twisted case different
from $A_{2n}^{(2)}$ with $k|r$. 
\rem{\adnd}
Remark that $x_{ijr}(A_{2n}^{(2)})=x_{ijr}(B_n^{(1)})$ if $(i,j)\neq(1,1)$
and that in the equations $(*)_{r,i}$ $x_{11r}$ never appears; then 
$(*)_{r,i}(A_{2n}^{(2)})=(*)_{r,i}(B_n^{(1)})$; in particular the
solutions of $(*)_{r,i}$ in case $A_{2n}^{(2)}$ are the same as the 
solutions of $(*)_{r,i}$ in case $B_n^{(1)}$.
\rem{\tkndr}
Let $X_{\tn}^{(k)}\neq A_{2n}^{(2)}$ and $r>0$ be such that $k\not|r$. Then
$\forall i\in I^r$ (that is $\forall i\in I$ such that $d_i=1$) the system
$(*)_{r,i}(X_{\tn}^{(k)})
=(*)_{r,i}\Big(A_{{\tn-n\over k-1}}^{(1)}\Big)$
(with a shift of indices in case $A_{2n-1}^{(2)}$, that is 
$A_{ij}^{(r)}(A_{2n-1}^{(2)})=A_{i-1j-1}^{(r)}\Big(A_{{\tn-n\over
k-1}}^{(1)}\Big)=A_{i-1j-1}^{(r)}(A_{n-1}^{(1)})$
$\forall i,j\in I^r$). 
\rem{\tkdr}
In the twisted case different from $A_{2n}^{(2)}$, remarking
that
$\dk_i\cnt_{ij}=\dk_j\cnt_{ji}$ when $a_{ij}\neq 0$, the system $(*)_{r,i}$
becomes the following (under the condition that $k|r$): 
$$\sum_{l\geq
i}{A_{il}^{(r)}\over[d_l]_q}(o(j)o(l))^{{r\over
\dk_j\cnt_{jl}}}\cnt_{jl}[a_{lj}^s]_{q^r}=0
\ \ \ {\rm{for}}\ j\in I_0\ {\rm{with}}\
i<j.\leqno{(*)_{r,i}^{{\rm{tw}}}}$$
Notice that if $i<j\in I_0$ are such that $\cnt_{jl}=1$ $\forall l\geq i$
such that $a_{lj}\neq 0$, then $\dk_j=k_r^{\prime}$ and each equation above
can be written as 
$$\sum_{l\geq
i}{A_{il}^{(r)}o(l)^{{r\over
k_r^{\prime}}}\over[d_l]_q}[a_{lj}^s]_{q^r}=0
\ \ \ {\rm{for}}\ j\in I_0\ {\rm{with}}\
i<j.\leqno{(*)_{r,i,j}^{k|r}}$$
This always happens when $X_{\tn}^{(k)}=D_{\tn}^{(k)}$ ($k=2,3$); on the
other hand when $X_{\tn}^{(k)}=A_{2n-1}^{(2)}$ or $E_6^{(2)}$
we have that $\cnt_{jl}\neq 1\Rightarrow j=2$,
which implies that if $i\neq 1$ then 
$(*)_{r,i}^{{\rm{tw}}}=\{(*)_{r,i,j}^{k|r}|j>i\}$. Moreover we 
see also that 
$$(*)_{r,1}^{{\rm{tw}}}=\{(*)_{r,1,j}^{k|r}|j>2\}\cup
\Big\{\sum_{l\in I_0}{A_{1l}^{(r)}\over[d_l]_q}(o(2)o(l))^{{r\over
\dk_2\cnt_{2l}}}\cnt_{2l}[a_{l2}^s]_{q^r}=0\Big\}.$$
Remark also that $A_{11}^{(r)}$ does not appear in the equation 
$(*)_{r,i,j}^{k|r}$ when $(i,j)\neq(1,2)$, so that the only condition on
$A_{11}^{(r)}$ is the one arising from the last equation, i.e. 
${A_{11}^{(r)}\over[d_1]_q}
(o(2)o(1))^{{r\over
\dk_2\cnt_{21}}}\cnt_{21}=
\sum_{l>1}{A_{1l}^{(r)}\over[d_l]_q}(o(2)o(l))^{{r\over
\dk_2\cnt_{2l}}}\cnt_{2l}[a_{l2}^s]_{q^r}$, where the $A_{ij}^{(r)}$'s with
$(i,j)\neq (1,1)$ are solutions of the system $\{(*)_{r,i,j}^{k|r}|i<j\
{\rm{and\ }}(i,j)\neq(1,2)\}$.
\prop{\sraij}
Let $r>0$ and $i\in I^r$. 
A solution of $(*)_{r,i}$ is given by: 
$${A_{ij}^{(r)}o(j)^{{r\over k_r^{\prime}}}\over[d_j]_q}=\cases
[n-j+1]_{q^{2r}}&{\rm{in\ cases\ }}B_n^{(1)}{\rm{
and\ }}A_{2n}^{(2)},{\rm{and}}\cr &{\rm{if\ }}(i,j)\neq(1,1)\ {\rm{in\
case\ }}F_4^{(1)}\cr 
[n-j+1]_{q^r}[2]_{q^r}&{\rm{if\ }}j>i=1\
{\rm{in\ case\ }}C_n^{(1)}{\rm{and}}\cr 
&{\rm{if\ }}j-1>i=1\ {\rm{in\ case\
}}D_n^{(1)}\cr 
[n-2]_{q^r}&{\rm{if\ }}i=1,j=2\ {\rm{in\ case\
}}D_n^{(1)}\cr 
[j-1]_{q^r}[n-3]_{q^r}&{\rm{if\ }}4>j>i=1\ {\rm{in\ case\
}}E_n^{(1)}\cr 
[n-j+1]_{q^r}[3]_{q^r}&{\rm{if\ }}j\geq 4,i=1\ {\rm{in\
case\ }}E_n^{(1)}\cr 
[2]_{q^{5r}}&{\rm{if\ }}i=j=1\ {\rm{in\ case\
}}F_4^{(1)}\cr
[3-j]_{q^{3r}}&{\rm{if\ }}i=1\ {\rm{in\ case\
}}G_2^{(1)}\cr
[{\tn-n\over k-1}-j+1]_{q^r}&{\rm{if\ }}k\not|r\ {\rm{in\ cases\
}}D_{\tn}^{(k)}\ {\rm{and\ }}E_6^{(k=2)}\cr
{(-1)^{r\over 2}\over 2}[n]_{q^r}&{\rm{if\ }}2|r\ {\rm{and\ }}i=j=1\ 
{\rm{in\ case\ }}A_{2n-1}^{(2)}\cr
(-1)^{r\over 2}[2]_{q^{3r}}&{\rm{if\ }}2|r\ {\rm{and\ }}i=j=1\ 
{\rm{in\ case\ }}E_6^{(2)}\cr
[n-j+1]_{q^r}&{\rm{otherwise,}}
\endcases$$
where $j\in I^r$ is such that $i\leq j$.
\dim
The only statements which do not follow immediately from proposition
\sntw\ and remarks \adnd, \tkndr\ and \tkdr\ are those concerning
the case when $X_{\tn}^{(2)}=A_{2n-1}^{(2)}$ or
$E_6^{(2)}$, $2|r$ and $i=j=1$, and it is a matter of straightforward
computations.\finedim
\vskip .2 truecm
{\bf{\S \kniv. The Killing form on the new imaginary root vectors.}}
\vskip .1 truecm
The next (and last) step is computing 
$(\bar E_{\a},\bar F_{\a})$ when $\a$ is an imaginary root. 
\lem{\kbef}
Let $r>0$ and $i\in I^r$. Then 
$$(\bar E_{(r\d,i)},\bar F_{(r\d,i)})=
-A_{ii}^{(r)}\sum_{j\in I^r:j\geq i}A_{ij}^{(r)}{x_{jir}\over
q_i-q_i^{-1}}=
-{A_{ii}^{(r)}\over q-q^{-1}}\sum_{j\in I^r:j\geq
i}{A_{ij}^{(r)}\over[d_j]_q}x_{ijr}.$$
\dim
The definition of $\bar E_{\a}$ implies 
that $(\bar E_{(r\d,i)},\bar F_{(r\d,i)})=A_{ii}^{(r)}(\bar
E_{(r\d,i)},F_{(r\d,i)})$, because 
$(\bar E_{(r\d,i)},\bar F_{(r\d,i)}-A_{ii}^{(r)}F_{(r\d,i)})=0$.
But 
$$(\bar E_{(r\d,i)},F_{(r\d,i)})=
\sum_{j\in I^r:j\geq i}A_{ij}^{(r)}( E_{(r\d,j)},F_{(r\d,i)})=
-\sum_{j\in I^r:j\geq i}A_{ij}^{(r)}{x_{jir}\over
q_i-q_i^{-1}}.$$\finedim
\prop{\kft}
Let $r>0$, $i\in I^r$. Then we have 
$${o(i)^{r\over
k_r^{\prime}}(q^{-1}-q)r[d_i]_q\over
\dk_i[r]_qA_{ii}^{(r)}}(\bar
E_{(r\d,i)},\bar F_{(r\d,i)})=$$
$$=\cases
[2]_{q^{r(2n-1)}}&{\rm{if\ }}i=1\ {\rm{in\ case\
}}B_n^{(1)}\cr
[n-i+2]_{q^{2r}}[2]_{q^r}&{\rm{if\ }}i\neq 1\ {\rm{in\ cases\
}}B_n^{(1)}{\rm{and}}\ A_{2n}^{(2)}{\rm{and}}\cr
&{\rm{if\ }}i>2\ {\rm{in\ case\
}}F_{n=4}^{(1)}\cr
[2]_{q^r}[2]_{q^{r(n+1)}}&{\rm{if\ }}i=1\ {\rm{in\ case\
}}C_n^{(1)}\cr
[2]_{q^r}[2]_{q^{r(n-1)}}&{\rm{if\ }}i=1\ {\rm{in\ case\
}}D_n^{(1)}\cr
{[9-n]_{q^r}[2]_{q^{3r(n-4)}}\over[2]_{q^{r(n-4)}}}&
{\rm{if\ }}i=1\ {\rm{in\ case\
}}E_n^{(1)}\ {\rm{with\ }}n=6,7\cr
{[2]_{q^r}[2]_{q^{15r}}\over
[2]_{q^{3r}}[2]_{q^{5r}}} &{\rm{if\ }}i=1\ {\rm{in\ case\
}}E_8^{(1)}\cr 
[2]_{q^{5r}}&{\rm{if\ }}i=2\ {\rm{in\ case\
}}F_4^{(1)}\cr 
{[2]_{q^{9r}}\over[2]_{q^{3r}}}&
{\rm{if\ }}i=1\ {\rm{in\ case\
}}F_4^{(1)}\cr 
[6]_{q^r}&{\rm{if\ }}i=2\ {\rm{in\ case\
}}G_2^{(1)}\cr 
{[2]_{q^{6r}}\over[2]_{q^{2r}}}&
{\rm{if\ }}i=1\ {\rm{in\ case\
}}G_2^{(1)}\cr 
[2]_{q^{r(2n+1)}}&
{\rm{if\ }}i=1\ {\rm{and\ }}2|r\ {\rm{in\ case\
}}A_{2n}^{(2)}\cr 
[2]_{q^r}[2n+1]_{q^r}&
{\rm{if\ }}i=1\ {\rm{and\ }}2\not|r\ {\rm{in\ case\
}}A_{2n}^{(2)}\cr 
[{\tn-n\over k-1}-i+2]_{q^r}&
{\rm{if\ }}k\not|r\ {\rm{in\ cases\
}}D_{\tn}^{(k)}\ {\rm{and\ }}E_{\tn=6}^{k=2}\cr 
{(-1)^{{r\over 2}}[2]_{q^{rn}}\over 2}&
{\rm{if\ }}i=1\ {\rm{and\ }}2|r\ {\rm{in\ case\
}}A_{2n-1}^{(2)}\cr 
[2]_{q^{rn}}&
{\rm{if\ }}i=1\ {\rm{and\ }}k=2|r\ {\rm{in\ case\
}}D_{\tn}^{(k)}\cr 
{[2]_{q^{3r}}\over[2]_{q^r}}&
{\rm{if\ }}i=1\ {\rm{and\ }}k=3|r\ {\rm{in\ case\
}}D_{\tn}^{(k)}\cr 
[2]_{q^{3r}}&
{\rm{if\ }}i=2\ {\rm{and\ }}2|r\ {\rm{in\ case\
}}E_6^{(2)}\cr 
(-1)^{{r\over 2}}{[2]_{q^{6r}}\over[2]_{q^{2r}}}&
{\rm{if\ }}i=1\ {\rm{and\ }}2|r\ {\rm{in\ case\
}}E_6^{(2)}\cr 
[n-i+2]_{q^r}&{\rm{otherwise.}}\endcases$$
\dim
For the non twisted case see \rmcen.
Then also the cases $A_{2n}^{(2)}$ with $i\neq 1$, $X_{\tn}^{(k)}$ with
$k\not|r$, $A_{2n-1}^{(2)}$ and $D_{\tn}^{(k)}$ with $k\neq 1|r$ and
$i\neq 1$, and $E_6^{(2)}$ with $2|r$ and $i\neq 1,2$ are clear.
\finedim
\vskip .2 truecm
{\bf{\S \mfrm. The expression of $R$ in terms of the $\bar E_{\a}$'s.}}
\vskip .1 truecm
We can now collect the results obtained till now and get an explicit
formula for $R$. 
\teo{\emfrm}
$$R=\bigg(\prod_{\a\in\tPhi_+}{\rm{exp}}_{\a}\Big(
\rapp{q_{\a}^{-1}-q_{\a}}{c_{\a}}
\bar E_{\a}\otimes\bar F_{\a}
\Big)\bigg)
q^{-t_{\infty}}$$
where $c_{\a}=(q_{\a}^{-1}-q_{\a})(\bar E_{\a},\bar F_{\a})=\cases
1&{\rm{if}}\ \a\in\Phi_+^{{\rm{re}}}\cr 
(q_{\a}^{-1}-q_{\a})(\bar E_{\a},\bar F_{\a})&{\rm{if}}\
\a\in\tPhi_+^{{\rm{im}}},\endcases$
and $c_{\a}$ is explicitly determined from proposition \kft.
\dim
The claim is a straightforward consequence of the preceding
results: corollary \rmaf\ and proposition \kft. 
\finedim

\vskip .5 truecm

{\bf{\ref. REFERENCES.}}

\vskip .3 truecm 
\beck\ Beck, J., {\it{Braid group action and quantum affine algebras}}, 
Commun. Math. Phys. {\bf 165} (1994), 555-568. 
\medskip 

\bebr\ Beck, J., {\it{Convex bases of PBW type for quantum affine
algebras}},  Commun. Math. Phys. {\bf 165} (1994), 193-199. 
\medskip 

\bourb\ Bourbaki, N., {\it{Groupes et alg\`ebres de Lie 4, 5, 6}}, 
Hermann, Paris, 1968. 
\medskip 

\cp\ Chari, V., Pressley, A., {\it{Twisted quantum affine algebras}}, 
Commun. Math. Phys. {\bf 196} (1998), 461-476. 
\medskip 

\dasl2h\ Damiani, I., {\it{A basis of type Poincar\'e-Birkhoff-Witt for the
quantum algebra of
$\hat{\sl}_2$}}, J. Algebra {\bf 161} (1993), no. 2, 291-310. 
\medskip

\rmens\ Damiani, I., {\it{La $R$-matrice pour les alg\`ebres quantiques
de type affine non tordu}}, Ann. Scient. \'Ec. Norm. Sup., $4^e$ s\'erie,
t. 31 (1998), 493-523.
\medskip

\rmcen\ Damiani, I., {\it{The 
Highest Coefficient of det$H_{\eta}$ and the Center of the Specialization 
at Odd Roots of Unity for Untwisted Affine Quantum Algebras}}, 
J. Algebra {\bf 186} (1996), 736-780. 
\medskip 

\dgz\ Delius, G.W., Gould, M.D., Zhang, Y.-Z., {\it{Twisted quantum affine
algebras and solutions to the Yang-Baxter equation}}, Int.J.Mod.Phys. A11
(1996) 3415-3438.
\medskip 

\drnry\ Drinfeld, V.G., {\it{A new realization of Yangians and quantized 
affine algebras}}, Soviet Math. Dokl. {\bf 36} (1988), no. 2, 212-216. 
\medskip 

\drqybe\ Drinfeld, V. G., {\it{Hopf algebras and quantum Yang-Baxter
equation}}, Soviet Math. Dokl., vol. 32, 1985, p. 254-258.
\medskip 

\drqg\ Drinfeld, V.G., {\it{Quantum groups}}, Proc. ICM, Berkeley (1986),
798-820. 
\medskip 

\gmw\ Gandenberger, G.M., MacKay, N.J., Watts, G.M.T.\!, {\it{Twisted
algebra $R$\!-matrices and $S$-matrices for $B_n^{(1)}$ affine Toda
solitons and their bound states}}, Nucl.Phys. B465 (1996) 329-349.
\medskip 

\hmph\ Humphreys, J.E., {\it{Introduction to Lie Algebras and
Representation  Theory}}, Sprin\-ger-Verlag, USA (1972). 
\medskip 

\mapi\ Iwahori, N., Matsumoto, H., {\it{On some Bruhat decomposition and
the structure of the Hecke rings of ${\goth p}$-adic Chevalley groups}}, 
Inst. Hautes ƒtudes Sci. Publ. Math. No. 25 (1965) 5-48. 
\medskip 

\jmb\ Jimbo, M., {\it{A q-difference analogue of} ${\Cal U(\g)}$ { and the 
Yang-Baxter equation}}, Lett. Math. Phys. {\bf 10} (1985), 63-69. 
\medskip 

\jin\ Jing, N., {\it{On Drinfeld realization of quantum affine algebras}}, 
in Monster and Lie Algebras, ed. J. Ferrar and K. Harada, OSU Math Res
Inst Publ. 7, de Gruyter, Berlin (1998),195-206.
\medskip 

\jinginv\ Jing, N., {\it{Twisted vertex representations of quantum affine
algebras}}, Invent. Math. 102 (1990), 663-690. 
\medskip 

\KM\ Kac, V.G., {\it{Infinite Dimensional Lie Algebras}}, 
Birkh\"auser Boston, Inc., 
USA (1983). 
\medskip 

\ktrua\ Khoroshkin, S.M., Tolstoy, V.N., {\it{The universal $R$-matrix for 
quantum nontwisted affine Lie algebras}}, Funkz. Analiz i ego pril. 
{\bf 26:1} (1992), 85-88. 
\medskip 

\khtsa\ Khoroshkin, S.M., Tolstoy, V.N., {\it{Universal $R$-matrix for 
quantized(super) algebras}}, Commun. Math. Phys. {\bf 141} (1991), 
599-617. 
\medskip 

\kirresh\ Kirillov, A.N., Reshetikhin, N., {\it{$q$-Weyl Group and a 
Multiplicative Formula for Universal $R$-Matrices}}, 
Commun. Math. Phys. {\bf 134} (1990), 421-431. 
\medskip 

\lsrm\ Levendorskii, S.Z., Soibelman, Ya.S., 
{\it{Quantum 
Weyl group and multiplicative formula for the $R$-matrix of a 
simple Lie algebra}}, Funct. Analysis and its Appl., 2 {\bf 25} 
(1991), 143-145. 
\medskip                        

\lswg\ Levendorskii, S.Z., Soibelman, Ya.S., 
{\it{Some applications of quantum Weyl group I}}, J. Geom. Phys. {\bf 7} 
(1990), 241-254. 
\medskip 

\lss\ Levendorskii, S.Z., Soibelman, Ya.S., Stukopin, V., {\it{The Quantum 
Weyl Group and the Universal Quantum $R$-Matrix for Affine Lie 
Algebra $A_1^{(1)}$}}, Lett. Math. Phys. {\bf 27} (1993), 253-264. 
\medskip 

\ld\ Lusztig, G., {\it{Affine Hecke algebras and their graded version}},
J. Amer. Math. Soc. {\bf 2} (1989), 599-625. 
\medskip

\lusz\ Lusztig, G., {\it{Finite-dimensional Hopf algebras arising from
quantized enveloping algebras}}, J. Amer. Math. Soc. {\bf 3} (1990),
257-296. 
\medskip 

\liq\ Lusztig, G., {\it{Introduction to quantum groups}}, Birkh\"auser
Boston,  USA (1993). 
\medskip 

\lqgr\ Lusztig, G., {\it{Quantum groups at roots of 1}}, Geom. Ded. 
{\bf 35} (1990), 89-113. 
\medskip 

\matsu\ Matsumoto, H., {\it{G\'en\'erateurs et r\'elations des groupes de
Weyl  g\'en\'eralis\'es}}, C. R. Acad. Sci. Paris {\bf T. 258 II} (1964), 
3419-3422. 
\medskip 

\rRm\ Rosso, M., {\it{ An Analogue of P.B.W. Theorem and the Universal 
R-Matrix for $U_hsl(N+1)$}}, Commun. Math. 
Phys. {\bf 124} (1989), 307. 
\medskip 

\tan\ Tanisaki, T., {\it{Killing forms, Harish-Chandra isomorphisms
and universal R-matrices for Quantum Algebras}}, in Infinite Analysis
Part B, Adv. Series in Math. Phys., vol 16, 1992, p. 941-962.
\medskip 
\bye